\documentclass[reqno]{amsproc}

\usepackage{psfig,latexsym,graphicx, epsfig}
\usepackage{amssymb,amsmath,amsthm,amsopn,pifont} 
\usepackage{subfigure}
\usepackage{multirow}
\usepackage{psfrag} 
\newcommand\Aut{ {\rm Aut}}

\newcommand\cl{{\mathcal C}}
\newcommand\hl{{\mathcal H}}
\newcommand\p{{\mathcal P}}
\newcommand\s{{\mathcal S}}
\newcommand\C{{\mathbb C}}
\newcommand\D{{\mathbb D}}
\newcommand\Z{{\mathbb Z}}
\newcommand\R{{\mathbb R}}
\newcommand\bbeta{{\boldsymbol\beta}}
\newcommand\f{{\bf f}}
\newcommand\w{{\bf w}}

\newcommand\B{{\mathcal B}}
\newcommand\U{{\mathcal U}}
\newcommand\aac{{a}}
\newcommand\bb{{b}}
\newcommand\bio{{\boldsymbol\iota}}
\newcommand\bsig{{\boldsymbol\sigma}}
\newcommand\cm{{\rm cm}}
\newcommand\ssm{{\smallsetminus}}
\newcommand\QED {$\quad\square$}
\newcommand\mixarr{\;\hbox{\raise 2pt \hbox{$\leftarrow$}}
        \hbox{\lower 1pt \hbox{\hskip -9pt $\to$}}\;}

\newcommand\G{{\mathcal G}}
\newcommand\biota{{\boldsymbol\iota}}  

\newtheorem  {theorem}{Theorem}[section]
\newtheorem  {lem}[theorem]{Lemma}

\newtheorem{coro}[theorem]{Corollary}

\theoremstyle{definition}
\newtheorem   {definition}[theorem]{Definition}
\newtheorem   {ex}[theorem]{Example}

\theoremstyle{remark}
\newtheorem  {rem}[theorem]{Remark}

\numberwithin{equation}{section}

\input xy
\xyoption{all}

\begin{document}

\title[Hyperbolic Components]{Hyperbolic Components}

\thanks{I want to thank the NSF for its support under grant DMS0757856.}
\author[John Milnor]{John Milnor\\ 
 With an Appendix by A. Poirier.} 
\address{Institute for Mathematical Sciences, Stony Brook University, 
Stony Brook, NY 11794-3660}
\email{jack@math.sunysb.edu}




\date{}

\begin{abstract}
Consider polynomial maps $f:\C\to\C$ of degree $d\ge 2$, 
or more generally polynomial maps from a finite union of copies of 
$\C$ to itself.\break  
In the space of suitably normalized maps of this type, the 
hyperbolic maps 
form an open set called the \textbf{\textit{~hyperbolic
locus}}. The various connected components of this hyperbolic locus are
called \textbf{\textit{~hyperbolic components}},
 and those hyperbolic components with compact closure (or equivalently
those contained in the ``connectedness locus'')
are called \textbf{\textit{~bounded~}} hyperbolic components.
 It is shown that each bounded hyperbolic
component is a topological cell containing a unique post-critically  
finite map called its \textbf{\textit{~center point}}. For each degree
$d$, the bounded hyperbolic  
components can be separated into finitely many distinct types,  
each of which is characterized by a suitable \textbf{\textit{~reduced mapping 
scheme}} $\overline S_f$.  
 Any two components with the same reduced mapping scheme are 
canonically biholomorphic to each other. There are similar statements
for real polynomial maps, for polynomial
 maps with marked critical points, and for rational maps. Appendix
A, by Alfredo Poirier, proves that every reduced mapping scheme
can be represented by some classical hyperbolic component, made up of
 polynomial maps of $\C$.
This paper is a revised version
 of \cite{M2}, which was circulated but not published in 1992. 
\end{abstract}

\maketitle

\def\IMSmarkvadjust{0 pt}
\def\IMSmarkhadjust{0 pt}
\def\IMSmarkhpadding{0 pt}
\def\IMSpubltext{Published in modified form:}
\def\SBIMSMark#1#2#3{
 \font\SBF=cmss10 at 10 true pt
 \font\SBI=cmssi10 at 10 true pt
 \setbox0=\hbox{\SBF \hbox to \IMSmarkhpadding{\relax}
                Stony Brook IMS Preprint \##1}
 \setbox2=\hbox to \wd0{\hfil \SBI #2}
 \setbox4=\hbox to \wd0{\hfil \SBI #3}
 \setbox6=\hbox to \wd0{\hss
             \vbox{\hsize=\wd0 \parskip=0pt \baselineskip=10 true pt
                   \copy0 \break%
                   \copy2 \break%
                   \copy4 \break}}
 \dimen0=\ht6   \advance\dimen0 by \vsize \advance\dimen0 by 8 true pt
                \advance\dimen0 by -\pagetotal
	        \advance\dimen0 by \IMSmarkvadjust
 \dimen2=\hsize \advance\dimen2 by .25 true in
	        \advance\dimen2 by \IMSmarkhadjust

%
%
  \openin2=publishd.tex
  \ifeof2\setbox0=\hbox to 0pt{}
  \else 
     \setbox0=\hbox to 3.1 true in{
                \vbox to \ht6{\hsize=3 true in \parskip=0pt  \noindent  
                {\SBI \IMSpubltext}\hfil\break
                \input publishd.tex 
                \vfill}}
  \fi
  \closein2
  \ht0=0pt \dp0=0pt
 \ht6=0pt \dp6=0pt
 \setbox8=\vbox to \dimen0{\vfill \hbox to \dimen2{\copy0 \hss \copy6}}
 \ht8=0pt \dp8=0pt \wd8=0pt
 \copy8
 \message{*** Stony Brook IMS Preprint #1, #2. #3 ***}
}

\SBIMSMark{2012/2}{March 2012}{}


\section{Introduction.}\label{s1}\smallskip

\begin{definition}\label{d-ms}
 A \textbf{\textit{~hyperbolic mapping scheme~}} $S$~
(or briefly a \textbf{\textit{~scheme\,}})
consists of a finite set $|S|$ of ``vertices'',
 together with a map $F=F_S:|S|\to|S|$, and
an integer valued \textbf{\textit{~critical weight function}}~
$s\mapsto \w(s)\ge 0$, satisfying two conditions:
\begin{quote}
$\bullet$ Any vertex of weight zero is the iterated
 forward image of some vertex of positive weight, $\w(s)\ge 1$.

\noindent$\bullet$ (Hyperbolicity.) Every periodic orbit
 under $F$  contains at least one vertex of positive weight.
\end{quote}
\noindent The number $d(s)=\w(s)+1\ge 1$ will be called the
 \textbf{\textit{\,degree\,}} of the vertex $s$.
The scheme is called \textbf{\textit{~reduced~}} if $\w(s)\ge 1$ (or
$d(s)\ge 2$) for every $s\in|S|$.
\end{definition}\medskip

First consider a polynomial
 map $f:\C\to\C$ of degree $d\ge 2$ with connected Julia set which
is \textbf{\textit{\,hyperbolic\,}} in the sense that every critical
orbit converges to an attracting cycle.

\begin{definition}\label{d-S_f}
The \textbf{\textit{\,full mapping scheme\,}} $S_f$ of such a map
 has one vertex $s_U$ corresponding to each component $U$ of the Fatou set
which contains a critical or post-critical point. The weight
$\w(s_U)\ge 0$ is defined to be the number of critical points in $U$,
counted with multiplicity, and the associated map $F_f:|S_f|\to|S_f|$
carries $s_U$ to $s_{f(U)}$.
\end{definition}

However for many purposes a slightly simpler structure is more useful.
Every mapping scheme can be simplified to an associated reduced scheme
(see Remark \ref{r-ars}). In particular:

\begin{definition}\label{d-red}
The \textbf{\textit{~ reduced mapping scheme~}}  $\overline S=\overline S_f~$
associated with a hyperbolic polynomial map $f$ can be described as follows:
\begin{itemize}
\item There is one vertex $s=s_U\in|\overline S|$ for
 each Fatou component \hbox{$U\subset K(f)$} which contains at least one 
critical point.
\smallskip

\item The weight $\w(s)$ is again the
number of critical points in $U$, counted with multiplicity.
\smallskip

\item The map~ $F:|\overline S|\to|\overline S|$
is defined by $F(s_U)=s_{U'}$, where\break
\hbox{$U'=f^{\circ n}(U)\,,~ n>0$}, is the first forward image which contains
 a critical point.
\smallskip
\end{itemize}\end{definition}

Poirier has shown that every reduced mapping scheme can be obtained in
this way, from some hyperbolic map from $\C$ to itself. (See Appendix A.)

\medskip

\subsection*{Outline of what follows.} Section \ref{s-PS}
will introduce the space $\p^{S_0}$ of suitably normalized polynomial maps,
associated with any mapping scheme $S_0$, and modify several classical
 definitions so that they apply in this more general context.
Section \ref{s2} will provide a graphical description of mapping
 schemes, and discuss symmetries. Sections \ref{s3} and \ref{s4}
 will provide a universal
 topological model, based on Blaschke products, for hyperbolic components with
a specified reduced mapping scheme, showing that each hyperbolic component
is a topological  cell with a preferred center point.
 Section \ref{s5} will sharpen this result
by providing a universal biholomorphic model.
Section \ref{s6} discusses analogous results for polynomial mappings with
real coefficients, and
 more generally for \textbf{\textit{~real forms~}} of complex polynomial
mappings. Section \ref{s7} studies polynomial mappings which have been
{\textbf{\textit{critically marked\/}}} by specifying an ordered list of their
 critical points. It is shown that all of the principal results carry over to
 the critically marked case. Section \ref{s-rat} proves analogous results
for rational maps.

 Appendix \ref{sa}, by Alfredo Poirier, shows that 
every reduced scheme actually occurs as the scheme $\overline S_f$ for
some critically finite hyperbolic map $f:\C\to\C$. Appendix \ref{sab}
studies the number of distinct reduced schemes with given total weight.

The present work is a fairly straightforward extension of ideas originated
by Douady, Hubbard, McMullen, Rees and others, and many of
the statements were probably known as folk theorems. I am
 particularly grateful to Branner and Douady for their considerable help
with the earlier version, and to Araceli Bonifant, Adam Epstein,
Alfredo Poirier, and Scott Sutherland for their help with the present version.

\section{The Affine Parameter
Spaces $\p^d$ and $\p^{S_0}$.}\label{s-PS}

First consider the classical case of polynomial maps $f:\C\to\C$.
(See for example \cite{D1}, \cite{DH1},  or \cite{M4}.)

\begin{definition}\label{d-mcp}
A complex polynomial map $$\,f(z)=\sum_{j=0}^da_jz^j\,$$ will be called 
 \textbf{\textit{monic\/}} and \textbf{\textit{centered\/}}\, if
$\,a_d=1\,$ and $\,a_{d-1}=0$. (In the degree one case, by definition,
only the identity map is monic and centered.)  For $\,d\ge 2$,~
let  $\p^{d}$ be the complex $(d-1)$-dimensional affine
space consisting of all polynomial maps $f:\C\to\C$ which are
monic and centered. For each such $\,f$,~
the \textbf{\textit{filled Julia set~}} $K(f)\subset\C$ is the union of
all bounded orbits, and the
 \textbf{\textit{connectedness locus~}} $\cl^d\subset\p^{d}$
 is the
compact set consisting of all polynomials $f \in \p^{d}$ for which the filled
Julia set is connected, or equivalently contains all critical points.
 A polynomial or rational map
 is \textbf{\textit{~hyperbolic~}} if the orbit of every
critical point converges to an attracting cycle. (Here convergence
to the attracting fixed point at infinity is allowed, although we will not
 be interested in that case.)
 The open set consisting
of all hyperbolic maps will be denoted by $\hl^d\subset\p^d$.
\smallskip

A connected component $H\subset\hl^d$ 
has compact closure if and only if it is contained in $\cl^d$, or if and
only if every critical orbit converges to a finite attracting cycle.
Those connected components of $\hl^d$ which are contained in $\cl^d$
will be called
\textbf{\textit{~bounded hyperbolic components}}.
It is not hard to see that all of the maps $f$ in such a bounded
hyperbolic component $H$ have isomorphic\footnote
{We will see in \S\ref{s4} that each such
 $H$ is simply-connected, so that these isomorphisms are uniquely defined.}
reduced mapping schemes $\overline S_f$, so we can use the alternate notation
$\overline S_H=\overline S_f$. Note that the  \textbf{\textit{~total weight~}}
\begin{equation}\label{e-tw}
 \w(\overline S_H)~=~\sum_{s\in\, |\overline S_H|} \w(s)
\end{equation}
associated with each $H\subset\hl^d$ 
is equal to the complex dimension $\,d-1\,$ of $\hl^d$.
The number of isomorphism classes of reduced schemes grows rapidly with
 the total weight $\w(S)$.
(See Table \ref{t-numbers} for small values of $\w(S)$, and see
Figure \ref{f-wt2} for the special case $\w(S)=2$.) For details, see
 Appendix \ref{sab}. 
\end{definition}

According to Poirier, every one of these reduced  schemes can be 
realized by a suitable hyperbolic component in ${\mathcal P}^{\w(S)+1}$.
\smallskip

\begin{table}[!ht]
\begin{center}\begin{tabular}{|lcccccc|}
\hline
$\w(S)$  & 1 & 2 & 3 & 4&  5 &6\cr 
number&   1&  4 &12 &42& 138& 494\cr
\hline
\end{tabular}
\vspace{.3cm}

\caption{\label{t-numbers} \it The numbers of distinct reduced schemes with 
$\w(S)\leq 6$.}
\end{center}
\end{table}

In order to obtain a canonical model for hyperbolic components with a specified
reduced mapping scheme, we need to extend the concept of polynomial map
by allowing maps from some disjoint union of finitely many copies of $\C$ to
itself. More explicitly, we will consider the following. Let $S_0$ be an
 arbitrary  mapping scheme. We will think of the product
 $|S_0|\times\C$ as a disjoint union of copies of $\C$, indexed by the
points $s\in|S_0|$.

\begin{definition}[{\sc The parameter space $\p^{S_0}$}]\label{d-p^S}
By a \textbf{\textit{\,generalized polynomial map\,}}
based on the scheme $S_0$ will be meant a map
$${\bf f}: |S_0|\times\C\to|S_0|\times\C $$
which sends each $~s\times\C~$ onto $~F(s)\times\C~$ by a polynomial map
 of degree $d(s)$, where $F=F_{S_0}$.
Such a map is {\textbf{\textit\,normalized\,}} if each of these polynomial
maps $~s\times\C\to F(s)\times\C~$ is monic
and centered. (Compare Remark \ref{r-normalize}.)
The complex affine space consisting of all such normalized maps
will be denoted by $\p^{S_0}$. Thus $\f\in\p^{S_0}$ if and only if $\f$
has the form
$${\bf f}(s,z)~=~\big(F(s),\,f_s(z)\big)\,,$$
where each $f_s:\C\to\C$ is a monic centered
 polynomial of degree $d(s)=\w(s)+1$.
There is a preferred base point $\f_0\in\p^{S_0}$ given by
\begin{equation}\label{e-f0}
 \f_0(s,\,z)~=~\big(F(s),~z^{d(s)}\big)\,.
\end{equation}
In the special case where $|S_0|$ consists of a single point of weight $\w=d-1$,
note that $\p^{S_0}$ can be identified with the space $\p^d$ of Definition
 \ref{d-mcp}. Many of the basic definitions and
 results in the case of a map $f:\C\to\C$ carry over easily to
this more general context. (Proofs will be omitted if they are completely
analogous to the proofs in the classical case, as given for example in
\cite{DH1} or \cite{M4}.)
\end{definition}

\begin{definition}\label{d-J}
First consider the ``dynamic space'' $|S_0|\times\C$.
The \textbf{\textit{~Fatou set~}} associated with any map $ {\bf f}
\in\p^{S_0}$ is defined to be the open subset of $|S_0|\times\C$
 consisting of all points    $(s,z)$
such that the iterates of $ {\bf f}$, restricted to some neighborhood
 of $(s,z)$, form a normal family. Each connected component of the
Fatou set is called a \textbf{\textit{~Fatou component}}. The map 
${\bf f}$ is \textbf{\textit{hyperbolic}} if every critical orbit converges to
 a periodic orbit.

There are two reasonable concepts of the ``Julia set'' in this context.
The complement of the Fatou set in $|S_0|\times\C$ will be called the
\textbf{\textit{~fully invariant Julia set~}} $J( {\bf f})$.
Alternatively, following Julia, one could consider
the closure of the set of repelling periodic orbits. This forms a compact
forward invariant set $J_{\rm rec}( {\bf f})$,
which can be called
the \textbf{\textit{~recurrent Julia set}}. Note that $J_{\rm rec}(\f)$ is
strictly smaller than $J( {\bf f})$ whenever the map $F:|S_0|\to |S_0|$ is not
surjective.\smallskip

The union of all orbits which are bounded (i.e., contained in a compact subset
of $|S_0|\times\C$) is a compact set $~K({\bf f})\subset|S_0|\times\C$
 called the  \textbf{\textit{~filled Julia set}}.
The boundary $\partial K(\f)$ is equal to $J(\f)$; and $K(\f)$
 can be described as the union of $J(\f)$ with all
bounded Fatou components. Just as in the classical case, every bounded
Fatou component is biholomorphic to the open unit disk; and if
 $\f$ is hyperbolic, then the boundary of each such component
 is a Jordan curve. (Compare \S\ref{s4}.)
\end{definition}

\begin{definition}\label{d-par-sp} Now consider the \textbf{\textit{~parameter
space~}} $\p^{S_0}$.
Evidently $\p^{S_0}$ is a complex affine space with complex
 dimension equal 
 to the total weight $$\w(S_0)~=~\sum_{s\in \vert S_0\vert}\, \w(s) $$
(or to the total number of critical points, counted with multiplicity).
The \textbf{\textit{~connectedness locus~}} is
 defined to be the compact set
$\cl^{S_0}\subset\p^{S_0}$ consisting of all maps $ {\bf f}\in\p^{S_0}$ for which all critical
points are contained in $K( {\bf f})$. Equivalently, 
$\cl^{S_0}$ can be described as the set of all $\f\in\p^{S_0}$ such that the
 intersection of $K(\f)$ with each  $s\times\C$ is connected.
 The notation $\hl^{S_0}$ will be used for the set of all hyperbolic maps
 in $\p^{S_0}$. Each connected component of $\hl^{S_0}\cap\cl^{S_0}$ will
 be called a \textbf{\textit{~bounded hyperbolic component~}} $H$.
 Just as in Definition \ref{d-S_f},
we can define the mapping scheme $S_H$ associated with
each bounded hyperbolic component $H$.
(In most cases this new scheme $S_H$ will {\bf not} be
 the same as the ambient scheme 
$S_0$, although there is a natural map
from $S_H$ onto $S_0$.)
\end{definition}
\medskip

\begin{rem}
 Both the statement that each hyperbolic component is a topological
cell, and the statement that it has a preferred center point,
are strongly dependent on the fact that we consider only hyperbolic components
within the connectedness locus---the structure of hyperbolic components
outside the connectedness locus is very different. For instance,
Blanchard, Devaney and Keen \cite{BDK}
show that the \textbf{\textit{~shift locus\,}},~ the
unbounded hyperbolic component consisting
of maps for which all critical orbits escape to infinity,
has a very complicated fundamental group when $d\ge 3$. (In the somewhat
analogous moduli space for quadratic rational maps with marked critical points,
 there is a similar ``shift locus hyperbolic component'' which
 contains a Klein bottle as retract,
and hence also has a non-abelian fundamental group. Compare \cite[\S8.7]{M3}.)
\end{rem}\smallskip

\section{Graphs and Symmetries.}
\label{s2}

It is often convenient to represent each scheme $S$ by a finite graph
$\Gamma(S)$, with
the points of $|S|$ as vertices, and with a directed edge leading from each
vertex $s$ to $F(s)$. By definition, the \textbf{\textit~degree~} of such
an edge is  equal to $~d(s)=\w(s)+1$.
\smallskip

In the figures, each vertex of positive critical weight
$\w>0$ is represented by a cluster of $\w$ heavy dots, while vertices of weight
zero (if any) are represented by much smaller dots.

\begin{figure}[!ht]
\centerline{\psfig{figure=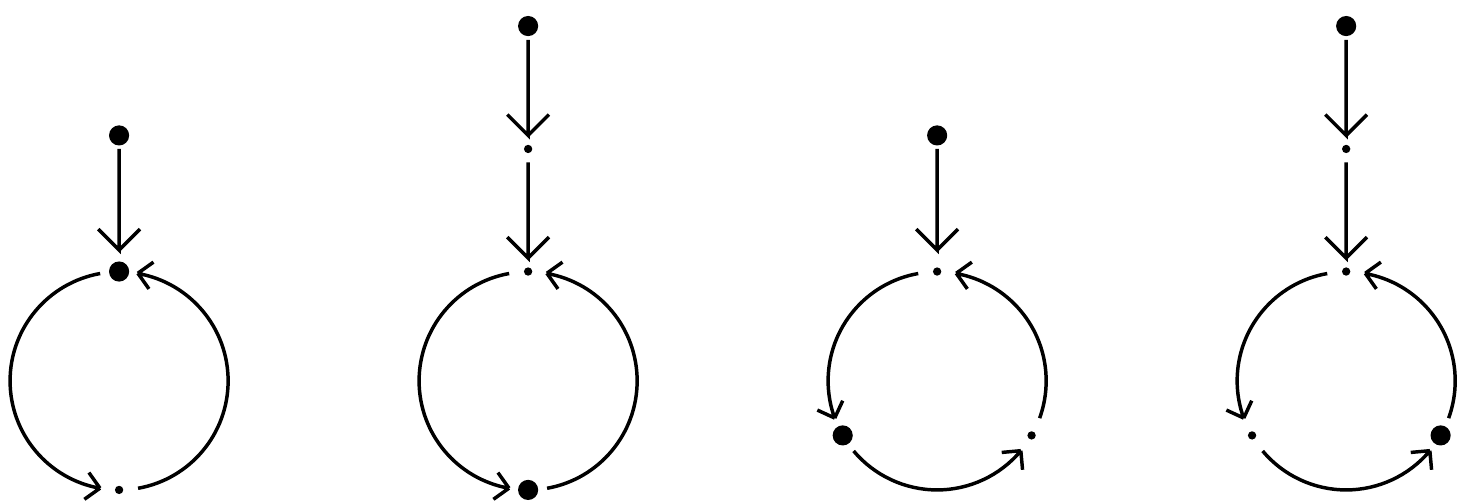,height=1.3in}}
\caption{\it \label{f-wt2p} Four different full schemes which give
rise to the same reduced scheme.}
\end{figure}

\begin{rem}[{\sc The Associated Reduced Scheme}]\label{r-ars}
Every mapping scheme $S$ 
 gives rise to an 
 \textbf{\textit{~associated reduced scheme~}} $\overline S$, as follows.
 By definition,
$|\overline S|$ is the subset of $|S|$ consisting of points of positive
 weight; and the critical
weight of each vertex of $\overline S$ is the same as its critical weight
 in $S$. The associated function $F_{\overline S}$ from $|\overline S|$
to itself is obtained by iterating $F_S:|S|\to|S|$ 
until we reach a vertex $s'$ of positive weight.
If we start with the graph $\Gamma(S)$ of an arbitrary mapping scheme,
then the graph $\Gamma({\overline S})$ of the associated
 reduced scheme can be obtained from $\Gamma(S)$ simply by shrinking
 each edge of degree one joining $s$ to $F(s)$ to its endpoint $F(s)$.
If we start with the full mapping scheme $S_f$ of a hyperbolic map,
 then clearly the associated reduced scheme $\overline S_f$ constructed
 in this way is identical to the object described in
 Definition \ref{d-red} above.
Note that for any $S_0$ the affine space $\p^{S_0}$
can be identified with $\p^{\overline S_0}$.
\end{rem}
\smallskip
\begin{figure}[!ht]
\centerline{\psfig{figure=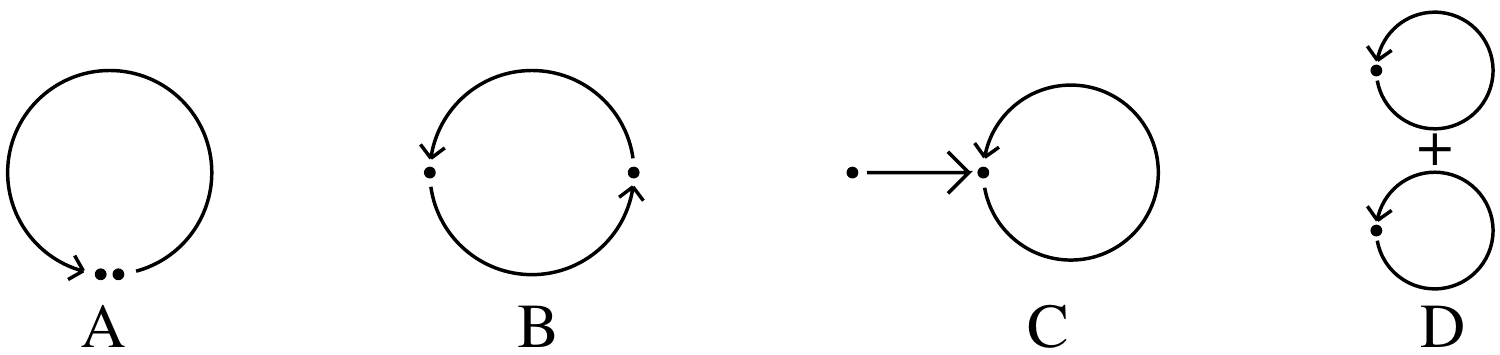,height=0.95in}}
\caption{\it \label{f-wt2} Graphs for  the four reduced schemes of weight 
\hbox{$\w(S)=2$}. (Compare the discussion of cubic polynomials of 
Example~\ref{ex-w=2}.)}
\end{figure}

The full scheme provides quite a bit of information
about any given hyperbolic map which is lost in this reduced
scheme. In fact there are infinitely many possible
 full schemes for each reduced scheme.
Figure \ref{f-wt2p} gives four different examples of schemes
of total weight two which are associated with cubic maps of $~\C$.
All of these correspond to the same reduced scheme, which is shown in
 Figure \ref{f-wt2}C.

Figure \ref{f-wt2} shows the graphs of all distinct
reduced schemes with total weight $\w(S)=2$. (Compare \cite{M1}.)
The numbers of distinct reduced schemes of given total weight are shown in
 Table \ref{t-numbers} of \S\ref{s-PS} for $\w(S)\le 6$.
(Compare Appendix \ref{sab}.) However, I don't know any
formula for the number of such schemes in general.\medskip

\subsection*{Symmetries.} 
Let $\G=\G(S)$ be the finite abelian group consisting of all maps
 $~ {\bf g}:|S|\times\C\to|S|\times\C~$
which send each $~s\times\C~$ linearly onto itself, and which commute with the
base map 
\begin{equation*}\label{e-f00}
  {\bf f}_0(s,z)=\big(F(s),\,z^{d(s)}\big)\,.
\end{equation*}
More explicitly, each $~ {\bf g}\in\G~$ must have the form
\begin{equation}\label{e-rho} {\bf g}(s,z)~=~ (s,\,\rho_sz)\,,
\end{equation}
where each vertex $s$ is assigned a root of unity $\rho_s$ 
 satisfying the condition 
\begin{equation}\label{e-sym}
 \rho_s^{d(s)}~=~ \rho_{F(s)}\,.
\end{equation}

\begin{rem}\label{r-normalize}
 The restriction to monic centered maps in Definition \ref{d-p^S}
can be justified as follows.
If we start with an arbitrary generalized polynomial map with scheme $S_0$
which is not required to be monic or centered,
then it is not difficult to find an automorphism 
${\bf h}$ of $|S_0|\times\C$ which carries each
 component holomorphically onto itself, so that 
${\bf h}^{-1}\circ\f\circ {\bf h}\in\p^{S_0}$. This $\bf h$ is uniquely
determined up to composition with some ${\bf g}\in\G(S_0)$.
\end{rem}

The order of this abelian group $\G(S)$ can be computed as follows. It suffices
to consider the connected case, since the automorphism group
$~\G(S+S')~$ of a disjoint union 
is clearly isomorphic to $~\G(S)\times\G(S')$.

\begin{lem}\label{l-|G|}
The order of the symmetry group $\G(S)$ of a connected mapping scheme
 is equal to the product
$$ \big(d_1\cdots d_k-1\big)\, d_{k+1}\cdots d_\ell\,,$$
where   $~d_1,\,\ldots,\,d_k~$
are the degrees of the vertices which belong to its unique cycle,
 and $~d_{k+1},\,\ldots,\,d_\ell~$ are
 the degrees of the remaining aperiodic vertices.
\end{lem}

\begin{proof}  First consider the case where all vertices are periodic,
so that $S$ consists only
of the cycle. Then we must have
 $~\rho_j^{d_j}=\rho_{j+1}\,$,~ where $j$
is understood to be an integer modulo $k$. It follows that
$$  \rho_j^{d_1d_2\cdots d_k}~=~\rho_j\,.$$
Thus $\rho_1$ can be an arbitrary $~(d_1d_2\cdots d_k-1)$-th root of unity,
and the remaining $\rho_j$ are then uniquely determined. Thus,
 in this case, $\G$ is cyclic of order $~d_1\cdots d_k-1$.

The proof now continues inductively, building up the scheme $S$
by adding one new vertex at a time
outside of the cycle. Evidently there are exactly $d_j$ possible
 choices for each new $\rho_j$, and the conclusion follows.
\end{proof}

In particular, vertices of degree one make no contribution
to the order of $\G(S)$; and in fact
it is easy to check that the symmetry group $\G(\overline S)$ for the associated
reduced scheme is isomorphic to $\G(S)$.

\begin{rem}
Each $~{\bf g}\in\G(S)~$ acts linearly on the affine space $\p^{S}$,~ sending 
each map $~{\bf f}:|S|\times\C\to|S|\times\C~$ to the map
 $${\bf f^g}~=~{\bf g}^{-1}\circ {\bf f}\circ {\bf g}\,.$$
In fact setting $~ {\bf f}(s,z)=\big(F(s),\, z^d+a_{d-2}z^{d-2}+
\cdots+a_0\big)~$ with $d=d(s)$, ~ it follows easily that
\begin{equation}\label{e-action}
{\bf f}^{\bf g}(s,\,z)~=~\Big(F(s),\, \frac{(\rho_sz)^d\,+\,a_{d-2}(
\rho_sz)^{d-2}\,+\, \cdots\,+\,a_1\rho_s z\,+\,a_0\,}{\rho_{F(s)}}\Big)
\end{equation}
is again an element of $\p^{S}.$
However, this action of $~\G(S)~$ on $\p^{S}$ is not always faithful.

\begin{definition} Let $\G_0(S)$ be the subgroup of $\G(S)$ consisting
of all ${\bf g}\in\G(S)$ which commute with every ${\bf f}\in\p^{S}$,
so that the action of $\bf g$ on $\p^{S}$ is trivial.
Thus the quotient group $\G(S)/\G_0(S)$ acts faithfully (i.e.,
 effectively) on $\p^{S}$.
\end{definition}

It will be convenient to define the ``free'' vertices of $S$ to be
 those which do not belong to the image $F(|S|)\subset|S|$.
It follows from Definition \ref{d-ms} that every free vertex has degree
$d(s)\ge 2$.

\begin{lem}\label{l-G_0} An element  ${\bf g}\in\G(S)$ belongs to this
 subgroup $\G_0(S)$ if and only if the associated roots of unity
 $\{\rho_s\}$ of equation $(\ref{e-rho})$ satisfy
$$ \rho_s~=~\begin{cases}
  \pm 1 & {\rm if}~ s~ {\rm is~ a ~free~ vertex~ with~ degree}~
 d(s)~{\rm equal~to}~2\,,~{\rm but}\\ 
+1 & {\rm in~ all~ other~ cases}\,.
\end{cases}$$
 Thus $\G_0(S)$ can be described as a direct sum of copies of the group
$\{\pm 1\}$, with one copy for each free vertex of degree two.
\end{lem}

\begin{proof} First note that $\rho_{F(s)}=1$ whenever
 $d(s)\ge 2$:
This follows by considering the constant term $a_0/\rho_{F(s)}$ in equation
(\ref{e-action}). Since every $s'\in F(|S|)$ is in the forward
orbit of some vertex
of degree $d(s)\ge 2$, it follows inductively from equation (\ref{e-sym})
that $\rho_{s'}=1$ for all $s'\in F(|S|)$. Now let $s$ be a free vertex. If
$d(s)>2$, then by considering the coefficient $\rho_sa_1/\rho_{F(s)}$
 of the linear term in
equation (\ref{e-action}) we see that $\rho_s=1$. But in the case $d(s)=2$
this linear term does not appear, and we can conclude only that
 $\rho_s^{\;2}=1$. This proves that any element of $\G$
which acts trivially must belong to $\G_0(S)$,
 and the converse statement follows similarly.
\end{proof} 
\end{rem}

\begin{rem}\label{r-Ghat} One can also consider the group $\Aut(S)$
consisting of all one-to-one maps  $\phi:|S|\to|S|$ which commute
 with $F$ and preserve the critical weight. This group acts
faithfully on the space $|S|\times\C$, mapping each pair $(s,\,z)$ to
 $(\phi(s),\,z)$. The groups $\G(S)$ and $\Aut(S)$
together generate a split extension
$$  1~\to~\G(S)~\to~\widehat\G(S)~\leftrightarrow~\Aut(S)\to 1\,,$$
consisting of all compositions 
$$ \phi\circ{\bf g}~:~(s,~z) ~\mapsto~ \big(\phi(s)\,,~\rho_s z\big)~.$$
This group acts faithfully on $|S|\times\C$,
 and hence acts (not always faithfully) on $\p^S$.
\end{rem}

\section{Blaschke Products and the Model
 Space $\B^S$.}\label{s3}

 This section will describe
 a topological model, based on Blaschke products, for hyperbolic components
 with mapping scheme $S$. First a review of some standard facts. 
Let $\D$ be the open unit disk in $\C$. For any $a \in \D$,
there is one and only one M\"obius transformation $\mu_a$ of the Riemann
sphere which maps $\D$ onto itself satisfying
$$ \mu_a(a)\,=\,0\,,\qquad{\rm and}\quad \mu_a(1)\,=\,1\,.$$
It is given by
\begin{equation}
 \mu_a(z) \;=\; k\,\frac{z-a}{1-\overline a\, z}\qquad{\rm with}\qquad k =
\frac{1-\overline a}{1-a}\,.
\end{equation}
\smallskip

\begin{lem}\label{l-bp}
Any proper holomorphic
map from $\D$ onto itself extends continuously over the closed disk
$\overline\D$, and can be written uniquely
as an $d$-fold \textbf{\textit{\,Blaschke product}}
\begin{equation}\label{e-bp}
\beta(z)\; =\; \beta(1)\;\mu_{a_1}(z)\,\cdots\,\mu_{a_d}(z)
\end{equation}
with $d\ge 1$, where $|\beta(1)|=1$, and
where  $~a_1\,\ldots,\,a_d~$ are the $($not necessarily distinct$)$ pre-images
of zero.
\end{lem}

\begin{rem}\label{r-bp}
It follows that every such map $\beta$ extends
uniquely as a rational map of degree $d$
 from the Riemann sphere $~\widehat\C=\C \cup \infty~$
 onto itself. It is not hard to check that this extended 
map commutes with the inversion $z\mapsto 1/\overline z \,=\,z/|z|^2$
 in the unit circle.  In particular, $z$ is a critical point if and only if 
$1/\overline z$ is  critical, and $z$ is periodic if and only if 
$1/\overline z$ is periodic.
\end{rem}

\begin{proof}[Proof of Lemma \ref{l-bp}]
Since $\beta$ is a proper map, the absolute value $|\beta(z)|$ tends to one
as $|z|\to 1$; and since $\beta$ is onto, it has at least one zero
 $a_1\in\D$.
It follows that the quotient $\beta_1(z)=\beta(z)/\mu_{a_1}(z)$ is a
 well defined holomorphic function on $\D$. Furthermore
 $|\beta_1(z)|\to 1$ as
$|z|\to 1$, and it follows from the maximum principle
that $|\beta_1(z)|\le 1$ everywhere in $\D$. If this function $\beta_1(z)$
 is constant, this completes the proof. But if $\beta_1(z)$ is 
non-constant, then it follows from the minimum principle that there is a zero
$\beta_1(a_2)=0$, and we can continue inductively, setting
$\beta_2(z)=\beta_1(z)/\mu_{a_2}(z)$, and so on. This induction must stop
after finitely many steps since a proper holomorphic map can have at most
finitely many zeros, counted with multiplicity.
\end{proof}

\begin{lem}\label{l-3.1}
A proper holomorphic map $\beta$ of degree $d \ge 2$ from the unit 
disk onto itself induces an $d$-to-one covering map from the circle 
$\partial \D$ onto
 itself. Such a map $\beta$ has at most one fixed point in the open disk $\D$.
 If there is an interior fixed point, then
there are exactly $~d-1~$ distinct boundary fixed points.
On the other hand, if there is no interior fixed point, then all
$d+1$ fixed points in the Riemann sphere, counted with multiplicity,
must lie on the unit circle. In any case, there are $d-1$ critical points,
counted with multiplicity, in the interior of $\D$, and none on the boundary.
\smallskip
\end{lem}

\begin{proof} If we set $z=e^{2\pi i\theta}$ and $\beta(z)=
e^{2\pi i\phi(\theta)}$, then $d\phi/d\theta$ can be identified with
 the (modified) logarithmic derivative
$$\frac{d\log\beta(z)}{d\log(z)}\; =\; \frac{\beta'/\beta}{z'/z}\;=\;
\frac{z\,\beta'}{\beta}	$$
evaluated on the unit circle (where $\beta'=d\beta/d z$ and $z'=1$).
 For $d=1$, since the circle maps
diffeomorphically onto itself, we have $d\phi/d\theta>0$ everywhere,
and the integral around the circle is given by $\oint d\phi=+1$. In
the case $d>1$, it follows from (\ref{e-bp})
that the logarithmic derivative is the sum of $d$ such
terms, hence we again have $d\phi/d\theta>0$, but with $\oint d\phi=d$.
 Thus $\beta$ induces an
$d$-fold covering map from the unit circle onto itself. In particular,
 any equation of the form 
$\beta(z)= {\rm constant}\in \partial \D$ has
exactly $d$ distinct solutions.

Now suppose that $\beta$ has a fixed point in the
open unit disk.
Then after conjugating by a conformal automorphism, we may assume that
the fixed point is $z=0$. In that case, since $\mu_0(z)=z$, we can write
$$\beta(z)\; =\; \beta(1)\;z\;\mu_{a_2}(z)\cdots\mu_{a_{d}}(z)
\qquad{\rm with}\qquad |\beta(1)|=1 \,.$$
It follows immediately from this formula that
 $|\beta(z)|<|z|$ for all $z \ne 0$ in the open disk. Thus this fixed
point is attracting, and is unique within $\D$. Furthermore,
expressing the logarithmic derivative on the unit circle
as an $d$-fold sum, as above,
the first term is $+1$, hence $d\phi/d\theta=
z\beta'(z)/\beta(z)>1$ whenever $d>1$.

 A similar argument shows that
there are exactly $d-1$ fixed points on the unit circle (all
repelling). In fact the difference function
$\theta\mapsto\phi-\theta$ is a covering map from
 the circle $\R/\Z$ onto itself with
degree $d-1$, and the zeros of this difference in $\R/\Z$
are exactly the fixed points.

Finally, there can be no critical points on the boundary. In fact if the
derivative vanishes at a boundary point $z_0$, then the map has local degree
 $\ge 2$ at $z_0$, hence no neighborhood of $z_0$ within $\overline\D$
can map into $\overline\D$.
\end{proof}

\begin{rem} More precisely, whenever there is an interior fixed point,
it is not difficult to show that the induced map on
$\partial \D$ is topologically conjugate to the linear map $t\mapsto td$ of the
circle $\R/\Z$. To prove this, construct a new homeomorphism
$$t\,:\,\partial\D~\to~\R/\Z$$ as follows. Choose one boundary fixed point
$z_0$ and assign it the coordinate $t(z_0)=0$. The $d$ immediate
pre-images of $z_0$ divide $\partial\D$ into $d$ disjoint half-open arcs
 $A_0,\,A_1,\,\ldots,\, A_{d-1}$, numbered in counterclockwise
 order starting and ending at the point $z_0$. Note that each of these
arcs maps bijectively onto the entire circle. Now define the function
$t:\partial\D\to\R/\Z$ by the formula
$$t(z)~=~\sum_{k=0}^\infty a_k(z)/d^{k+1}\,,$$
 where  the integers $0\le a_k(z)<d$ are defined by the condition
\hbox{$\beta^{\circ k}(z)\,\in\,A_{a_k(z)}$}. Then evidently
 $a_k\big(\beta(z)\big)=a_{k+1}(z)$, and hence
$$ t\big(\beta(z)\big)~=~\sum_{k=0}^\infty a_{k+1}(z)/d^{k+1}~=
~t(z)\,d -a_0(z)~\equiv~ t(z)\,d~~ ({\rm mod}~\Z)\,.$$
Using the condition that
$d\phi/d\theta$ is greater  than some constant $c>1$, it is not
difficult to prove that the first $k$ terms of this series determine
the point $z$ to an accuracy of $2\pi/c^k$. It then follows easily that
the function $z\mapsto t(z)\in\R/\Z$ is a homeomorphism, as required.
\end{rem}
\medskip

In studying polynomial maps, we concentrated on those which are monic 
and centered. As a substitute for the monic condition, let us say
 that a Blaschke product $\beta$
is \textbf{\textit~1-anchored~} if $\beta(1)=1$. However, we will need
two different concepts of centering, depending on whether we are dealing with
 a periodic point or an aperiodic point of $|S|$.

\begin{definition}\label{d-cent}
Let $\beta:\D\to \D$ be a proper holomorphic map of
degree \hbox{$d\ge 1$.} We will say that $\beta$ is
{\textbf{\textit{fixed point centered~}}} if $~\beta(0)=0$,~
and 
\textbf{\textit{~zeros centered~}} if the sum
 $$a_1\,+\,\cdots\,+\,a_d$$ of the points in $\beta^{-1}(0)\,$ (counted with
multiplicity) is equal to zero.
(In the case $d=1$, note that the only Blaschke product which is 1-anchored,
 and centered in either sense, is the identity map.)
\end{definition}


 In order to construct  an appropriate topological model
 for hyperbolic components with a given reduced mapping scheme,
we will need three lemmas.


\begin{lem}\label{l-fc}
Let $\beta$ be a Blaschke product of degree $d\ge 2$ which has
a fixed point $z_0$ in the open disk ${\mathbb D}$, and let $h$ be a M\"obius
automorphism of the unit disk. Then the conjugate  $\beta'=h^{-1}\circ\beta
\circ h$  is fixed point centered if and only if $h(0)=z_0$, and
is 1-anchored if and only if $h(1)$
is one of the $d-1$ fixed points of $\beta$
 on the boundary circle $\partial\D$.
Thus, for each such $\beta$, there are $d-1$ possible choices for $h$.
\begin{equation}\label{e-conj}
\xymatrix{{\mathbb D} \ar[d]^{\beta'} \ar[r]^h & {\mathbb D} \ar[d]^{\beta} \\
{\mathbb D} \ar[r]^h & {\mathbb D}}
\end{equation}
\end{lem}

The proof is immediate.\qed

\begin{lem}\label{l-cc}
Let $\beta$ be an arbitrary Blaschke product of degree $d\ge 1$,
and  let $h$ be a M\"obius automorphism. Then the composition
$\beta\circ h$ is 1-anchored if and only if $h(1)$ is  one of the $d$
points $z_1\in\partial\D$ for which $\beta(z_1)=1$.
For each such $z_1$, there is a unique choice of $h$ so that $\beta\circ h$
is zeros centered. 
\end{lem}

\begin{equation}
\xymatrix{{\mathbb D} \ar[dr]_{\beta\circ h}\ar[r]^h & {\mathbb D} \ar[d]^{\beta} \\
&{\mathbb D}}
\end{equation}

 The proof will depend on the following.

\begin{definition}\label{l-cbc}
 Given points $z_1\,,\,\ldots\,,\,z_k$ in a Riemann surface $W$ isomorphic
 to $\D$, it follows from Douady and Earle \cite[\S2]{DE} that
 there exists a conformal isomorphism $\eta: W \to \D$, unique up to a
rotation of $\D$, which takes the  $z_j$ to points with sum
 $~\eta(z_1)+\cdots+\eta(z_k)~$ equal to zero. By definition,
the pre-image \hbox{$\widehat z=\eta^{-1}(0) \in W$} is called the 
\textbf{\textit{~conformal barycenter~}} of the points \hbox{$z_1,\,\ldots,\,z_k
\in W$}. Evidently this conformal barycenter is uniquely defined.
\end{definition}
\smallskip

\begin{proof}[Proof of Lemma \ref{l-cc}] If $h:\D\to\D$ is a conformal
automorphism, note that $h$ maps the zeros of $\beta\circ h$
to those of $\beta$, and hence maps the conformal barycenter of the
zeros of $\beta\circ h$ to the corresponding barycenter
$\widehat a$ of the zeros $a_1,\ldots, a_n$ of $\beta$. In particular,
 it follows that $\beta\circ
h$ is zeros centered if and only if $h(0)=\widehat a$.
Using these facts, the proof is straightforward.
\end{proof}

 The notations $\B^d_{\rm fc}$ and $\B^d_{\rm zc}\,$
will be used for the topological
space consisting of all 1-anchored Blaschke products of degree $d$
which are respectively fixed point centered or
 zeros centered. For the special case $d=1$,
evidently $\B^1_{\rm fc}=\B^1_{\rm zc}$ consists of
 a single point, namely the identity map.

\begin{lem}\label{l-cell} Each of the two model spaces
 $~\B^d_{\rm fc}~$ and $~\B^d_{\rm zc}\,$,
is homeomorphic to an open cell of real dimension $2(d-1)$.
\end{lem}

\begin{proof} Let $\s_d(\C)$ be the 
$d$-{\textbf{\textit{fold symmetric product\/}}},
consisting of unordered $d$-tuples
$\{a_1\,,\,\ldots\,,\,a_d\}$ of complex
numbers. This can be identified with the complex affine space consisting of
all monic polynomials of degree $d$, under the correspondence
$$\{a_1\, ,\,\ldots\,,\,a_d\} \mapsto (z-a_1)\cdots (z-a_d)= z^d -\sigma_1 
z^{d-1} +\sigma_2z^{d-2}-\cdots +(-1)^d\sigma_d \,,  $$
where the $\sigma_j$ are the elementary symmetric functions of $\{a_1\,,
 \,\ldots\,,\,a_d\}$. Thus $\s_d(\C)$ is homeomorphic to \hbox{$\C^d\cong
\R^{2d}$}. Since $\C$ is homeomorphic to the 2-cell $\D$,
it follows that $\s_d(\D)$ is also homeomorphic to $\R^{2d}$.

Now consider the space $\B_{\,\rm fc}^{\;d}$
 consisting of 1-anchored
Blaschke products $\beta$ of degree $d$ which fix the origin. We
can write
$$\beta(z) \;=\; z\, \mu_{a_2}(z)\cdots \mu_{a_d}(z) $$
(taking $a_1=0$).
Evidently this space is homeomorphic to the symmetric product $\s_{\w}(\D)$
where $\w=d-1$, and hence is a topological cell, homeomorphic to $\R^{2\w}$.
\smallskip

To determine the topology of $\B^d_{\rm zc}$, we proceed as follows.
We show first that the subspace $\s_d(\overline \D)\subset \s_d(\C)$,
consisting of unordered $d$-tuples $\{a_1,\,a_2,\,\ldots,\,a_d\}$ with
 ${\rm max}_j(|a_j|)\leq 1$ is a closed topological $2d$-cell
 with interior equal to $\s_d(\D)$. In fact, for each
 $\{a_1\,,\,\ldots\,,\,a_d\}\in\s_d(\C)$ such that the maximum of the
 $|a_j|$ is equal to one, consider the half-line consisting
of points $\{ta_1\,,\, \ldots\,,\,ta_d\}$ with $t \ge 0$.
 The image of each such half-line in the space of elementary symmetric
 $n$-tuples is a curve consisting
of points $(t\sigma_1\,,\, t^2\sigma_2\,,\,\ldots\,,\,t^d\sigma_d)\in\C^d$,
which crosses the unit sphere of $\C^d$
exactly once,  since the function \hbox{$~t 
\mapsto \vert t\,\sigma_1\vert^2+ \cdots + \vert t^d\,\sigma_d\vert^2~$}
 is strictly monotone. Hence, stretching by an appropriate
factor along each such ray, we obtain the required homeomorphism from $\s_d
 (\overline \D)$ to the closed unit ball in $\C^d$. Using this construction we
 see also that the subspace of $\s_d(\D)$ consisting of unordered $d$-tuples 
with sum $~\sigma_1=a_1+\cdots+a_d~$ equal
 to zero is an open topological $2(d-1)$-cell.
Thus the set $\B^d_{\rm zc}$ consisting of Blaschke products of the form
$\beta(z)=\mu_{a_1}(z) \cdots \mu_{a_d}(z)$ with $a_1+\cdots+a_d=0$
is an open topological $2(d-1)$-cell.
\end{proof}

Combining the three previous lemmas, we can construct  an
appropriate topological model for hyperbolic components with a given
reduced mapping scheme.

\begin{definition}\label{d-b^S}
To any mapping scheme $S=(|S|\,,\,F\,,\,\w)$
we associate the {\textbf{\textit{model space\/}}}~ $\B^S$~ consisting of all
 proper holomorphic maps
$$ \bbeta:|S|\times \D\;\to\; |S|\times \D $$ such that $\bbeta$ carries
each $s\times \D$ onto $F(s)\times \D$ by a 1-anchored 
 Blaschke product
$$ (s,\,z)~\mapsto \big(F(s)\,,~\beta_s(z)\big) $$
of degree $d(s)=\w(s)+1\;$ which  is either fixed point centered
or zeros centered 
according as $s$ is periodic or aperiodic under $F$. (Thus,
in the special case of
a vertex $s$ of weight zero,  we require $\beta_s$ to be the identity map.)
\end{definition}\smallskip

\begin{lem}\label{l-top-cell}
 If the scheme $S$ has total weight $\w(S)$, then the model space
$\B^S$ is homeomorphic to an open cell of dimension $2\w(S)$. Furthermore,
 $\B^S$ is canonically homeomorphic to $\B^{\overline S}$, where
$\overline S$ is the associated reduced mapping scheme.
\end{lem}
\smallskip

\begin{proof} The first statement follows immediately from Lemma~\ref{l-cell}
since, as a topological space, $\B^S$ is simply a Cartesian product of
spaces of the form $\B^{\,d(s)}_{\rm fc}$ and $\B^{\,d(s)}_{\rm zc}$ of
dimension $2\w(s)=2(d(s)-1)$, where
 $\sum \w(s)=\w(S)$. The second statement follows
since the space $\B^1_{\rm fp}=\B^1_{\rm zc}$ is a single point. \end{proof}

We will show that the various maps in $\B^S$ serve as models for the
dynamics of all possible hyperbolic components $H$ with $S_H\cong S$.
As a preliminary step, given a mapping scheme $S$, first consider more
 general maps
$$\bbeta:|S|\times\D\to |S|\times\D$$
 which carry each $s\times\D$
onto $F(s)\times\D$ by a Blaschke product of degree $d(s)$, but with no
other restriction. Evidently, each such 
$\bbeta$ extends uniquely over the union
$|S|\times\overline\D$ of closed disks.

\begin{definition}\label{d-bm}
 By a \textbf{\textit{~boundary marking~}} $q~$ for
$\bbeta$ we will mean a function $~s\mapsto q(s)\,\in\, s\times\partial\D~$
which assigns a boundary point to each $s\times\overline\D$,~
and which satisfies the condition that
\begin{equation}\label{e-bm}
 q\big(F(s)\big)~=~ \bbeta\big(q(s)\big)\qquad{\rm
for~ every}\quad s\,.
\end{equation}
\end{definition}

\begin{lem}\label{l-bm}
Given $\bbeta$ as above, the number of possible boundary markings $q$
is equal to the order of the automorphism group $\G(S)$, as computed in 
Lemma~$\ref{l-|G|}$. In particular, such boundary markings always exist.
\end{lem}

\begin{proof}
First consider a vertex $s$ which is periodic,
 $F^{\circ k}(s)=s$. If $d_1d_2\cdots d_k$ is the product of the degrees
around this periodic orbit, then we can choose $q(s)$ to be any one of the
$d_1\cdots d_k-1$ fixed points of $\,\bbeta^{\circ k}$ on 
$s\times\partial\D$.
(Compare Lemma~\ref{l-3.1}.) The choice of $q\big(F(s)\big)$ is then
determined by equation (\ref{e-bm}), and we can continue inductively
around the cycle. In the case of an aperiodic vertex, suppose
inductively that $q\big(F(s)\big)$ has already
been chosen, but $q(s)$ has not. Then there are $d(s)$ possible
choices for $q(s)$, again by Lemma~\ref{l-3.1}.
 Further details of the proof are straightforward.
\end{proof}

\begin{theorem}\label{t-hbp}
As above, let $\bbeta:|S|\times\D\to |S|\times\D$ carry each $s\times\D$
onto\break $F(s)\times\D$ by a proper holomorphic map of degree $d(s)$. 
Suppose that $\bbeta$
is \textbf{\textit{~hyperbolic~}} in the sense that every orbit in
 $|S|\times\D$ converges to an attracting cycle in $|S|\times\D$. Then
for every boundary marking $q$ there
exists a unique automorphism $\bf h$ of $|S|\times\overline\D$ such that
the conjugate map $~{\bf h}^{-1}\circ\bbeta\circ{\bf h}~$ belongs to the space
 $\B^S$ of Definition $\ref{d-b^S}$, with ${\bf h}(s,1)=q(s)$.
\end{theorem}

\begin{proof}
If we choose $\bf h$ so that ${\bf h}(s,1)=q(s)$ for every $s$, then it
is straightforward to check that $~{\bf h}^{-1}\circ\bbeta\circ{\bf h}$
 is 1-anchored. We must show
that there is then a unique choice of the values ${\bf h}(s,\,0)$ so that
${\bf h}^{-1}\circ\bbeta\circ{\bf h}$
 also satisfies the appropriate centering conditions.
If $s$ is periodic under $F$, then since $\bbeta$ is hyperbolic, it
follows that $s\times\D$ contains a necessarily unique
attracting periodic point $(s,\,z_s)$ and we can
 choose ${\bf h}$ so that ${\bf h}(s,0)=(s,\,z_s)$. On the
other hand, if $s$ is not periodic, then by Lemma~\ref{l-cc},
 assuming inductively that $\bf h$
has already been defined on $F(s)\times\D$, we
can choose the automorphism ${\bf h}$ on $s\times\D$ so that
 $~{\bf h}^{-1}\circ\bbeta\circ {\bf h}~$ is zeros centered on $s\times\D$.
 The rest of the argument is straightforward.\end{proof}

\begin{definition}\label{d-Uf} 
If  $\f$ is  a hyperbolic map in the connectedness locus of the space $\p^{S_0}$
for some scheme $S_0$, let $~\U_\f\subset|S_0|\times\C~$ be the union of
 those Fatou components of $\f$ which contain critical or postcritical points.
\end{definition}

\begin{coro}\label{c-U=model}
Let $\f$ and $\U_\f$ be as above. 
Then the map $\f$ restricted to $~\U_\f~$
is conformally conjugate to some map $~\bbeta:|S_\f|\times\D\to|S_\f|\times\D~$
belonging to the model space $\B^{S_\f}$. A similar assertion holds for
 rational maps which are hyperbolic, with connected Julia set.
\end{coro}

\begin{proof} Since each connected component of $\U_\f$ is conformally
 isomorphic to the unit disk, we can choose a conformal isomorphism
from $\U_\f$ to $|S_\f|\times\D$. The conclusion then follows easily from
Lemma \ref{l-bm} and Theorem \ref{t-hbp}.
\end{proof}


We will also need the following elementary result.
\smallskip

\begin{lem}\label{l-cf}
 Each model space $\B^S$ contains one and only one map $\bbeta_{0}$
 which is critically finite, given by the formula
$$\bbeta_{0}(s,z)~=~\big(F(s),\,z^{d(s)}\big)\,. $$
\end{lem}

\begin{proof}
  First consider a Blaschke product $\beta : \D
 \to \D$ of  degree $n\ge 2$, which is critically finite. Then $\beta$
certainly has a periodic point, say of period $k$. If $k>1$ then the
$k$-fold iterate of $\beta$ would have $k$ distinct fixed points,
which contradicts Lemma \ref{l-3.1}; therefore, $k=1$. After conjugating by a
M\"obius automorphism, we may assume that this fixed point lies
at the origin. We must then prove that the origin is the only critical
point. If $\beta'(0)$ were non-zero, then we could choose a K{\oe}nigs
 coordinate, defined
throughout a maximal open set $U\subset\D$ which maps diffeomorphically
onto a round disk. The boundary $\partial U$ would then have to contain
 a critical point with infinite orbit, which contradicts the hypothesis.
On the other hand, if $\beta'(0)=0$, then one can choose a B\"ottcher
coordinate mapping a maximal open set $U$ diffeomorphically onto a round disk.
If $U\ne \D$, then $\partial U$ would again contain a critical point with
 infinite orbit, contradicting the hypothesis. Thus the origin is the only
critical point. After conjugating by a rotation,
 it follows that $\beta(w)=w^n$.

Now consider a critically finite $\bbeta\in\B^S$. If $s\in|S|$ is a periodic
vertex of period $k$, then the argument above shows that
 $\bbeta^{\circ k}(s,\,w)=(s,\,w^n)$, where $n$ is the degree of
 $\bbeta^{\circ k}$ on $s\times\D$. It follows
that $\bbeta$ has no critical points on $s\times\D$ other than $(s,\,0)$;
therefore $\bbeta(s,\,w)=(s,\,w^{d(s)})$. Now consider an aperiodic vertex $s$.
Assuming inductively that the only periodic or preperiodic point in
 $~F(s)\times\D~$ is $~\big(F(s),\,0\big)$, it follows that $\bbeta$ must map
 every critical point in $s\times\D$ to $\big(F(s),\,0\big)$. Let $s\times X$
be the finite set consisting of all pre-images of $\big(F(s),\,0\big)$ in
$s\times\D$. Then $s\times(\D\ssm X)$ is an unbranched $d(s)$-fold
 covering space of $F(s)\times(\D\ssm\{0\})$. Therefore $\D\ssm X$ is
 conformally isomorphic to a punctured disk. This implies that $X$ is
a single point, which must be the origin since this map is zeros centered.
This completes the induction, and hence completes the proof.\end{proof}

\section{Hyperbolic Components are Topological Cells.}\label{s4}

Given any mapping scheme $S_0$, consider the associated space
$\p^{S_0}$ of polynomial maps (Definition \ref{d-p^S}).
Let $H\subset\p^{S_0}$ be a hyperbolic component in the connectedness
locus of $\p^{S_0}$.

As in Definition~\ref{d-Uf} for each $\f\in H$,
let ${\U}_\f\subset K(\f)$ be the union of those Fatou components of $\f$
 which contain critical or postcritical points.\smallskip

The object of this section is to prove the following result.

\begin{theorem}\label{t-top-mod}
Let $S$ be the full mapping scheme associated with some representative
map in $H$. Then there exists a diffeomorphism $~H\stackrel{\cong}
{\longrightarrow}\B^S~$
which sends each $~\f\in H~$ to a map
 $\bbeta(f):|S|\times\D\to |S|\times\D$ in $\B^S$
which is conformally  conjugate
to the restriction $~\f|_{\U_\f}:\U_\f\to \U_\f$.
\end{theorem}

\noindent In fact the proof will show that there exist only finitely many
 such diffeomorphisms, where the number of possible choices is equal to the
order of the automorphism group $\G(S)/\G_0(S)$.

Combining this statement with Lemmas \ref{l-top-cell} and \ref{l-cf},
we immediately obtain the following result, which generalizes an
 unpublished theorem of McMullen.\smallskip

\begin{coro}\label{c-top-mod}
Every such hyperbolic component $H$ is a topological cell of dimension
$2\,\w(S)$; and every such $H$ contains a unique critically finite map.
\end{coro}

We must be careful in the proof of Theorem~\ref{t-top-mod} since it is
not a priori clear that $H$ is simply connected. As an example to
illustrate the difficulty, each
$\f\in H$ has a mapping scheme $S_\f$. Following a path from
$\f$ to $\f'$, we obtain a well defined isomorphism from $S_\f$ to
$S_{\f'}$. However, we must check that this isomorphism does not
depend on the choice of path.\smallskip




In analogy with Definition \ref{d-bm}, we introduce the concept
of boundary markings for hyperbolic maps. Recall that
 the boundary of each component of ${\U}_\f$ is a Jordan curve.\footnote
{{\sc Proof.} (For a more general result, see \cite{RY}.)
It suffices to consider the classical case of a polynomial map
 $f:\C\to\C$. Since $f$ is hyperbolic with connected Julia set,
 its Julia set is locally connected. Therefore, for any bounded Fatou component
$U$, a conformal equivalence $\D\stackrel{\cong}{\longrightarrow} U$, extends
to a continuous map $\overline\D\to \overline U$. If two points
$e^{i\theta}$ and $e^{i\theta'}$ in $\partial\D$ mapped to the same point of
 $z\in\partial U$, then the
broken line from $e^{i\theta}$ to $0$ to $e^{i\theta'}$ would map
 to a simple closed curve $\Gamma\subset\C$. By the maximum modulus principle,
the bounded component of the complement of $\Gamma$ must lie in the interior
of the filled Julia set, and hence must be contained in $U$. But this would
imply that there is an entire interval of angles, say with
 $\theta\le\phi\le\theta'$, so that $e^{i\phi}$ maps to $z$. This is
 impossible by a theorem of Riesz and Riesz. (See for example
\cite[\S\S17.14, 19.2 and A.3]{M4}.) }

\begin{definition}\label{d-bm2}
A \textbf{\textit{~boundary marking for a hyperbolic map~}}
$\f\in\cl^{S_0}$ will mean a function $q$ which assigns to each 
connected component
$U\subset{\U}_\f$ a boundary point $~q(U)\in\partial U~$
so as to satisfy the identity $$q\big(\f(U)\big)= \f\big(q(U)\big)\,.$$
\end{definition}

\begin{definition}\label{d-Htilde} Let $S$ be the full mapping scheme for some
representative map in $H$, and let
 $\widetilde H$ be the set of all triples consisting of
\begin{itemize}
\item a map $~\f\in H$,
\item a boundary marking $~q~$ for $~\f$,~ and
\item an isomorphism $~\biota:S_\f~\cong~S$.
\end{itemize}
\end{definition}

\begin{lem}\label{l-tildeH} This set $\widetilde H$ has a natural topology
so that every point of $H$ has a neighborhood $N$ which is evenly
 covered\,\footnote{By definition,  $N$ is \textbf{\textit{\,evenly covered\,}}
if each connected component of $p^{-1}(N)$ maps homeomorphically onto $N$.}
under the projection $\widetilde H\to H$.
\end{lem}

\begin{proof} This is straightforward. In fact, 
each point $~q(U)~$ is preperiodic and eventually repelling,
and therefore deforms continuously as we deform the map $\f$.
Similarly the isomorphism $\biota$ deforms continuously with $\f$.
\end{proof}

It follows that every connected component of $\widetilde H$ is a
(possibly trivial) covering space of $H$.
(It also follows that we can lift the complex structure from $H$, so that
the projection map is locally biholomorphic.)

\smallskip

Next we project this space $\widetilde H$ onto the model space $\B^S$.

\begin{lem}\label{l-UcongB}
To every $(\f,\,q,\,\biota)\in\widetilde H$, there is uniquely associated a
map $$\bbeta=\pi(\f,\,q,\,\biota)\in \B^S\,,$$ together with a conformal
 conjugacy between the restriction $~\f|_{\U_\f}:\U_\f\to \U_\f~$ and the
 map  $~\bbeta:|S|\times\D\to|S|\times\D$.
\end{lem}

\begin{proof}
We use the isomorphism $\biota$ to identify $S_\f$ with $S$.
Start with some arbitrary conformal isomorphism which carries each component
 $U\subset \U_\f$ onto the corresponding $s_U\times\D$. Then the boundary
 marking $q(U)\in\partial U$ will correspond to a boundary marking in
$s_U\times\partial\D$. We can then use Theorem~\ref{t-hbp}
 to obtain a corrected conformal isomorphism $\U_\f\to|S|\times\D$
which is actually a conformal conjugacy between $\f|_{\U_\f}$ 
 and a corresponding map $\bbeta=\pi(\f,\,\biota,\,q)\in\B^S$.
\end{proof}

\begin{theorem}\label{t-ev-cov} The resulting projection $\pi:\widetilde H
\to \B^S$ is continuous. Furthermore, every point
 $\bbeta\in\B^S$ has a neighborhood $N$ which is evenly covered.
\end{theorem}

\begin{rem}\label{r-proof} If we assume this theorem for the moment,
then the main results of this section follow easily.
Since $\B^S$ is a topological cell by Lemma~\ref{l-top-cell},
it is certainly simply-connected. Thus Theorem \ref{t-ev-cov} implies
 that every connected component
of $\widetilde H$ maps homeomorphically onto $\B^S$. Since $\B^S$
has a unique critically finite point by Lemma~\ref{l-cf}, it follows that
each connected component of $\widetilde H$ also has a unique critically finite
point. On the
other hand, each connected component of $\widetilde H$ is a covering
space of $H$ by Lemma~\ref{l-tildeH}. In fact it must actually map
homeomorphically onto $H$; for if it were a non-trivial covering space, then
$\widetilde H$ would have more than one critically finite point.

Thus, choosing any section $H\to\widetilde H$, the composition
 $$H~\longrightarrow~\widetilde H~\stackrel{\pi}{\longrightarrow}~\B^S$$
 maps $H$ homeomorphically
onto the topological cell $\B^S$. This shows that Theorem~\ref{t-top-mod} and
 Corollary \ref{c-top-mod}, as stated at the beginning of this section,
follow immediately from Theorem  \ref{t-ev-cov}.\end{rem}\medskip

The proof of Theorem \ref{t-ev-cov} will make use of the following.

\begin{lem}\label{l-choose-r}
Given $\bbeta\in\B^S$, it is possible to choose a radius $0<r(s)<1$
for each $s\in |S|$ satisfying the following two conditions: 

\begin{enumerate}
\item Every critical point of $\bbeta$ in $s\times\D$
is contained in $s\times\D_{r(s)}$, where \break
 \hbox{$\D_r~=~\{w\in\C~;~|w|<r\}$}
denotes the open disk of radius $r$.
\item The image of the closure of this disk under $\bbeta$  satisfies
\begin{equation}\label{e-r-cond}
\bbeta\big(s\times\overline\D_{r(s)})~\subset ~F(s)\times\D_{r(F(s))}\,.
\end{equation}
\end{enumerate}\end{lem}


\begin{proof} 
Start with the aperiodic vertices. If $s$ belongs to the set\break
 \hbox{$S'=S\ssm F(S)$~} of ``free'' vertices, then  any $r(s)$
sufficiently close to one will do. Next choose $r(s)$ for the aperiodic vertices
in $F(S')$, and continue inductively. Once $r(s)$ has been chosen for all
aperiodic vertices, the remaining choices are  not difficult. In fact all
 of the maps around a cycle are fixed point centered. Therefore, for
 $s$ periodic,
$${\rm if}\qquad \bbeta(s,\,w)=\big(F(s),\,w'\big)\,,\qquad{\rm then}
\qquad |w'|~\le~|w|\,,$$
 with strict inequality whenever $w\ne 0$ and $d(s)\ge 2$ by the Schwarz
Lemma. Further details are straightforward, since every cycle contains at
 least one vertex of degree $\ge 2$.\end{proof}

It is now easy to choose radii $R(s)$ slightly larger than $r(s)$ so that
\begin{equation}\label{e-R-cond}
\bbeta\big(s\times\overline\D_{R(s)})~\subset ~F(s)\times\D_{r(F(s))}\,.
\end{equation}
The annuli
$$ A(s)~=~s\times\big(\D_{R(s)}\ssm\overline\D_{r(s)}\big) $$
will play an important role. Note the
crucial property that no orbit under $\bbeta$ can pass through
the union  $~\bigcup_s A(s)~$ more than once.\smallskip

\begin{figure}[htbp]
\begin{center}
\psfrag{B}[c][c][1.2]{${\mathcal B}$}
\psfrag{i}[c][c][1.2]{${\biota}$} 
\psfrag{b}[c][c][1.2]{${\bbeta}~$}
\centerline{\psfig{figure=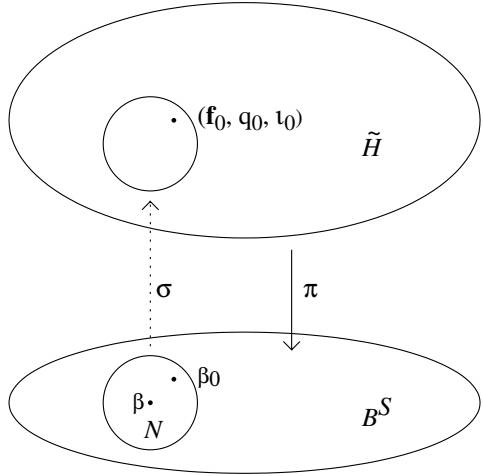,height=2.5in}}
\caption{\label{f-5.8}\it Proof of the even covering property.}
\end{center}
\end{figure}

{\sc Proof of Theorem \ref{t-ev-cov}.} Recall that $H\subset\p^{S_0}$
 is a hyperbolic component with mapping scheme $S$; that $\widetilde H$
is the finite covering space of $H$ consisting of triples
 $(\f,\,q,\,\biota)$; and that
$$\pi:\widetilde H\to\B^S$$
is the associated projection map. Given any $\bbeta\in\B^S$, we must find a
neighborhood $N$ of $\bbeta$ which is evenly covered. This means that, given
any $(\f_0,\,q_0,\,\biota_0)\in\pi^{-1}(N)$, and setting
$\bbeta_0=\pi(\f_0,\,q_0,\,\biota_0)\in N$, we must find a
 section $\sigma :N\to \widetilde H$ such that
$\sigma(\bbeta_0)=(\f_0,\,q_0,\,\biota_0)\,$,
with $\pi\circ\sigma$ equal to the identity map of $N$.
To achieve this, we will impose two conditions on $N$:
\medskip

{{\sc Condition 1.} \it This neighborhood $N$ must be small enough so that the
 conditions $(1)$ and $(2)$ of Lemma $\ref{l-choose-r}$, as
well as inequality $(\ref{e-R-cond})$, hold 
 with the same choice of radii $r(s)$, $R(s)$ for all 
 $\bbeta'\in N$.}\medskip

Let ${\U}_{\f_0}$ be the union of all critical
and postcritical Fatou components for $\f_0$.
Using Lemma \ref{l-UcongB}, we can identify
 $|S|\times\D$ with this open set ${\U}_{\f_0}\subset|S_0|\times\C$.
 Furthermore, under this
identification, the map $\bbeta_0$ from $|S|\times\D$ to itself
corresponds to the map $\f_0$ from $\U_{\f_0}$ to itself.\smallskip

For any $\bbeta_1\in N$, the
 image $\sigma(\bbeta_1)\in\pi^{-1}(\bbeta_1)\subset\widetilde H$ will be
 constructed by quasiconformal surgery. (Compare \cite{DH2}.)
The first step is to construct a
preliminary map $~\widehat\f_1~$ from $|S_0|\times\C$  to itself as  follows. Set
$~\widehat\f_1(s,\,w)=\bbeta_1(s,\,w)~$ whenever $(s,\,w)$ belongs to the small disk
$$ s\times\D_{r(s)}~\subset~|S|\times\D~\cong~{\U}_{\f_0}~
\subset~|S_0|\times\C\,.$$
On the other hand, let $~\widehat\f_1~$ coincide with  $\bbeta_0$
 outside the union of larger disks 
$$ \bigcup_s\; s\times\D_{R(s)}~\subset~|S|\times\D~\cong~{\U}_{\f_0}~
\subset~|S_0|\times\C\,.$$
 Within the intermediate closed annuli $~\overline A(s)
=s\times\big(\overline\D_{R(s)}\ssm\D_{r(s)}\big)$, we interpolate linearly, setting
$$ \widehat\f_1(s,\,w)~=~t\,\bbeta_0(s,\,w)\,+\,(1-t)\,\bbeta_1(s,\,w)\,,
\qquad{\rm where}\qquad t~=~\frac{|w|-r(s)}{R(s)-r(s)}\,.$$
We can now impose the second condition: \medskip

{{\sc Condition 2.} \it The neighborhood $N$ must be small enough so that,
for each $\bbeta_1$ in $N$,
the map $\widehat\f_1$ defined in this way has Jacobian determinant bounded
away from zero throughout each $\overline A(s)$.}\medskip

Next we will use quasiconformal surgery to construct a new conformal
structure on $S_0\times\C$ which is  $\widehat\f_1$ invariant. To do this,
start with the standard (quasi-) conformal structure on the small disks
$s\times\D_{r(s)}$, and also on all points of \hbox{$|S_0|\times\C$} which are
 not in the iterated pre-image of
${\U}_{\f_0}$. Now pull this quasiconformal structure back to the rest
of ${\U}_{\f_0}$ under the action of $\widehat\f_1$ and its iterates. This
will yield a well defined quasiconformal structure on $|S_0|\times\C$,
which has bounded dilatation since an orbit can pass through the union
of annuli $A(s)$ at most once. Using the measurable Riemann mapping theorem,
we can choose a straightening map  $\eta$ which carries each $s\times\C$
 to itself, and
which carries our exotic quasiconformal structure to the standard conformal
structure.  This implies that the map 
$~\f_1=\eta^{-1}\circ\widehat\f_1\circ\eta~$
is holomorphic with respect to the standard conformal structure. Now,
after composing $\eta$ with suitable component-wise affine transformations,
 we may assume that
$\f_1$ is  1-anchored and centered.
 Using the Ahlfors-Bers measurable Riemann
mapping theorem with parameters \cite{AB}, we can choose these affine
transformations so that $\f_1$ varies continuously
as $\bbeta_1$ varies over the neighborhood $N$. Now we define the section
$\sigma:N\to\widetilde H$ by setting $\sigma(\bbeta_1)=
(\f_1,\,q_1,\,\biota_1)$ where $\biota_1$ is constant and where the boundary
marking $q_1$ also varies continuously with $\bbeta_1$.\medskip

Next, we must prove that $\pi(\f_1,\,q_1,\,\biota_1)=\bbeta_1$. The proof
will depend on the following.

\begin{lem}\label{l-con-cl}
The conformal conjugacy class of a map $\bbeta_1\in N$ is uniquely determined
by the conformal conjugacy class of the restriction of $\bbeta_1$ to the union
 $\bigcup_s s\times\D_{r(s)}$ of subdisks.
\end{lem}

\begin{proof} Define sets $~\D_k(s)\subset s\times\D~$ inductively by setting
$~\D_0(s)=s\times\D_{r(s)}$ and
$$\D_{k+1}(s)~=~(s\times\D)~\cap~\bbeta_1^{-1}\D_k\big(F(s))~.$$
Then it is not hard to check that each $\D_k(s)$ is an open topological disk
with smooth boundary, and that
 $~\D_0(s)\subset\D_1(s)\subset\cdots$, with union $~s\times\D$. Because all
the critical values are well inside the disks,
 each $\D_{k+1}(s)$ can be described conformally as a $d(s)$-fold
branched covering of  $\D_k\big(F(s)\big)$, where the nontrivial branching
already occurs in the subset $\D_0(s)$. Thus we can build these sets up
inductively, starting only with the restriction of $\bbeta_1$ mapping
 $~\bigcup_s\D_0(s)~$ into itself.
 The union is the required conformal dynamical
system, conformally conjugate to $~\bbeta_1:|S|\times\D\to|S|\times\D$.
\end{proof}

\begin{proof}[ Proof of Theorem \ref{t-ev-cov}, conclusion.] Applying 
Lemma \ref{l-con-cl}
to the map $\f_1$ as constructed above, it follows that
 $\pi(\f_1,\,q_1,\,\biota_1)\in\B^S$ is conformally conjugate to $\bbeta_1$.
Since the boundary marking varies continuously with $\bbeta_1$, 
 this implies that $\pi(\f_1,\,q_1,\,\biota_1)=\bbeta_1$ as required.

Finally, we must prove that the projection $~\pi:\widetilde H\to \B^S~$ is 
also continuous. But the map $\sigma:N\to\sigma(N)\subset\widetilde H$
constructed above is known to be continuous and one-to-one; hence it maps
any compact subset of $N$ homeomorphically. Since the projection $\pi$ can
be identified locally with $\sigma^{-1}$, it follows that $\pi$ is continuous.
 This completes the proof
of Theorem \ref{t-ev-cov}, and hence of Theorem \ref{t-top-mod}
 and Corollary \ref{c-top-mod}.
\end{proof}

\begin{figure}[ht!]
\centerline{\psfig{figure=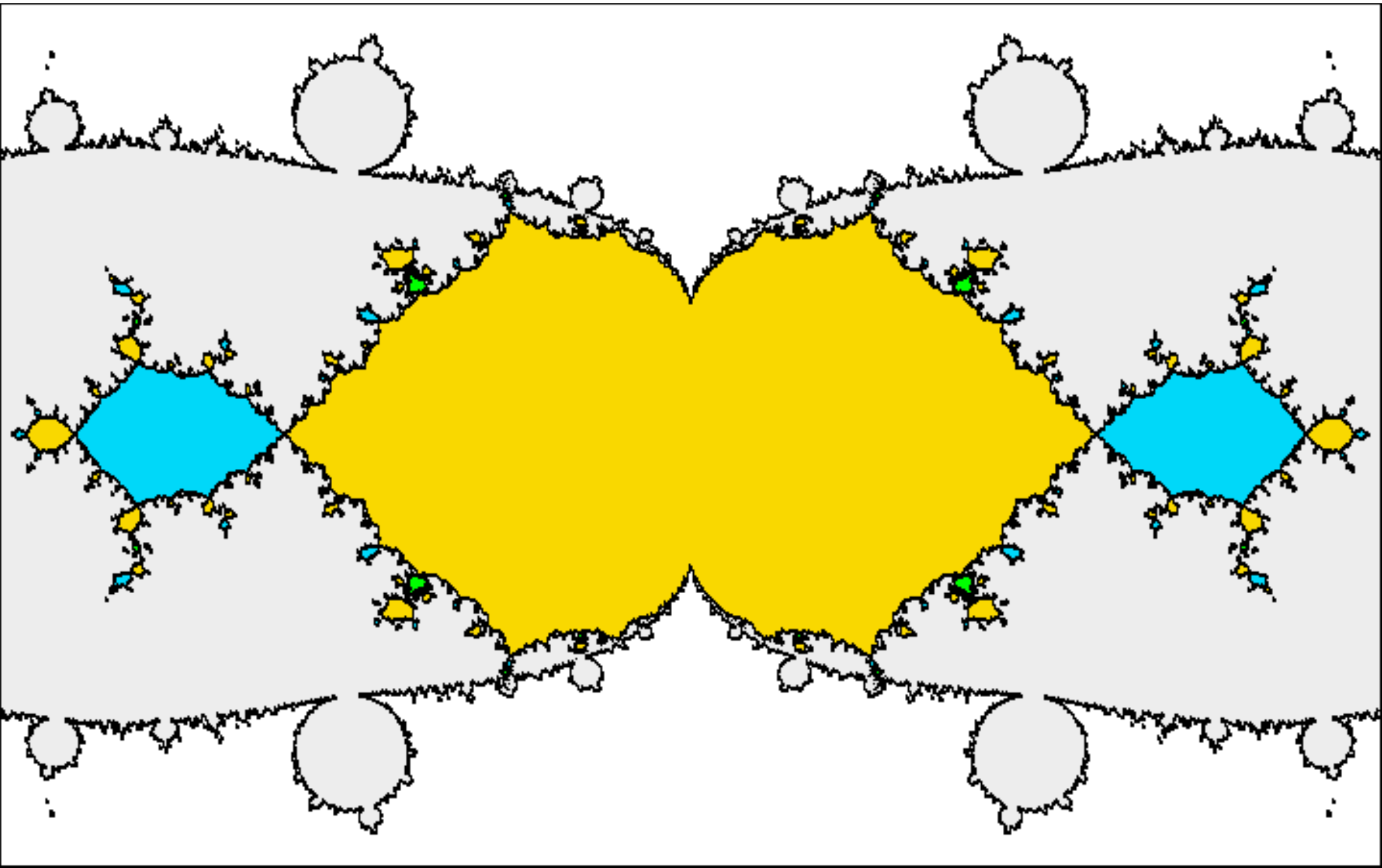,height=1.9in}}
\caption{\label{f-cubslice} \it The $b$ parameter plane for the family
of cubic maps 
$z\mapsto z^3-1.5\,z+b$.
The large central region corresponds to a hyperbolic component with
full mapping scheme of the form \hbox{$~~\bullet\leftrightarrow\bullet~~$} of
bitransitive type, while the two side regions have schemes of the form ~~
\hbox{~~$\bullet$\;$\rightarrow$\;{\large$\cdot$}\;$\rightarrow$\;$\bullet\;
\leftrightarrow$\;{\large$\cdot$}~~}
of capture type. $($At either of the two common boundary points
$b=\pm 0.4202\cdots$,  the attracting
period two orbit persists, but one critical point becomes preperiodic.$)$
In the surrounding light grey region, one critical orbit escapes to infinity,
while in the outer white region both critical orbits escape.}
\end{figure}

\begin{rem}\label{r-bndry}
In degrees $\ge 3$, one cannot expect the boundaries of hyperbolic
components to be smooth manifolds.
 (See Figure \ref{f-cubslice}, and compare \cite{PTL}.) In fact, it is not at
all clear that the boundary of every hyperbolic component  
must be a topological sphere. (Compare the discussion in Example \ref{ex-sc}.)
\end{rem}

\begin{rem}\label{r-subcurve}
It is often useful to restrict attention to some holomorphic subvariety
 of the parameter space $\p^{S_0}$. In general, one can't
expect hyperbolic components in such a parameter subspace to satisfy
Corollary \ref{c-top-mod}. For example, the parameter slice shown in
 Figure \ref{f-cubslice} contains no critically finite points.
 However, in the special case where the subspace
is defined by requiring one or more critical points to be periodic of specified
period, the proof can be adapted as follows. 
Let $H\subset\p^{S_0}$ be a hyperbolic component with mapping scheme $S$.
For each periodic vertex $s$ of $S$,
 write the weight $\w(s)$ as a sum $\w'(s)+\w''(s)$ where $\w'(s)$
is to be the number of free critical points  in the
 corresponding Fatou component, and where the periodic point in this component
is required to be a critical of point multiplicity at least $\w''(s)$
whenever $\w''(s)>0$. In the
Blaschke product model for $H$,
this means that the maps $\bbeta$ must have the form
$$ \bbeta(s,\,z)~=~\big(F(s)\,,
~z^{\w''(s)+1}\mu_{a_1}(z)\cdots \mu_{a_{\w'(s)}}(z)\big)~.$$
It is straightforward to check that the subspace of
${\mathcal B}^{\w(s)}_{\rm fc}$ defined in this way is a topological cell
of dimension $2\w'(s)$ with a preferred center point. Since
the full model space ${\mathcal B}^S$ is homeomorphic
 to a cartesian product of cells
 ${\mathcal B}^{\w(s)}_{\rm fc}$ with $s$ periodic, together with cells
${\mathcal B}^{\w(s')}_{zc}$ with $s'$ aperiodic, it follows that the
subspace of ${\mathcal B}^S$ defined by all of  these conditions is a
 topological cell with dimension twice the number of free critical points.
(Here all of the critical points associated with an aperiodic vertex
are free by definition.) Furthermore, this cell contains a unique
critically finite point. The analogous statements for the corresponding
 subspace of $H$ then follow from Theorem \ref{t-top-mod}.

Note however that this argument works only within the given 
 hyperbolic component. If we try to form an analogous global subvariety of
$\p^{S_0}$, we must first mark one or more critical points, which 
will usually change 
the global topology. (Compare Remark \ref{r-slice}.)
\end{rem}

\section{Analytic Isomorphism between
 Hyperbolic Components.}\label{s5}

If $H_\alpha \subset \cl^{S_1}$ and $H_\beta \subset \cl^{S_2}$ are
two different hyperbolic components with reduced mapping scheme
isomorphic to $S$, then by Theorem~\ref{t-top-mod} there are diffeomorphisms
$$H_\alpha~ \stackrel{\cong}{\longrightarrow}~\B^S
~ \stackrel{\cong}{\longleftarrow}~ H_\beta \,,$$
uniquely defined up to a choice among finitely many boundary markings, or
equivalently up to the action of the group $\G(S)/\G_0(S)$ on $\B^S$. The
composition mapping $H_\alpha$ to $H_\beta$ will be called a
{\textbf{\textit{canonical diffeomorphism\/}}} between these two sets. We
will prove the following.

\begin{theorem}\label{t-5.1}
 This canonical diffeomorphism $H_\alpha \to
 H_\beta$ between open subsets of complex affine spaces is biholomorphic.
\end{theorem}

\begin{definition}\label{d-s-mod}
As \textbf{\textit{~standard model~}} for hyperbolic components with
scheme $S$ we can take the hyperbolic component $H_0^S\subset\p^S$
which is centered at the map $\f_0(s,\,z)=\big(F(s),\,z^{d(s)}\big)$.
\end{definition}
In particular, it follows from  Theorem \ref{t-5.1}
 that the canonical diffeomorphism from $H_\alpha$
to the standard model $H_0^S$ is biholomorphic. Note that this
diffeomorphism is unique up to the action of the finite group
$\G/\G_0$ of linear automorphisms of $H_0^S$. The proof of Theorem \ref{t-5.1}
 will be based on the following.\smallskip

\begin{definition}
  We will say that a map $\f \in \p^{S_1}$ satisfies
a {\textbf{\textit{critical orbit relation\/}}} if either 
\begin{enumerate}
\item the $\w$ critical points of $\f$ are not all distinct, or 
\item the associated critical orbits are not disjoint from each
 other, or 
\item some critical orbit is periodic or eventually periodic.
\end{enumerate}
 It is not difficult to show that
the set of all $\f$ which satisfy some critical orbit
relation forms a countable union of algebraic varieties in the affine space
$\p^{S_1}$. However, we can make a sharper statement for
 the hyperbolic subset of $\p^S$.
\end{definition}

\begin{lem}\label{l-open}
Let $Q_\alpha$ be the subset consisting of maps in $H_\alpha$ which have no
critical orbit relation. Then $Q_\alpha$ is a dense open subset of $H_\alpha$.
\end{lem}

\begin{proof} Given $\f\in H_\alpha$, as in Definition \ref{d-Uf}
let $\U_\f$ be the union of all 
Fatou components of  $\f$ which contain critical or postcritical points.
First consider the simplest case in which $\U_\f$ is connected.
In other words, assume that all of the
critical points of $\f$ lie in $\U_\f$, which must be the immediate attracting
basin of an attracting fixed point $p_\f$. This means that
 the associated full mapping scheme $S$  consists of a single vertex of
 weight $\w$. Note first, for any $\f\in Q_\alpha$, that the multiplier
 $\lambda_\f$ at $p_\f$ must be non-zero. For otherwise $p_\f$ would be a
fixed critical point. Thus we can choose a K{\oe}nigs
linearizing function $~\kappa_\f:\U_\f\to\C~$ which maps a neighborhood of
$p_\f$ biholomorphically onto a neighborhood of the origin, and satisfies
\begin{equation}\label{e-koen}
 \kappa_\f\big(\f(z)\big)~=~ \lambda_\f\,\kappa_\f(z)\,.
\end{equation}

Let  $~c_\f^{\;1},\,\ldots,\,c_\f^{\;\w}~$ be the critical points of $\f$.
(Here the superscripts are just labels; not exponents.)
 If $\f$ has no critical orbit
relations, then these critical points 
must be distinct, and can be chosen as functions which vary
holomorphically as $\f$ varies through a small
neighborhood. Furthermore, the values $\kappa_\f(c_\f^{\;j})$ must
all be distinct and non-zero. The condition (\ref{e-koen}) determines
the K{\oe}nigs function $\kappa_\f$ only up to a multiplicative constant.
It will be convenient to normalize this function 
so that $$\kappa_\f(c_\f^{\;1}\big)~=~+1\,.$$
 It is then not
hard to show that $\kappa_\f(z)$ is holomorphic as a function of both
$\f$ and $z$  throughout this small neighborhood.

For any $\lambda\in\D\ssm\{0\}$, let ${\mathbb T}_\lambda$ denote the compact
 torus
which is obtained by identifying each $z\in\C\ssm\{0\}$ with $\lambda\,z$,
and hence with all multiples of the form $\lambda^k\,z$.
 Then it is not difficult to check that a map $\f\in H_\alpha$ has no
critical orbit relations if and only if
\begin{quote}
\begin{enumerate}
\item $\lambda_\f~\ne~0\,$, and

\item the images of the
numbers $\kappa_\f(c_\f^{\;j})$ under projection to ${\mathbb T}_{\lambda_\f}$
are all distinct.
\end{enumerate}
\end{quote}

\noindent
Since these are both open and dense conditions, this proves 
Lemma~\ref{l-open}  in the
special case. The proof in the general case is completely analogous.
Just choose one periodic point in each attracting cycle, and work with the
associated multipliers, K{\oe}nigs functions, and compact tori. 
Details are left to the reader.
\end{proof}

\begin{lem}\label{l-candiff}
 For any $\f_1$ with no critical orbit relations
the canonical diffeomorphism from $H_\alpha$ to $H_\beta$ is biholomorphic
throughout a neighborhood of $\f_1$.
\end{lem}

\begin{proof}
Again, we first consider the special case where all of the critical points
lie in the immediate attracting basin of a single attracting fixed point.
With $\lambda_\f$ and the $\kappa_\f(c_f^{\;j})$ as above,
we will show that the mapping 
\begin{equation}\label{e-coords} \f~\mapsto~\big(\;
 \lambda_\f\,,~ ~  \kappa_\f(c_\f^{\;2})\,,~~\kappa_\f(c_\f^{\;3}))\,,~~
\cdots\,,~~ \kappa_\f(c_\f^{\;\w})\;\big)~\in~\C^\w~,
\end{equation}
constitutes a local holomorphic coordinate system
 as $\f$ varies over a small neighborhood of $\f_1$ in $Q_\alpha$.
(As above, we assume that the K{\oe}nigs function has been
normalized so that $\kappa_\f(c_\f^{\;1})=1\,$.)
In fact given $\f_1$, and given these $\w$ coordinate values,
we will show how to
reconstruct the map $~\f:\U_\f\to \U_\f~$ up to conformal conjugacy.

 Choose a sequence of connected open sets
$$N_\f(0)~\subset~N_\f(1)~\subset~N_\f(2)~\subset~\cdots~~\subset~\U_\f $$
with union equal to the entire space $~\U_\f~$ as follows. 

\begin{definition}\label{d-nhds}
Let $N_\f(\ell)$ be the connected component which contains the fixed point 
$p_\f$ in the open set
$$ \{~z\in \U_\f~;~|\kappa_\f(z)|<a/\lambda_\f^{~\ell}~\}~.$$
Here the constant $a$ should be small enough so that $\f$ has no critical points
in the closure
$\overline N_\f(1)$, and should be carefully chosen so that no critical orbit
of $\f$ hits the boundary $\partial N_\f(1)$.
\end{definition}

Here are three easily verified properties.

\begin{enumerate}
\item[\bf(a)] If this condition is satisfied for some given map $\f_1$, then it
 will also be satisfied for any map $\f$ in a sufficiently small
 neighborhood of $\f_1$.\smallskip

\item[{\bf(b)}] Each such $\f$ maps $N_\f(1)$ biholomorphically
onto the proper subset $N_\f(0)$ of itself.\smallskip

\item[{\bf(c)}] Each $N_\f(\ell+1)$ maps onto $N_\f(\ell)$ by a branched 
covering which is branched only over those critical values of $\f$
 which lie in $N_\f(\ell)$. The  topological 
pattern of this branching remains the same for all $\f$ close to $\f_1$.
\end{enumerate}\smallskip

\noindent(In this special case where ${\mathcal U}_\f$ is connected,
 the $N_\f(\ell)$ are all connected sets; but this will
no longer be true in the general case considered later. But in all cases,
each connected component of $N_\f(\ell)$ will be simply connected, with smooth
boundary.)\smallskip

Given $\f_1\in Q_\alpha$, we must show that any $\f$ sufficiently close to $\f_1$
is uniquely determined by $\lambda_\f$, together with the numbers
 $\kappa_\f(c_\f^{\;j})$. We will first show by induction on $\ell$ that
the conformal conjugacy class of the restriction
\begin{equation}\label{e-fN}
\f:N_\f(\ell)\to N_\f(\ell)
\end{equation} is uniquely
determined by this data. To begin the induction,  we need only
 $\lambda_\f$ to determine the conformal conjugacy class of $\f$ restricted
to $N_\f(1)$. Assuming that we have constructed $N_\f(\ell)$ and the restriction
of $\f$ to this Riemann surface, we need only to know the precise branch points
in order to construct a Riemann surface isomorphic to $N_\f(\ell+1)$ as
 a branched covering. But
locally, for $\f$ near $\f_1$, each branch point $\f(c_\f^{\;j})$ is uniquely
determined by the K{\oe}nigs coordinate $\kappa_\f(c_\f^{\;j})$. The inclusion
map of $N_\f(\ell)$ into $N_\f(\ell+1)$ is then determined inductively.
In fact the required branched covering $N_\f(\ell+1)\to N_\f(\ell)$ can be
constructed as an extension of the branched covering $N_\f(\ell)\to 
N_\f(\ell-1)$, which is known by the induction hypothesis.

Now passing to the union as $\ell\to\infty$, we conclude that the conformal
conjugacy class of $\f:\U_\f\to \U_\f$ is uniquely determined. Passing to the
Blaschke product model, this means that the associated point in $\B^{\w+1}$
is uniquely determined, up to a choice of boundary markings. But the boundary
marking must vary smoothly with $\f$, hence it is uniquely determined by
the boundary marking for $\f_1$. Finally, using Theorem~\ref{t-top-mod},
it follows that $\f$ is uniquely determined. Since a holomorphic map
 which is one-to-one must be biholomorphic, this completes the proof for the
case that the full mapping scheme for $H_\alpha$ has only one vertex.
\smallskip

The proof for an arbitrary connected mapping scheme $S$ is similar. Again
let $\,\U_\f$ be the union of the Fatou components which contain critical
or postcritical points. Let $m$ be the period of the unique attracting orbit,
and let $\lambda_\f$ be its multiplier. Choosing some $m$-th root
 $\lambda_\f^{\,1/m}$, the modified K{\oe}nigs equation
$$\kappa_\f\big(\f(z)\big)~ =~\lambda_\f^{1/m}\kappa_\f(z) $$
has a solution $\kappa_\f:\U_\f\to\C$
which is unique up to a multiplicative constant. As before
we can normalize so that $\kappa(c_\f^{\;1})=+1$. Constructing open sets
$$N_\f(0)\subset N_\f(1)\subset\cdots$$
 with union $\U_\f$ as before,
we can again prove inductively that the conformal conjugacy class of $\f$
restricted to each $N_\f(\ell)$ is determined by the $w$ coordinates
$$~\lambda_\f^{1/m}\,,~\kappa_\f(c_\f^{\;2})\,,~\ldots\,,~\kappa_\f(c_\f^{\;\w})\,,$$
and therefore that $\f$ is uniquely determined. The result in the case of a
scheme $S$ with several components then follows easily, applying this argument
to one component of $S$ at a time. This completes the proof of Lemma
 \ref{l-candiff}.\end{proof}

\begin{proof}[ Proof of Theorem~\ref{t-5.1}]
It now follows that the diffeomorphism $H_\alpha \to H_\beta$ is
holomorphic everywhere. In fact the Cauchy-Riemann equations, which are
necessary and sufficient conditions for a $C^1$-smooth map to be holomorphic,
are satisfied throughout a dense open subset of $H_\alpha$. Hence
Theorem \ref{t-5.1} follows by continuity.
\end{proof}

\begin{rem} \label{r-hol-fam}
In fact, 
  the entire holomorphic dynamics of the family
of maps $~\f_{\U_\f}:\U_\f\to\U_\f~$ with $\f\in H$ depends only on the mapping
scheme $S$. More precisely, let $~\U_H\subset H\times|S_0|\times\C~$ be the
 set of triples $~(\f,~ s,~z)~$ with $\f\in H$  and $(s,~z)\in\U_\f$. Then
 the biholomorphic conjugacy class of the dynamical system \hbox{$(\f,~s,~z)
\mapsto \big(\f,~\f(s,~z)\big)~$} depends only on the mapping scheme of $H$. 
The proof is  similar to the proof of Theorem \ref{t-5.1}.
\end{rem}\smallskip

\begin{rem} The above discussion doesn't discuss boundary  behavior.  In fact,
the diffeomorphism of Theorem \ref{t-5.1} cannot always extend
continuously over the boundary. (Compare Example \ref{ex-sc}.) 
 \smallskip

Here is an even deeper question.
In the Douady-Hubbard theory of the Mandelbrot set $M$, every hyperbolic
component $H\subset M$ embeds in a small copy $M_H\subset M$, where $M_H$
is homeomorphic to $M$ under a homeomorphism which carries $H\subset M_H$
to the cardioid component $H_0\subset M$. More generally, we can ask
the following.
\smallskip

\begin{quote}\it Under what conditions does the canonical biholomorphic map
 from \hbox{$H_0^S\subset\p^S$} to a given hyperbolic component
$H\subset\p^{S_0}$ extend to an embedding of the entire connectedness
locus $\cl(\p^S)$ into $\cl(\p^{S_0})$?
\end{quote}
\smallskip

\noindent As a simplest example, let $H$ be a hyperbolic component
 of type $D$ in the cubic connectedness locus $\p^3$.
 (Compare Figure~\ref{f-wt2}D.) When is $H$ contained in a 
complete 
Cartesian product $M\times M$ of two copies of the Mandelbrot set?
\end{rem}

\section{Real Forms.}\label{s6}

First consider a real polynomial map $f_\R:\R\to\R$, of degree $d\ge 2$. 
We can extend $f_\R$  uniquely to a complex polynomial map $f:\C\to\C$.
 This extended map  will 
 commute with the complex conjugation operation 
$z\mapsto\overline z$. 
If $f$ is hyperbolic, 
and if $~ {\mathcal U}_f~$ is the union of those Fatou components which contain
 critical points, then complex conjugation
 carries this set ${\mathcal U}_f$ onto itself, carrying each component which
intersects the real axis onto itself, but interchanging the remaining
components in pairs.
In order to find an appropriate universal
model for such behavior, we consider the following construction.\medskip

Given a scheme $S_0$, consider 
 antilinear\footnote{A function $\gamma$ between complex vector
spaces is called \textbf{\textit{~antilinear~}} (or ``conjugate-linear'')
if $~\gamma(c\,v)=\overline c\,\,\gamma(v)~$ for every $c\in\C$.}
 involutions  {$~\gamma:|S_0|\times\C\to|S_0|\times\C~$}
 which commute with the base map
\begin{equation}\label{e-basem}
 \f_0(s,\,z)~=~(F(s),\,z^{d(s)})~.
\end{equation}
In other words, for each $s\in|S_0|$ we assume that
 $\gamma$ is antilinear
 as a function from $s\times\C$ to some $s'\times\C\,$, and further
we assume that
$$\gamma \circ\gamma={\rm identity}\,,\qquad{\rm and}
\qquad \gamma\circ \f_0=\f_0\circ\gamma\,.$$

\smallskip

\begin{definition}\label{d-rf}
For any $\gamma$ satisfying these conditions, the subset
$\p^{S_0}(\gamma)$ consisting of all $\f\in\p^{S_0}$ which commute with
$\gamma$ will be called a \textbf{\textit{~real form~}} of $\p^{S_0}$.
Two such real forms $\p^{S_0}(\gamma_1)$ and $\p^{S_0}(\gamma_2)$
will be called \textbf{\textit{~isomorphic~}}
if there is an automorphism $\eta$ of $|S_0|\times\C$ which is linear on
each $s\times\C$
such that the subset $~\p^{S_0}(\gamma_1)\subset\p^{S_0}~$ maps isomorphically
 onto $~\p^{S_0}(\gamma_2)~$ under the correspondence
 $~\f\,\mapsto\,\eta^{-1}\circ \f\circ\eta\,$.
As an example, if  $~\gamma_2=\eta^{-1}\circ\gamma_1\circ\eta$,~ then the
corresponding real forms will certainly be isomorphic.
\end{definition}

The following two lemmas will help to clarify this definition.

\begin{lem}\label{l-gam}
Any antilinear involution $\gamma$ of $|S_0|\times\C$
which commutes with $\f_0$ is given by the formula
$$ \gamma(s,\,z)~=~\big(s'\,,~\aac(s)\,\overline z\big)\,,$$
where each $~\aac(s)~$ is a root of unity satisfying
\begin{equation}\label{e-alphaF}
\aac\big(F(s)\big)~=~\aac(s)^{d(s)}\,,
\end{equation}
and where $~s\leftrightarrow s'~$
is an involution of $|S_0|$ $($or the identity map$)$ satisfying
$$ d(s)=d(s')\,,\quad \aac(s)=\aac(s')
\,,\quad{\it and}\quad F(s')=F(s)'\,.$$
 For each scheme $~S_0$, there are only finitely many such involutions
 $~\gamma\,$.
\end{lem}

\begin{proof} The
 fact that the $\aac(s)$ are roots of unity depends on following
the relation (\ref{e-alphaF}) around each cycle contained in $S_0$.
Further details are straightforward, and will be left to the reader.
\end{proof}

We have described $\p^{S_0}$ as a complex affine space. However, since
it has a preferred base point $\f_0$, there is a closely related complex
vector space consisting of all differences $\f-\f_0$.

\begin{lem}\label{l-fps}
As in Lemma \ref{l-gam}, let $\gamma$ be an  antilinear involution 
 of $|S_0|\times\C$  which commutes with $\f_0$. Then $\gamma$ acts on 
$\p^{S_0}$ by  an involution
\begin{equation}\label{e-conj-invo}
 \f~\mapsto~\gamma\circ \f\circ\gamma\,,
\end{equation}
which acts antilinearly on the vector spaces of differences $\f-\f_0$.
The fixed point set $~\p^{S_0}(\gamma)\subset\p^{S_0}~$
of the involution $(\ref{e-conj-invo})$
is a real affine space with real dimension $w(S_0)$.
\end{lem}

Thus $\p^{S_0}(\gamma)$ is a real affine space whose elements are complex maps.
The real dimension of $\p^{S_0}(\gamma)$ is equal to the complex dimension
 of $\p^{S_0}$,  or in other words to half the real dimension of $\p^{S_0}$.
\medskip

\begin{rem}\label{r-caution} Two different
antilinear involutions $\gamma$ and $\gamma'$
of $|S_0|\times\C$ may give rise to the same antilinear
involution (\ref{e-conj-invo}) of $\p^{S_0}$,
 and hence to the same real form. This occurs if and only if
 $\gamma\circ\gamma'$ belongs to the subgroup $\G_0(S_0)$ of linear 
automorphisms which commute with all elements of $\p^{S_0}$, as described
in Lemma~\ref{l-G_0}.\end{rem}

\begin{proof}[\sc Proof of Lemma~\ref{l-fps}.] If
$$ \f(s,\,z)~=~\Big(F(s)\,,~z^d+\sum_0^{d-2} c_jz^j\Big)\,,$$
where $d=d(s)$, then a brief computation using Lemma~\ref{l-gam} shows that
$$ \gamma\circ \f\circ\gamma(s,\,z)~=~\Big(F(s)\,,~z^d+
\sum_0^{d-2} \aac\big(F(s)\big)\,\overline\aac(s)^j\,\overline c_j\, z^j\Big)\,.$$
Evidently the correspondence
\hbox{$~c_j\mapsto \aac\big(F(s)\big)\,\overline\aac(s)^j\,\overline c_j~$}
between coefficients is antilinear, as asserted.

Now note that for any antilinear
involution of a complex vector space $V$,
 the fixed point set is a real vector space with real
 dimension equal to exactly half of the complex dimension of $~V$.
In fact the underlying real vector space of $V$ splits as the direct sum of
the $(+1)$-eigenspace, which is precisely the fixed point set, and the
 $(-1)$-eigenspace.
But antilinearity guarantees that multiplication by $\sqrt{-1}$ must
interchange the $+1$ and $-1$ eigenspaces, and the conclusion follows.
\end{proof}

If we generalize our spaces of
polynomial maps by allowing not only monic polynomials with leading
coefficient $+1$, but also polynomials with leading coefficient $-1$,
then there is an alternative standard model for real forms which is perhaps
 easier to work with, since it eliminates the complex roots of unity
of Lemma \ref{l-gam}.

\begin{definition}\label{d-signed}
Given any choice of signs $\bsig:|S_0|\to\{\pm 1\}$, let $\p_\bsig^{S_0}$ be
the complex affine space consisting of all maps
 $\f:|S_0|\times\C\to|S_0|\times\C$ which carry
each $s\times\C$ to $F(s)\times\C$ by a centered polynomial of degree $d(s)$
with leading coefficient $\bsig(s)$.
\end{definition}

\begin{definition}\label{l-std-inv}
To any automorphism $\bio$ of $S_0$ which satisfies $~\bio\circ\bio=
{\rm identity}\,$,\break so that  $\bio$ is either 
 the identity map or an involution,
 there is associated a\break \textbf{\textit{~standard antilinear
involution~}}  $\gamma_\bio:|S_0|\times\C\to|S_0|\times\C$, given by
$$ \gamma_\bio(s,\,z)~=~\big(\bio(s),\,\overline z\big)\,.$$
Combining this with the previous definition, we can form the
real affine space $\p_\bsig^{S_0}(\gamma_\bio)$ consisting of
all $\f\in\p_\bsig^{S_0}$ with $$\f\circ\gamma_\bio~=~\gamma_\bio\circ \f\,.$$
\end{definition}

\begin{lem}\label{l-alt-form}
Let $\p^{S_0}(\gamma)$ be a real form of $\p^{S_0}$, and let
$\bio:s\leftrightarrow s'$ be the automorphism of $S_0$ associated
with $\gamma$. Then there exists at least one
 choice of signs $\bsig:|S_0|\to\{\pm 1\}$
so that the real form  $\p_\bsig^{S_0}(\gamma_\bio)\subset\p_\bsig^{S_0}$
is isomorphic to the real form $\p^{S_0}(\gamma)\subset\p^{S_0}$.
\end{lem}

\medskip

\begin{proof}
We will construct an automorphism
$\eta$ of $|S_0|\times\C$ of the form
$$\eta(s,\,z)~=~\big(s,\,\bb(s)\,z\big)$$
and an associated choice of signs $\bsig$
so that the conjugation $\f\mapsto\eta^{-1}\circ \f\circ\eta\,$ maps  
$\p^{S_0}(\gamma)$
isomorphically onto $\p_\bsig^{S_0}(\gamma_\bio)\subset\p_\bsig^{S_0}$.
Start with the expression
$$ \gamma(s,\,z)~=~\big(s',~\aac(s)\,\overline z\big) $$
of Lemma~\ref{l-gam}, and choose each $\bb(s)$ so that
$~\bb(s)^2\,=\,\aac(s)\,$ with $~\bb(s') = \bb(s)$.\, Then a brief computation,
 using the fact that $|\bb(s)|=1$, shows that
$$ \eta^{-1}\circ\gamma\circ\eta~=~\gamma_\bio\,.$$
Now define $\bsig(s)$ by the equation
\begin{equation}\label{e-sig}
\bsig(s)\, \bb\big(F(s)\big)~=~\,\bb(s)^{d(s)}\,.
\end{equation}
Squaring this equation (\ref{e-sig}), and using the identity
 $~\aac\big(F(s)\big)=\aac(s)^{d(s)}$, we see that $\bsig(s)=\pm 1$.
For the base map of equation (\ref{e-basem}),
a straightforward computation  shows that 
$$ \eta^{-1}\circ \f_0\circ\eta(s,\,z)~=~\big(F(s)\,,~ \bb(F(s))^{-1}
\bb(s)^{d(s)}z^{d(s)}\big)~=~\big(F(s),\,\bsig(s)\,z^{d(s)}\big)\,.$$
A similar argument shows that the correspondence
 $~\f\mapsto\eta^{-1}\circ \f\circ \eta~$ carries $\p^{S_0}$ isomorphically onto
$\p^{S_0}_\bsig$; and hence maps $\p^{S_0}(\gamma)$ isomorphically onto
$\p^{S_0}_\sigma(\gamma_\bio)$.\end{proof}
\smallskip

\begin{figure}[!ht]
\centering
\subfigure[\label{f3a}]{\psfig{figure=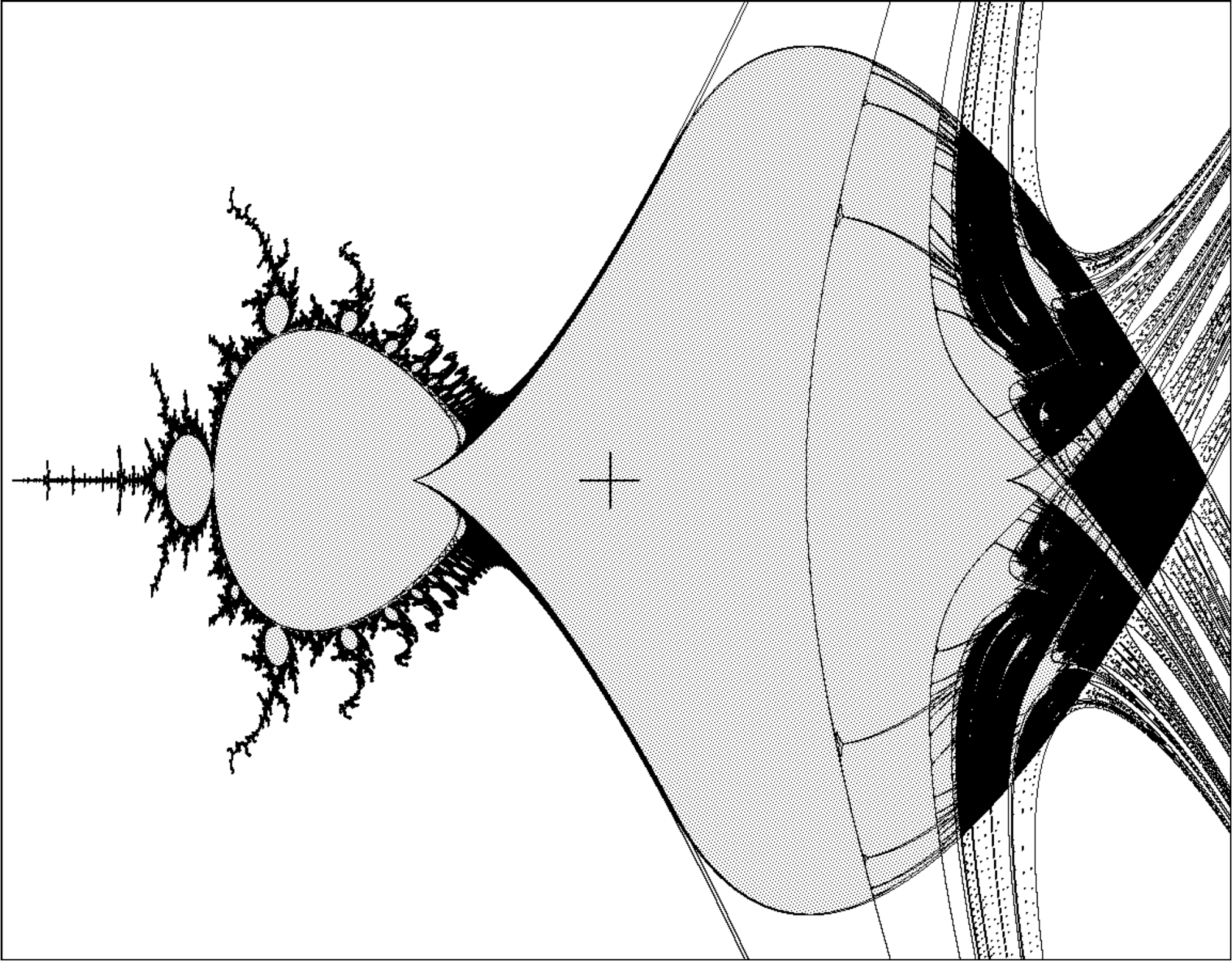,height=2.2in}}
\vspace{.2cm}

\subfigure[\label{f3b}]{\psfig{figure=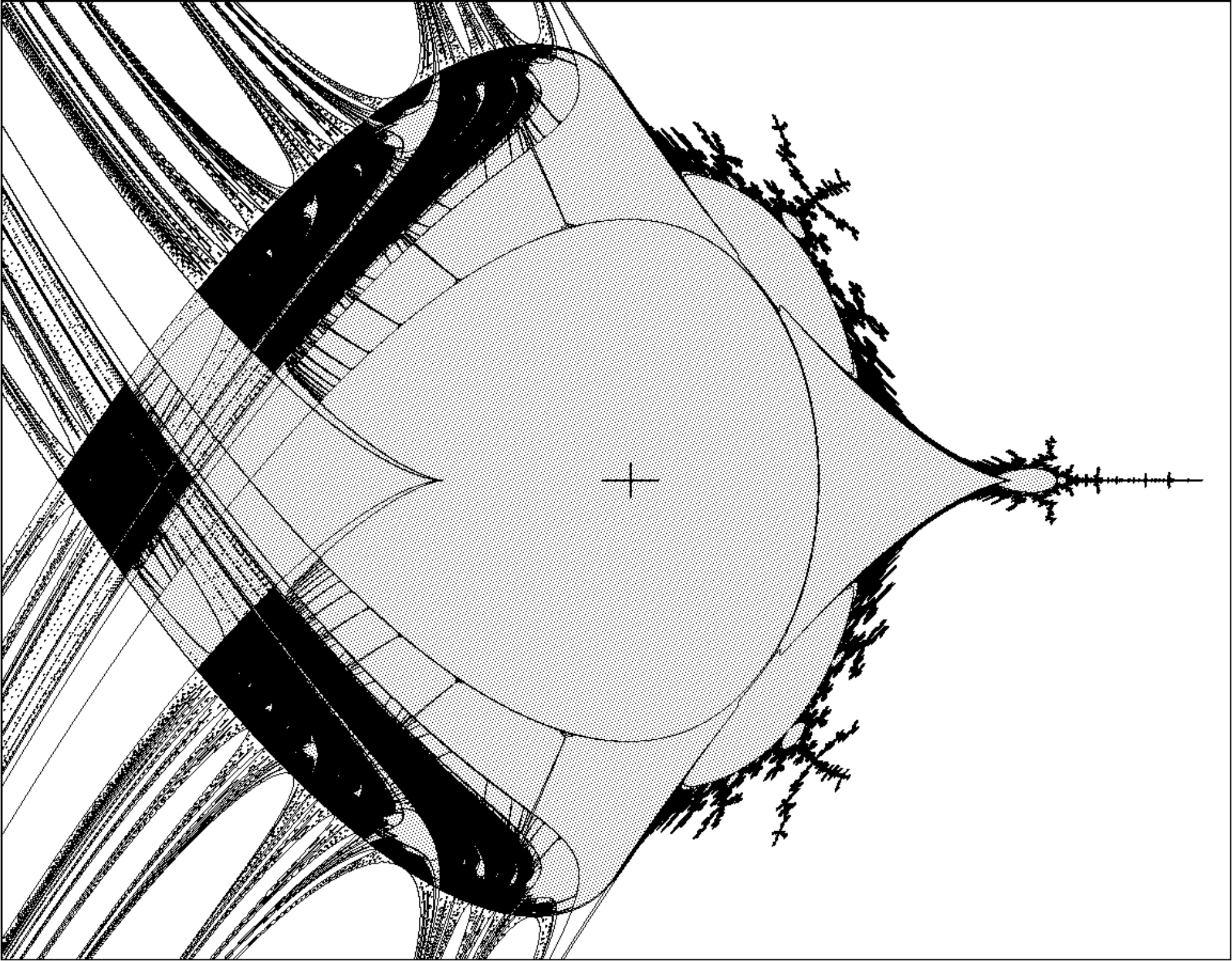,height=2.2in}}
\caption{\it \label{f3} The spaces $\p_+^3(\gamma_{\bio_0})$ and
 $\p_-^3(\gamma_{\bio_0})$ of real
cubic maps, intersected with the complex connectedness locus $\cl^{3}$.
More precisely, these pictures show the real
 $(A,b)$-plane where \hbox{$x\mapsto \pm x^3-3Ax+b$.}
$($Compare {\rm\cite{M1}}.$)$ In both figures there is a compact real
connectedness locus, containing Mandelbrot-like subsets in the left
(respectively right) half-plane
and also a Cantor set's worth of curves reaching off to infinity in the other
half-plane. These curves represent maps
for which just one of the two critical points has bounded orbit.}
\medskip
\end{figure}

\begin{rem}\label{r-sgn-ch}
The signs $\bsig(s)$ are far from uniquely determined, since we are free to
switch the signs of the $\bb(s)$. Examining the defining equation
(\ref{e-sig}), we see that replacing any given $\bb(s_0)$ by $-\bb(s_0)$
will have the following effect:

\begin{itemize}
\item If $F(s)=s_0$ with $s\ne s_0$, then $\bsig(s)$ will change sign.

\item If $F(s_0)\ne s_0$ and $d(s_0)$ is odd,  or if
$F(s_0)=s_0$ and $d(s_0)$ is even, then $\bsig(s_0)$ will change sign.
\end{itemize}
In all other cases, the signs remain unchanged. Here is an easy example.

\begin{lem}\label{l-easy}
If $S_0$ is a union of cycles, so that $F$ maps $S_0$ bijectively onto itself,
and if all the degrees $d(s)$ are even, then we can choose $\bsig(s)$ to be
identically $+1$, so that every real form $\p^{S_0}(\gamma)$ is isomorphic to
the associated  $\p^{S_0}(\gamma_\bio)$.\end{lem}

The proof is straightforward. \qed

\end{rem}

As another example,
consider the space $\p^d$ of monic centered polynomial maps
of degree $d$.

\begin{lem}\label{l-rf}
If the degree $d$ is even, then every real form of $\p^d$ is isomorphic to the
standard real form consisting of monic centered polynomials with real
coefficients. But if $d\ge 3$ is odd then there are two equivalence classes,
represented by real polynomials either with leading coefficient $+1$ or with
leading coefficient $-1$.\end{lem}

\begin{proof} The associated scheme $S$ has only one point, so the
 correspondence $s\leftrightarrow s'$ must be the identity map $\bio_0$.
 If $d$ is even, then it follows from Lemma \ref{l-easy} that there is only
 one real form $\p^d_+(\gamma_{\bio_o})$ consisting of monic centered
 polynomials with real coefficients. However, if $d$ is odd there is 
an additional real form $\p^d_-(\gamma_{\bio_0})$,
 consisting of centered  polynomials with leading coefficient $-1$
and with real coefficients. To see that these two are not isomorphic,
note that each can be considered as a family of maps from $\R$ to $\R$.
In either case, we can compactify $\R$ by adding  points at 
 $+\infty$ and  $-\infty$, and then extend
 to a self-map of $\R\cup\{+\infty\}\cup\{-\infty\}$. In the 
 $\p^d_+(\gamma_{\bio_o})$ case we obtain two fixed points at infinity;
 but in the $\p^d_-(\gamma_{\bio_o})$ case we obtain a
 period two orbit at infinity.\end{proof}

 Compare Figure \ref{f3} for pictures of the connectedness locus
in  the two real forms of $\p^3$, denoted by
$\,\p_+^3(\gamma_{\bio_0})$ and 
$\p_-^3(\gamma_{\bio_0})$.
These are 2-dimensional pictures, since the total weight is $\w=2$.
\medskip

\begin{figure}[!ht]
\centering
\subfigure[\label{f4a}]{\psfig{figure=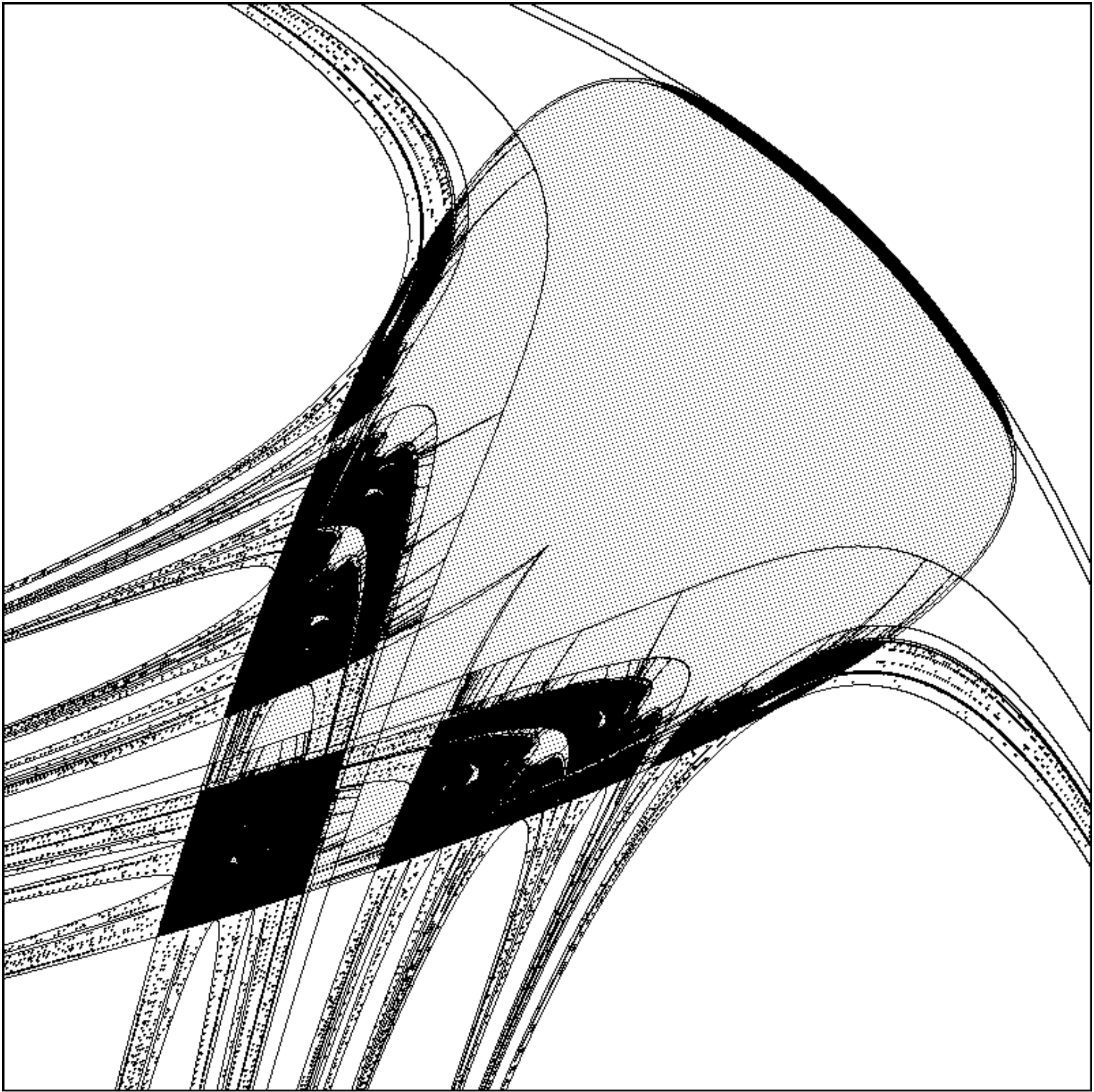,height=2in}}\qquad
\subfigure[\label{f4b}]{\psfig{figure=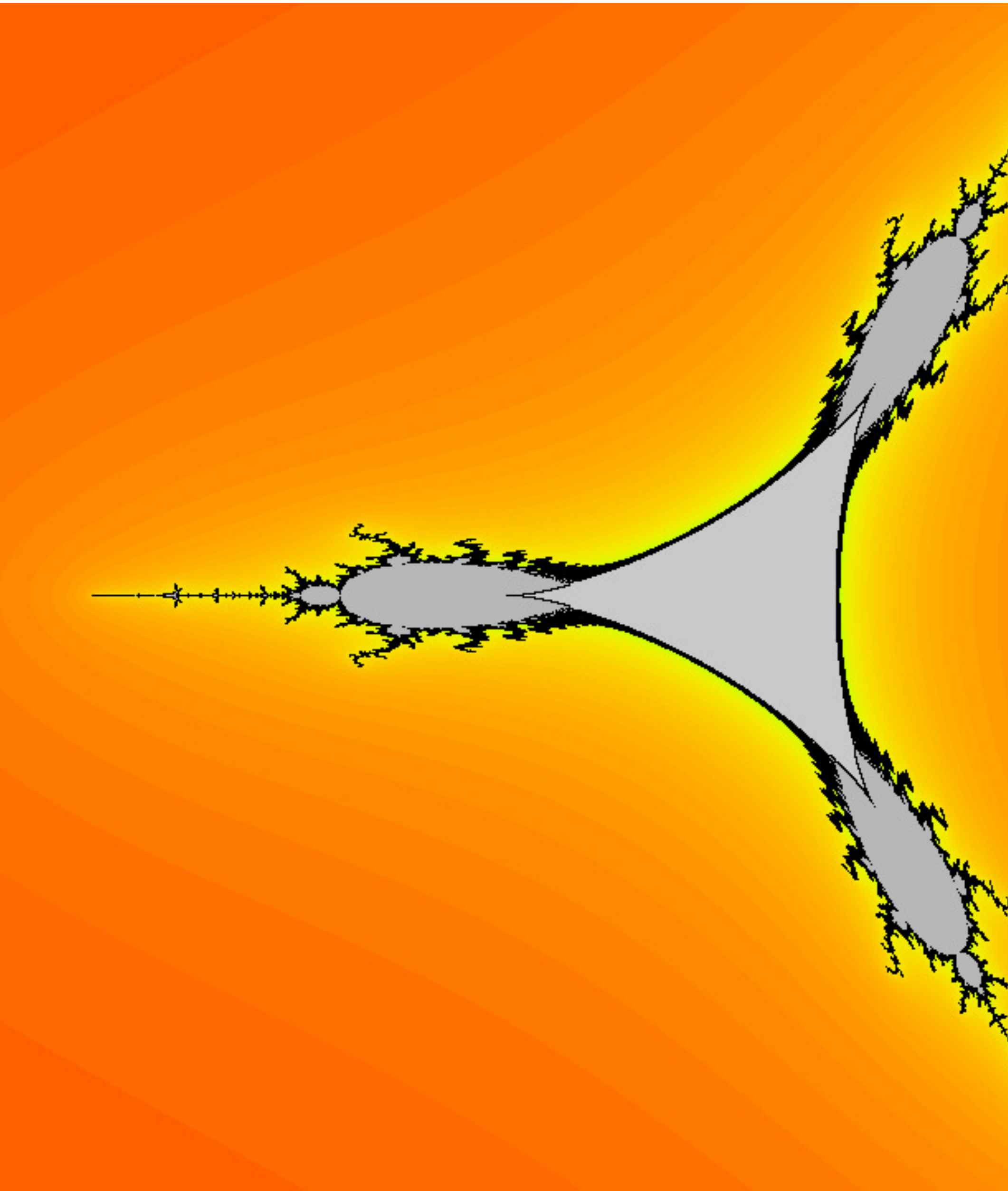,height=2in}} 
\vspace{.2cm}

\subfigure[\label{f4c}]{\psfig{figure=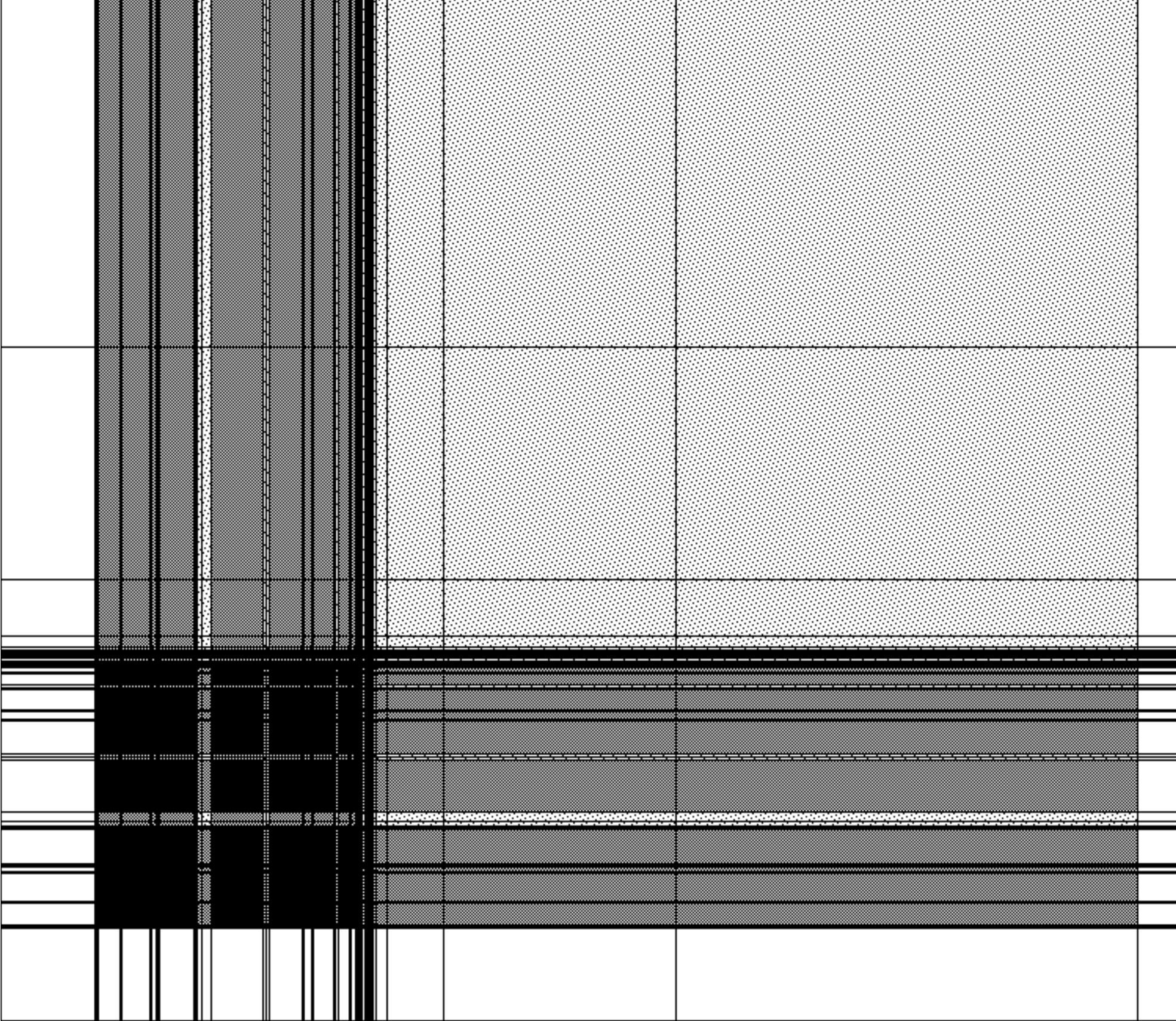,height=2in}}\qquad
\subfigure[\label{f4d}]{\psfig{figure=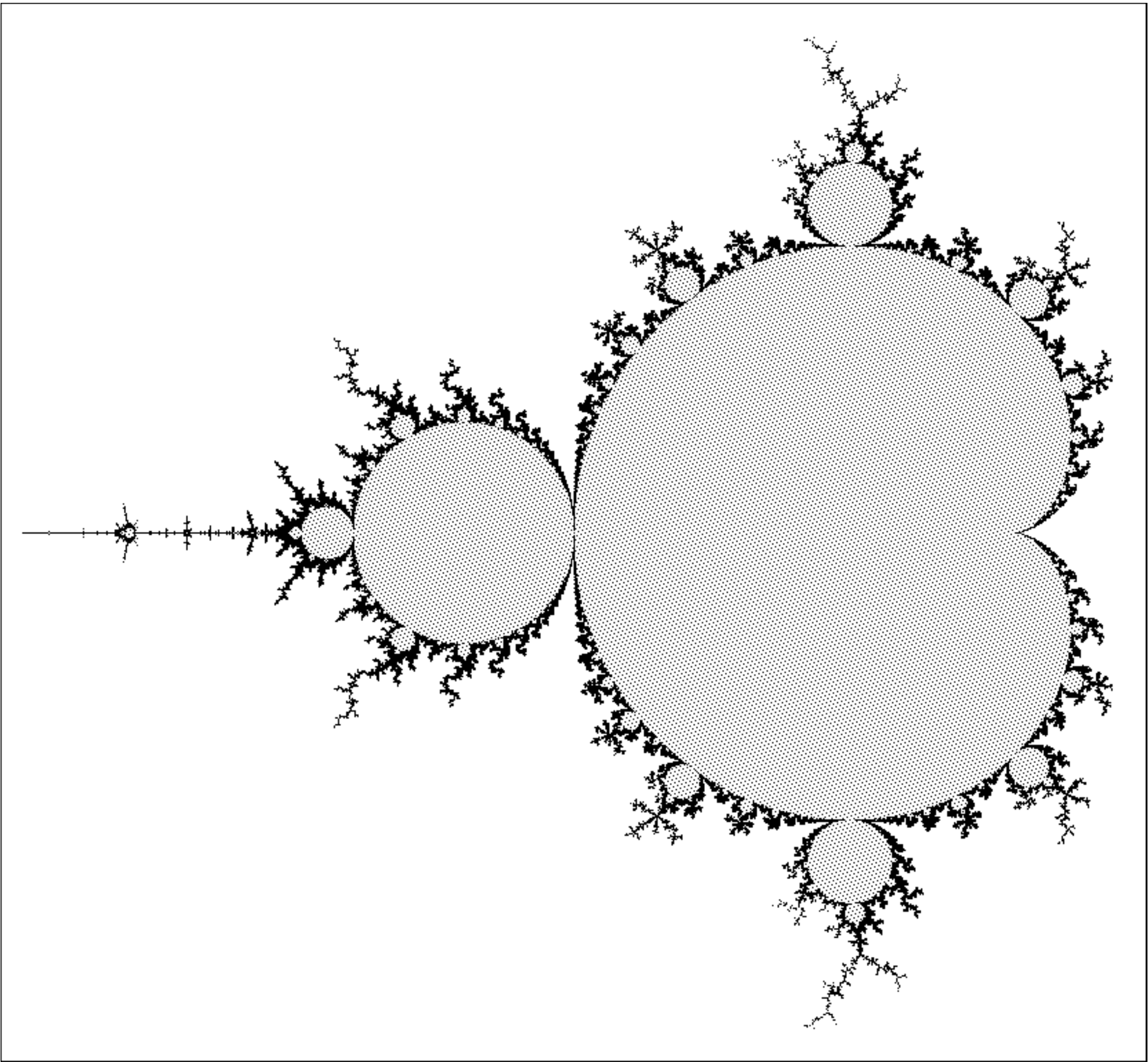,height=2in}} 
\caption{\it \label{f4} Connectedness loci for four real forms of weight 
two.}
\end{figure}

\begin{ex}[{\sc Other schemes of weight two}]\label{ex-w=2} Let $S_{B}$
 be the \textbf{\textit{~bitransitive~}}
 scheme, as illustrated in Figure \ref{f-wt2}B.
Thus the two points $s_1$ and $s_2$ of $\vert S_{B} \vert$ are mapped to
 each other:
 $F:s_1\leftrightarrow s_2$. It follows from Lemma \ref{l-easy} that there are
just two real forms $\p^{S_B}(\gamma_\bio)$ corresponding to the two possible
choices of $\bio:s\leftrightarrow s'$. In either case, every periodic orbit
must have even period.
 \medskip

{\bf The Top.\footnote{``Top'' in the sense of a children's spinning toy
(trompo, toupie, Kreisel). (This figure might also be interpreted as an
alien face, with pointy chin and shiny eyes.)}}
 Choosing the identity map $\bio_0(s)=s$,
we obtain the antilinear map $\gamma_{\bio_0}$ with
\begin{equation}\label{e-indep}
 (s_1,\,z)\leftrightarrow(s_1,\,\overline z)\,,\qquad
(s_2,\,z)\leftrightarrow(s_2,\,\overline z)\,.
\end{equation}
 The corresponding real parameter space $\p^{S_{B}}(\gamma_{\bio_0})$ can
 be described as the set of all maps of the form
$$ \f(s_1,\,z)~=~(s_2,\, z^2+c_1)\,\qquad \f(s_2,\,w)~=~(s_1,\,w^2+c_2)$$
where $c_1,\,c_2$ are real parameters. Note that the 2-fold iterate
has the form
$$ \f^{\circ 2}(s_1,\,z)~=~\big(s_1,\,z^4+2c_1z^2+(c_1^2+c_2)\big)\,,$$
 where the expression on the right varies precisely over all even real 
 polynomials of degree 4. (Compare \cite{Ra}.)
The connectedness locus in the $(c_1,\,c_2)$-plane is shaped like a
children's top, as shown in Figure~\ref{f4a}. 
 For each $\f$ in the
connectedness locus, the union of all real bounded orbits consists of
a non-trivial closed interval in each $s_j\times\R$,
$$K_\R(\f)~=~\big(s_1\times[-a_1,\,a_1]\big)~\cup
~\big( s_2\times[-a_2,\,a_2]\big)~,$$
where the right hand endpoints form a periodic orbit
 $~(s_1,\,a_1)~\leftrightarrow~(s_2,\,a_2)~$ with
multiplier $\lambda\ge 1$. The  boundary of the connectedness locus consists
 of three real analytic pieces, which are defined respectively by the
conditions that $c_1=-a_2$, or $c_2=-a_1$, or that $\lambda=1$.
 (In the first two cases, one critical orbit is preperiodic, with
 $\f^{\circ 2}(s_j,\,0)=(s_j, a_j)$), while in the third case
both critical orbits converge to the same parabolic orbit.) At the
common endpoint of any two of these three
real analytic pieces, two of these three conditions are satisfied.

 As in Figure \ref{f3}, the parameter picture also shows uncountably
many curves reaching off to infinity, representing maps for which just one of
the two critical points has bounded orbit. \medskip

{\bf The Tricorn.} The other real form of $\p^{S_{B}}$
corresponds to the non-trivial involution
 $\bio_1:s_1\leftrightarrow s_2\,$,~ with
\begin{equation}\label{e-flip}
\gamma_{\bio_1}\,:\,(s_1,\,z)~\longleftrightarrow~(s_2,\,\overline z)\,.
\end{equation}
The form $\p^{S_{B}}(\gamma_{\bio_1})$ then consists of all maps of the form
\begin{equation}\label{e-tric}
 \f(s_1,\,z)~=~(s_2,\, z^2+c)\,,\qquad \f(s_2,\,w)~=~(s_1,\,w^2+\overline c)
\end{equation}
where the parameter $c$ is complex.
The corresponding  connectedness locus, 
as shown in Figure \ref{f4b}, is known as the \textbf{\,\textit{~tricorn\,}}
(or \textbf{\textit{\,Mandelbar set\,}}).  (Compare \cite{CHRS}, \cite{M1},
 \cite{NS2}.)
The tricorn is invariant under 120$^\circ$ rotation. To prove this
note that the equation (\ref{e-tric}) remains valid if $z,\,w,\,c$
are replaced by $~\eta\,z,~\overline\eta\,w~$ and $~\overline\eta\,c~$
respectively, where $\eta^3=1$. The
central hyperbolic component of the tricorn, to be denoted by
${\mathcal H}_{\rm tric}$, is bounded by a deltoid
curve. In fact the closure $\overline{\mathcal H}_{\rm tric}$ can be
parametrized as
$$ c(t)~=~ t/2~-~\overline t^2/4~,$$
where $t$ varies over the closed unit disk $\overline\D$. (Note that
$c(\eta\,t)=\eta\,c(t)$ when $\eta^3=1$.) There are cusps at the three
points where $t^3=-1$. For each $t\in\overline\D$ with $t^3\ne -1$,
 the Julia set consists
of two simple closed curves, each mapped to the other with degree two.
Thus we can parametrize one of these curves by the circle
 $|z|=1$ so that the second iterate maps $z$ to $z^4$, It follows that there are
three period two orbits in the Julia set, corresponding to the three cube roots
of unity. 
 For $|t|<1$ there is also an attracting period two orbit 
 $$(s_1,\, \overline t/2)~~\leftrightarrow~~(s_2, t/2)  $$
in the Fatou set, with multiplier
$t\,\overline t\in[0,\,1)$. As $|t|\to 1$, this attracting orbit tends to
one of the repelling period two orbits in the Julia set, and these
 become parabolic in
the limit. The three edges of $\partial{\mathcal H}_{\rm tric}$ correspond to
the three cube roots of unity in the discussion above. 
However, at the three cusp points, two of the three boundary period two
orbits crash together
so that each component of the Julia set
 becomes a copy of the ``fat basilica'' Julia set $~J(z\mapsto z^2-z)$. 

\begin{figure} [ht]
\centerline{\psfig{figure=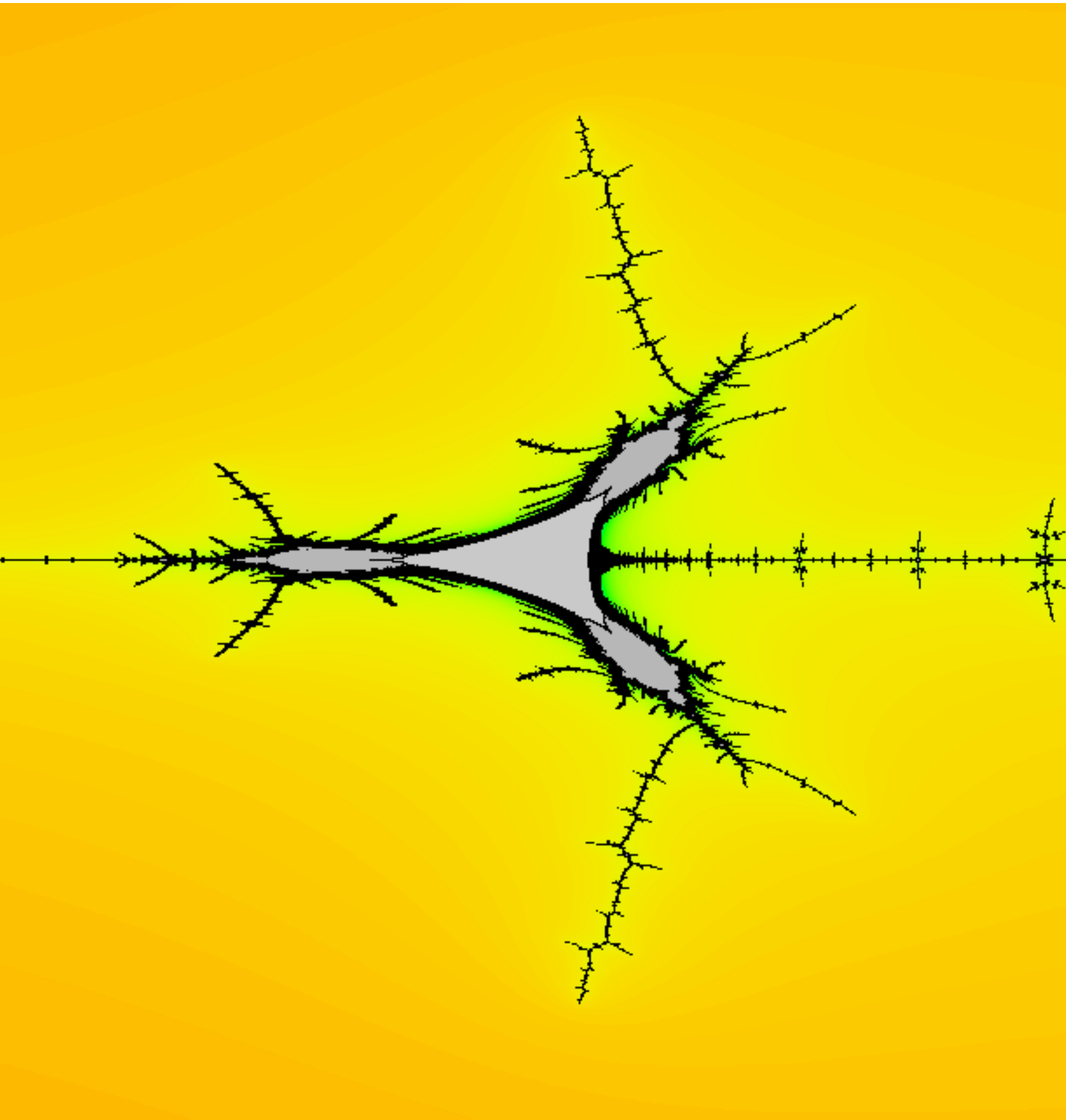,height=2.5in}}
\caption{\label{f-stric}\it
 A small copy of the tricorn contained in the tricorn,
 centered at $c=-1.7548\cdots$ along the real axis.}
\end{figure}
 
At each of the three cusps of $\overline{\mathcal H}_{\rm tric}$,
 there is an attachment which
resembles a distorted copy of the $1/2$-limb of the Mandelbrot set.
In fact the intersection of the left hand attachment
with the real axis is identical to the real part of the $1/2$-limb.
However, the distortion is so extreme that the tricorn
 is not locally connected. (See \cite{NS1}, \cite{HS}.)
Furthermore, the resemblance is not perfect. In particular, for any
primitive small copy of the Mandelbrot set with odd period within the
Mandelbrot $1/2$-limb there is a corresponding small copy\footnote{Here the
word ``copy'' is used loosely, and does not necessarily mean homeomorphic
copy.}
 of the tricorn
within the tricorn. (Compare Figure \ref{f-stric}.)

It is interesting to note that both Figures \ref{f3a} and \ref{f3b}
also contain small distorted 
copies of
the tricorn, which would be visible only under high magnification. (Compare
 \cite{M1}, \cite{NS2}, and see Figure \ref{f-cubtric}.)

\begin{figure} [ht]
\centerline{\psfig{figure=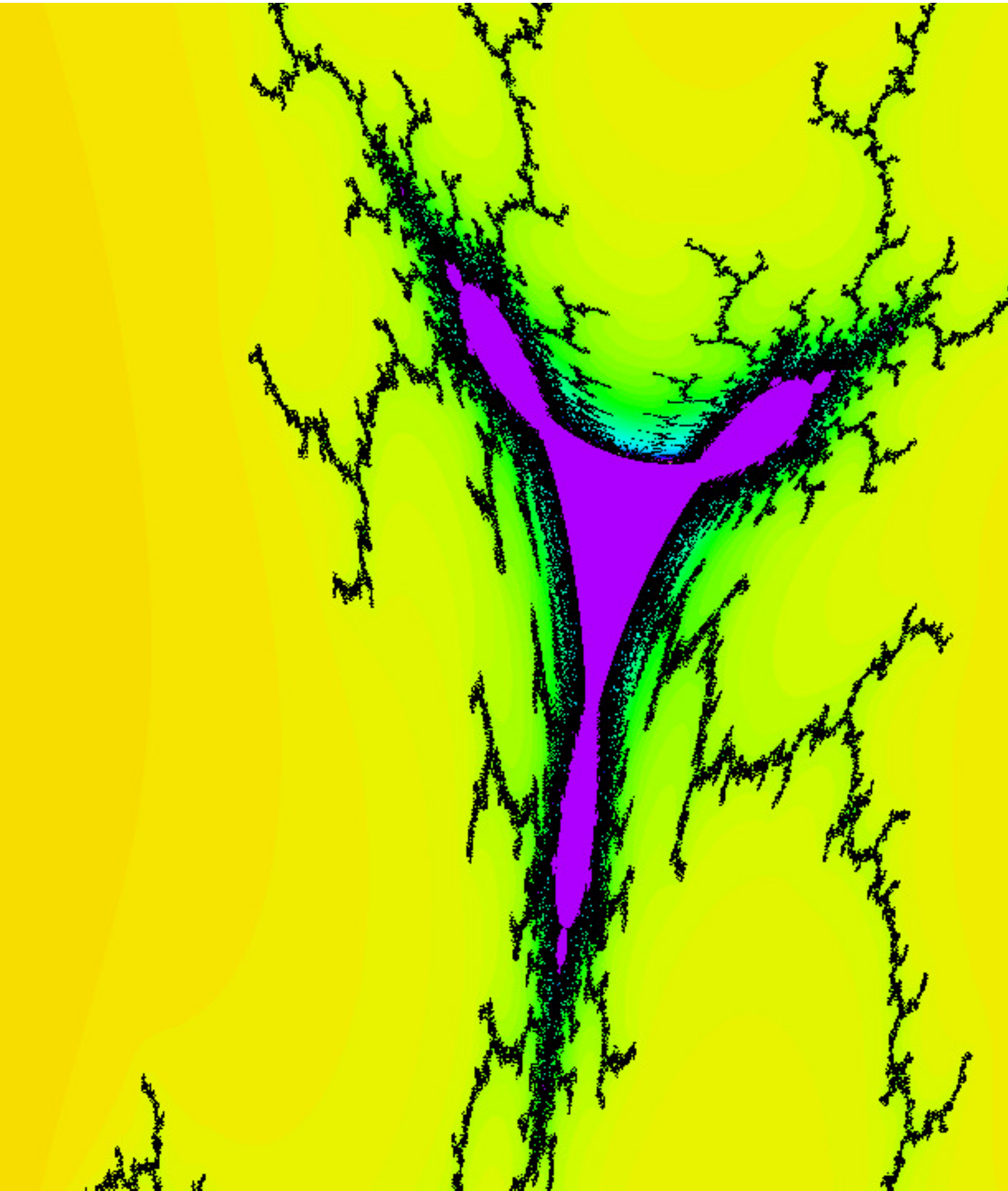,height=2.5in}}
\caption{\label{f-cubtric}\it
Detail of Figure $\ref{f3a}$, showing a small tricorn. $($Window:
$[-.627,\,-.62]\times[.47,\,.477]$ in the real $(A,\,b)$-plane.$)$}
\end{figure}
\medskip

{\bf Decomposable Cases.}
Next let $S_D$ be the \textbf{\textit{~decomposable~}}
 scheme of weight two, as illustrated in Figure \ref{f-wt2}D. The identity 
correspondence $\bio_0(s)=s$ gives rise to a 
real form $\p^{S_{D}}(\gamma_{\bio_0})$ consisting of pairs of non-interacting 
quadratic maps
$$ (s_1,\,z)~\mapsto~(s_1,\,z^2+c_1)\,,\quad(s_2,\,w)~\mapsto~(s_2,\,w^2+c_2)$$
with real coefficients $c_1$ and $c_2$.
 The corresponding
connectedness locus in the real $(c_1,\,c_2)$-plane is just the Cartesian
product $[-2,\,1/4]\times[-2,\,1/4]$ of two copies of the real part of the
Mandelbrot set: See Figure \ref{f4c}, where the chaotic region is
shown in black. Similarly, the nontrivial involution
 $\bio_1:s_1\leftrightarrow s_2$
gives rise to a real form $\p^{S_{D}}(\gamma_{\bio_1})$ consisting of pairs of 
complex holomorphic maps which are complex conjugates of each other
$$ (s_1,\,z)~\mapsto~(s_1,\, z^2+c)\,,\qquad (s_2,\,w)~\mapsto~(s_2,\,w^2
+\overline c)\,.$$
In this case, the connectedness locus in the $c$-plane
is just the Mandelbrot set.\medskip

\begin{figure}[ht!]
\centerline{\psfig{figure=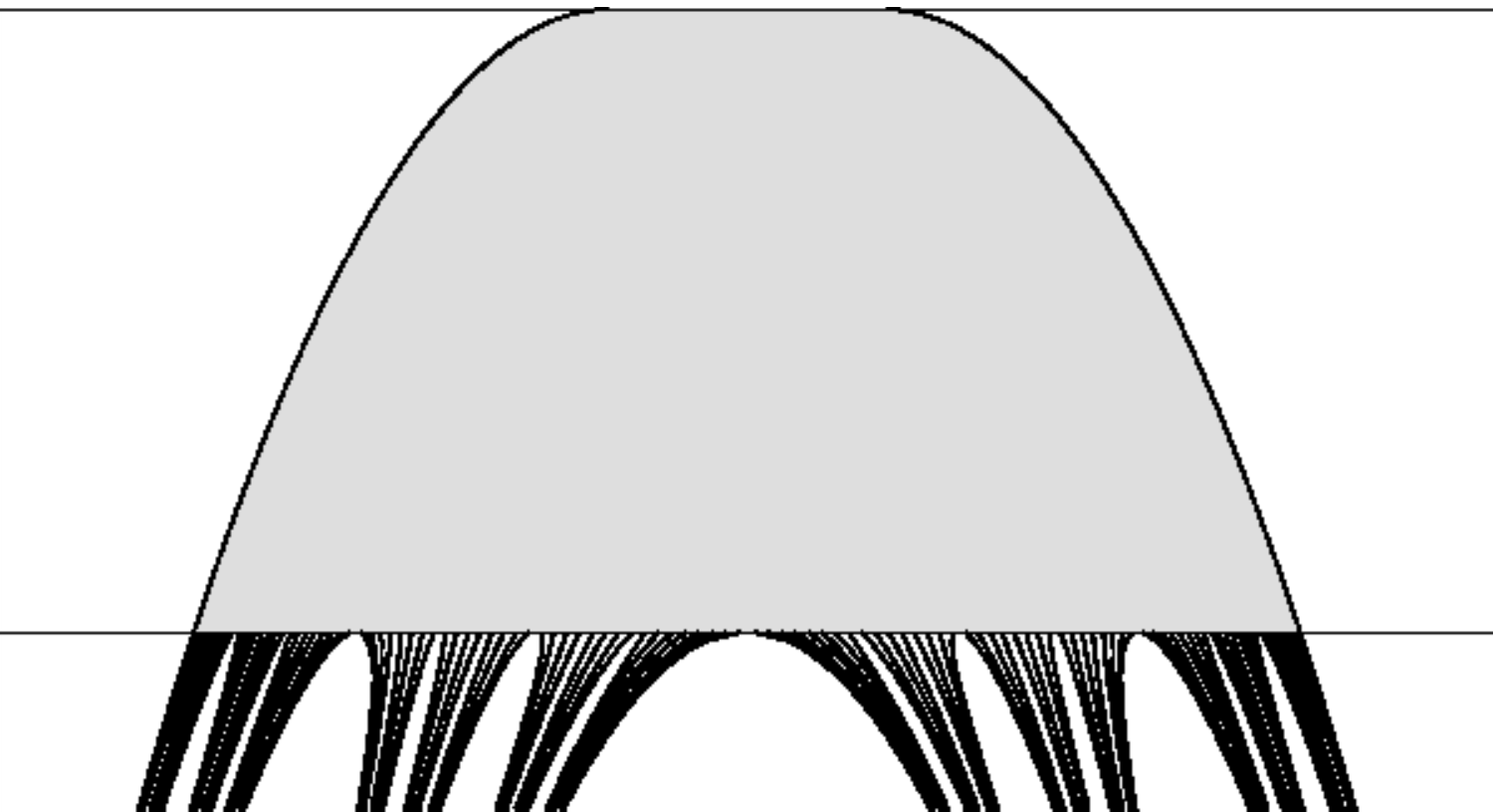,height=2.1in}}
\caption{\it\label{f5}
 Connectness locus for the real form $\p^{S_C}(\gamma_{\bio_0})$}
\end{figure}

{\bf The Capture Component.}
Finally, the \textbf{\textit{~capture~}} scheme $S_C$ of Figure \ref{f-wt2}C,
with $$F:s_1\mapsto s_2\mapsto s_2\,,$$
has a unique real form $\p_{+,+}^{S_C}(\gamma_{\bio_0})$
consisting of maps
\begin{equation}\label{e-cap}
(s_1,\,z)~\mapsto(s_2,\,z^2+c_1)\,\qquad(s_2,\,w)~\mapsto~(s_2,\,w^2+c_2)\,,
\end{equation}
with real $c_1,\,c_2$.  The corresponding connectedness locus 
in the $(c_1,\,c_2)$-plane is the shaded region in Figure \ref{f5}.
It can be described as the set of all
real pairs $(c_1,\, c_2)$ such that both $c_1$ and $c_2$ have bounded orbit
under the map $w\mapsto w^2+c_2$.
In fact $c_2$ has bounded orbit if and only if $-2\le c_2\le 1/4$
(corresponding to the region between the two parallel lines in the figure),
 while $c_1$ has bounded orbit 
for points in the shaded region, and also for points along the uncountable
 family  of curves leading off to infinity below this shaded region.

\begin{rem}\label{r-cap}
If we rely only on Remark \ref{r-sgn-ch}, then we would expect to find
a second real form $\p^{S_C}_{-,+}(\gamma_{\bio_0})$ consisting of maps
\begin{equation}\label{e-alt}
(s_1,\,z)~\mapsto(s_2,\,-z^2+c_1)\,\qquad(s_2,\,w)~\mapsto~(s_2,\,w^2+c_2)
\end{equation}
with real $c_1,\,c_2$.
However, the two forms (\ref{e-cap}) and (\ref{e-alt}) are actually isomorphic
under the complex change of coordinates
$$ \eta(s_1,\,z)=(s_1,\,iz)\,,\quad \eta(s_2,\,w)=(s_2,\,w)\,.$$
On the other hand, if we
consider (\ref{e-cap}) and (\ref{e-alt}) simply as defining
maps from $|S|\times\R$ to itself, and allow only real changes of coordinate,
then these two forms really are non-isomorphic. This discrepancy
between real coordinate changes and complex coordinate changes is closely
related to Remark \ref{r-caution}.
\end{rem}\end{ex}\medskip

To conclude, note the following analogue of the results
 in \S\ref{s4} and \S\ref{s5}.

\begin{theorem}\label{t-6.4}
 Every hyperbolic component in a real connectedness locus of weight $\w$ is a
 topological $\w$-cell with a unique ``center'' point, and is real analytically
 homeomorphic to a uniquely defined principal hyperbolic component $H_{0}^S(
 \gamma)$, or to a suitably defined space $B^S(\gamma)$ of Blaschke
products, under a homeomorphism which is uniquely determined up to the action
of the appropriate finite symmetry group. 
\end{theorem}
\medskip

\noindent The proof involves going through the arguments in previous
sections, keeping track of the extra involution, and is not difficult.
\QED
\medskip

\section{Polynomials with Marked Critical Points.}\label{s7}

By a {\textbf{\textit{critically marked\/}}} polynomial map of degree $\w+1$ we
 will mean a polynomial map $f$ together with an ordered list
 $(c_1\,,\,\ldots\,,\, c_{\w})$ of its critical points. Even if $f$ is a real
polynomial, this list must include all complex critical points, with a critical
point of multiplicity $m$ listed $m$ times, so that the
derivative is given by $$f'(z)~=~(\w+1)(z-c_1)\cdots(z-c_{\w})\,,$$
whenever $f$ is a monic.
As an example, Branner and Hubbard \cite{BH} studied critically marked cubic
polynomials, using the monic centered normal form
$$ f(z)~=~z^3-3a^2z+b\,,$$
with ordered list of critical points $(a,\,-a)$.

Similarly, we can define the concept of a critically
marked Blaschke product.
All of the principal results of the previous sections
extend to the marked case. In particular, for any mapping scheme
$S$ we can define a space $\p_{\rm cm}^S$ of marked polynomial maps
and a space $\B_{\rm cm}^S$ made up out of marked Blaschke products.
Then any hyperbolic component of type $S'$ in the marked connectedness
locus $\cl_{\rm cm}^{S}\subset \p_{\rm cm}^{S}$ is canonically homeomorphic
to $\B_{\rm cm}^{S'}$. The one step in this program which causes additional
difficulty is the analogue of Lemma~\ref{l-cell}, showing that the appropriate
spaces of critically marked Blaschke products are
 topological cells. For this we need the following result.
(Compare \cite{Bo}, \cite{Z}.)
\medskip

\begin{theorem}[Bousch, Zakeri] A Blaschke product
 $\beta:\overline\D\to\overline\D$ of degree $d=\w+1$ which fixes the points
$0$ and $1$ is uniquely determined by its critical points $c_1,\,\ldots,\, c_\w$,
which can be arbitrary points of the open unit disk. Hence the space of all
such maps 
 is diffeomorphic to the $\w$-fold symmetric product of $\D$ with itself.
In particular, this space is a topological cell of dimension $2\w$.
\end{theorem}

As a substitute for zeros centered maps in the critically marked case,
it seems natural to work with the space of 1-anchored critically
marked Blaschke products which are \textbf{\textit{\,critically centered\,}}
in the sense that $c_1+c_2+\cdots+c_\w=0$. However, in order to determine
such a Blaschke product $\beta$
 uniquely, we need one more piece of information, namely the value  $\beta(0)$.
(Compare the statement that a monic polynomial is uniquely determined
by its critical points together with its constant term.) Using the
 Bousch-Zakeri theorem, is not difficult
to check that the space of all such critically marked and centered
Blaschke products of degree $\w+1$ is a topological $2\w$-cell. (For the
special case $\w=0$, this definition doesn't make sense, so we simply define
the corresponding space of Blaschke products to consist of the identity
map only.)

It then follows easily that the  corresponding model space $\B_{\rm cm}^S$,
made out of critically marked Blaschke products which are either
fixed point centered or critically center is also a
 topological cell.\medskip
\medskip

The concept of a \textbf{\textit{~real form~}}
 for the space $\p^S_\cm$ can be defined in analogy
with the discussion in \S\ref{s6}. In general, the space $\p^S_\cm$ has more
real forms than $\p^S$. (This is closely related to the fact that $\p^S_\cm$
has more symmetries than $\p^S$.) As an example,
the space $\p^3_\cm$ of critically marked cubic maps has four distinct real
forms, which can be put in the form
$$ f(z)~=~\pm(z^3-3a^2z+b)\,.$$
Here the critical points $\{a,\,-a\}$ are either real or pure imaginary,
and the sign is either $+$ or $-$, while $b$ is always real.
 More generally, for $2d-1\ge 3$ the space
$\p^{2d-1}_\cm$ has $2d$ distinct real forms which can be labeled
 by the initial sign together
with the number of real critical points. Similarly $\p^{2d}_\cm$ has $d$
real forms.

\section{Rational Maps.}\label{s-rat}

Many of the constructions from this note can be 
applied also to hyperbolic rational maps
$f:\widehat\C\to\widehat\C$. However, there is a difficulty with the statements,
since there is no convenient normal form which works for all rational maps,
and a difficulty with the proofs since the boundary of a critical hyperbolic
 Fatou component need not be a Jordan curve.
We can deal with the first problem by introducing
 a suitable moduli space. (Compare \cite{M3}.)

\begin{definition}\label{d-fm}
By a \textbf{\textit{\,fixedpoint-marked rational map\,}}
$(f;\,z_0,\,z_1,\,\ldots,\,z_d)$
 will be meant a  
rational map $f:\widehat\C\to\widehat\C$ of degree $d\ge 2$, together with
an ordered list of its $d+1$ (not necessarily distinct)
fixed points $z_j$.
\end{definition}

\begin{lem}\label{l-fm} 
The space $~{\rm Rat}^d_{\rm fm}~$ of all such
fixedpoint-marked maps of degree $d$ is a smooth manifold of complex
dimension $2d+1$.
\end{lem}

\begin{proof}
First consider the open subset consisting of all points
 of $~{\rm Rat}^d_{\rm fm}~$ with $f(\infty)\ne\infty$. Each such $f(z)$
can be written uniquely as a quotient
$p(z)/q(z)$ of two polynomials with $q(z)$ monic of degree $d$.
The fixed point equation then takes the form
\begin{equation}\label{e-fm}
z\,q(z)\,-\,p(z)~=~(z-z_0)(z-z_1)\cdots(z-z_d)~,
\end{equation}
where the $z_j$ are the (not necessarily distinct) fixed points of $f$.
Given $q(z)$ and the $z_j\in\C$,
we can solve uniquely for $p(z)$. Here $p(z)$ will be relatively prime to
$q(z)$ if and only if $q(z_j)\ne 0$ for all $j$. Thus we have a well
behaved coordinate neighborhood. Similarly, for each integer $0\le n\le d$
the set of $f$ satisfying $f(n)\ne n$ is a coordinate neighborhood.
The entire space is covered by these $d+2$ 
coordinate neighborhoods, since a map of degree $d$ can have at most
$d+1$ fixed points.
\end{proof}

If we conjugate such a fixedpoint-marked rational map by a M\"obius 
automorphism
$~g:\widehat\C\to\widehat\C$, then we obtain a new fixedpoint-marked map
$$ (g\circ f\circ g^{-1} \,  ;~ g(z_0),\,g(z_1),\,\ldots,\,g(z_d))~. $$
The quotient space of ${\rm Rat}^d_{\rm fm}$ under this action of the
M\"obius group will be called the \textbf{\textit{~moduli space~}}
 $~{\mathcal M}^d_{\rm fm}~$ for fixedpoint-marked maps. This moduli
 space is a non-compact complex algebraic variety of
dimension $2d-2$. The action of the M\"obius group is clearly free on the open
subset consisting of points of ${\rm Rat}^d_{\rm fm}$ with at least three
distinct fixed points. Thus ${\mathcal M}_{\rm fm}^d$ has
possible singularities only within 
the subvariety consisting of 
 conjugacy classes with at most two distinct fixed points.

By a \textbf{\textit{~hyperbolic component~}} in
${\mathcal M}_{\rm fm}^d$ will be meant a connected component
in the open subset consisting of all conjugacy classes of hyperbolic
fixedpoint-marked maps. Each such hyperbolic component is a smooth
complex manifold, since every hyperbolic map has $d+1$ 
 distinct fixed points.  By definition, such a hyperbolic component
belongs to the \textbf{{\textit~connectedness locus~}} if
its representative maps have connected Julia set.


\begin{theorem}\label{t-rat}
Every hyperbolic component $\mathcal H$ in the connectedness locus of 
${\mathcal M}^d_{\rm fm}$
 is canonically homeomorphic to the model space ${\mathcal B}^S$,
where $S=S_f$ is the mapping scheme for a representative map $f$.
In particular, every such $\mathcal H$ is simply connected,\footnote
{Here it is essential that we work with fixedpoint-marked maps. Compare
\cite{Mc2}, which exploits the fact that hyperbolic components in the
space ${\rm Rat}^d$ of all (unmarked) degree $~d~$ rational maps often have an
interesting fundamental group.} with a unique critically finite point.
Similarly, ${\mathcal H}$ is biholomorphic to the standard model of Definition
\ref{d-s-mod}. 
\end{theorem}

The proof will be based on the following preliminary result.
Consider rational maps with just three
marked fixed points. If these three points are distinct, then there
is a unique M\"obius conjugate with these fixed points respectively at
zero, one, and infinity. The resulting map can be written uniquely
as a quotient \hbox{$f(z)=p(z)/q(z)$} of two relatively prime polynomials,
with $p(0)=0$ and $p(1)=q(1)$, where $p(z)$ is monic of degree $d$,
and where $q(z)$ has degree at most $d-1$.
It follows easily that polynomials of this form can be parametrized by an
open subset of the coordinate space  $\C^{2d-2}$.

\begin{theorem}\label{t-rat2} Let $H$ be a connected component in the
space of all hyperbolic rational maps of degree $d$ which have this
normal form, with fixed points at zero, one and infinity, and which
have connected Julia set.
Then $H$ is homeomorphic to the model space $\B^S$, where
 $S$ is the mapping scheme for a representative map in $H$. 
\end{theorem}

 For the proof of this preliminary theorem, we will
 need a concept of ``{\it boundary marking\/}'' for each map
$f:\widehat\C\to\widehat\C$ belonging to  $H$. However,
since the topological  boundary of a Fatou component $U$ for $f$
 may not be a Jordan curve, we must work with a modified
 concept of boundary.

\begin{lem}\label{l-fbd} 
Any Riemann surface $U$ which is conformally isomorphic to the open disk $\D$
can be canonically embedded into an \textbf{\textit{\,ideal compactification}}
$\,\widehat U\,$ which
is diffeomorphic to the closed disk $\overline\D$, and which is natural
in the sense that any proper holomorphic map $U\to U'$ between two
such Riemann surfaces extends to a smooth map $~\widehat U\to\widehat U'$.
In the special case where $~U~$ is an open subset of the Riemann sphere
 $\widehat\C$, then the boundary
 $\partial\widehat U$ $($the \textbf{\textit{ideal boundary\,}} of $U)$
 can be identified with the  set of prime ends of $U$.
In particular, whenever the topological closure
 $\overline U\subset\widehat\C$ is  locally
 connected, the identity map of $~U~$ extends to a continuous
map from $\widehat U$ onto $\overline U$.
\end{lem}

\begin{proof} (Compare \cite{Mc1}.)
Choosing any conformal isomorphism $\phi:U\to\D$, we can pull
the Euclidean metric of the unit disk back to $U$, and then form the metric
completion $\widehat U$. It is not difficult to check that this completion
is independent of the choice of $\phi$. In fact, using Lemma~\ref{l-bp},
we see easily that the differentiable structure of $\overline\D$, pulled back to
$\widehat U$, is independent of the choice of $\phi$. Similarly, if
$U'$ is another Riemann surface conformally equivalent to $\D$, it follows
from Lemma~\ref{l-bp}
that any proper holomorphic map $U\to U'$ extends to a smooth map between ideal
compactifications. The final properties, in the case where $U$ is an open
subset of the Riemann sphere follow easily from Carath\'eodory theory.
(See for example \cite[\S17.13]{M4}.)
\end{proof}

\begin{proof}[Proof of Theorem \ref{t-rat2}]
The argument now proceeds just as in the polynomial case. Choose a basepoint
$f_1\in H$ and form the
space $\widetilde H$ consisting of triples $(f,~q,~\biota)$, where
$f\in H$, where $q$
is a boundary marking sending each critical or postcritical Fatou
component $U$ to a point $q(U)\in\partial\widehat U$ with
 $q\big(f(U)\big)=f\big(q(U)\big)$, and where $\biota$ is an
 isomorphism from $S=S_{f_1}$ to $S_f$. 

We must first show that every connected component of $\widetilde H$ is
a (possibly trivial) covering space of $H$. In particular, we must
show that the boundary marking $q$ deforms continuously as we deform
$f$. Clearly each repelling periodic point deforms continuously as we
deform $f$. Ma\~n\'e, Sad and Sullivan \cite{MSS} showed that this extends to a
continuous deformation of the entire Julia set. Then Slodowski
\cite{Sl} showed that this deformation can be extended to an isotopy of the
entire Riemann sphere. In particular, the closure of 
each Fatou component deforms
continuously, and it follows that each prime end deforms continuously.
Thus the boundary marking $q(U)$ also deforms continuously.
It then follows easily that every point of $H$ has a neighborhood
which is evenly covered under the projection $\widetilde H\to H$.

Just as in \S\ref{s3}, each component of $\widetilde H$ is also a
covering of the simply connected model space ${\mathcal B}^S$, and
hence projects homeomorphically onto 
${\mathcal B}^S$.
It follows that each component of $\widetilde H$ contains a unique critically
finite point, and hence must map homeomorphically onto $H$,
 as required.
\end{proof}

\begin{proof}[Proof of Theorem \ref{t-rat}]
Given a point in the component ${\mathcal H}\subset{\mathcal M}^d_{\rm fm}$,
note that the $d+1$ fixed points for a representative map
 are necessarily distinct. Using only the
first three marked points, we can obtain a unique representative $f$
 in the corresponding component $H\subset{\rm Rat}^d$ of Theorem \ref{t-rat2}.
 This defines a projection ${\mathcal H}\to H$ which is clearly
 a covering map, since the distinct
fixed points vary smoothly as we deform the map $f$. Since $H$ is simply
connected, it follows that this covering map is a homeomorphism.
 Finally, the proof of Theorem \ref{t-5.1} extends easily
to this more general context. 
\end{proof}

\begin{rem}\label{r-antirat}
It is  also interesting to study real forms of rational maps. (Compare 
\cite{M3}.) There are just two antiholomorphic involutions of the
Riemann sphere up to M\"obius conjugation, namely the complex conjugation
operation with $\R\cup\infty$ as fixed point set, and the
\textbf{\textit{~antipodal map~}} $\gamma(z)= -1/\overline z~$ which has
no fixed points. A rational map commutes with complex conjugation if
and only if it is a quotient of polynomials with real coefficients.
The family of rational maps commuting with the antipodal map is more
 interesting. (See \cite{BM}.)
It includes many odd degree rational functions, such as
$z\mapsto z^{2n+1}$, but no even degree functions. In fact any continuous map
of the sphere
which commutes with the antipodal map must have odd degree by a classical
theorem of Borsuk.
\end{rem}

\medskip

{\bf Quadratic Rational Maps.} 
The quadratic moduli space ${\mathcal M}^2_{\rm fm}$ can be identified with the
affine variety consisting of all $(\alpha,\,\beta,\,\gamma)\in\C^3$ satisfying
the equation
\begin{equation}\label{e-m2fm}
\alpha\beta\gamma\,-\,\alpha\,-\,\beta\,-\,\gamma\,+\,2~=~0~ .
\end{equation}
Here $\alpha,\,\beta$ and $\gamma$ are the multipliers at the three marked
fixed points. (Compare \cite{M3}. For other work on quadratic rational
maps, see for example \cite{R1}, \cite{R2}.)

If, for example, $~\alpha\beta\ne 1$,~ then we can solve equation (\ref{e-m2fm})
for $\gamma$ as a smooth function of $\alpha$ and $\beta$. On the other hand,
if $\alpha\beta=1$ then equation (\ref{e-m2fm}) reduces to the equality
$\alpha+\beta=2$, which implies easily
that $\alpha=\beta=1$; this corresponds to the case where the two
 corresponding fixed points crash together. {\it It follows
 that the affine variety defined by equation $(\ref{e-m2fm})$ is smooth except
at the point $\alpha=\beta=\gamma=1$ where all three fixed
 points crash together.} For a description of the singularity at this triple
fixed point class, see the discussion following
 Theorem \ref{t-qtm} below. 

In the case $\alpha\beta\ne 1$, putting the first two fixed points at zero
and infinity, a linear change of coordinates will
 put the corresponding rational map into the normal form
\begin{equation}\label{e-nf2}
 f(z)~=~z\frac{z+\alpha}{\,\beta\,z+1}~,
\end{equation}
with the third marked fixed point at $(1-\alpha)/(1-\beta)$. 
On the other hand, in the case
 $\alpha=\beta=\gamma=1$, if we put the triple fixed point
at infinity, then an affine change of coordinates will yield the normal form
$~f(z)=z+1/z$.

\begin{rem} A hyperbolic component in the connectness locus
need not have compact closure within 
${\mathcal M}^d_{\rm fm}$. In the quadratic case, Adam 
Epstein \cite{Ep} has shown that a hyperbolic component $H$ consisting of maps
with two disjoint attracting cycles has compact closure if and
only if neither attracting cycle is a fixed point. This closure $\overline H$
 can be very
difficult to visualize. (See the following non-compact Examples.)
\end{rem}

\begin{ex} \label{ex-sc}
({\bf The Simplest Case}, although it is not very simple).
Let $H_0\subset{\mathcal M}^2_{\rm fm}$ be the hyperbolic component
consisting of all $(\alpha,\,\beta,\,\gamma)$  satisfying (\ref{e-m2fm})
with $|\alpha|<1$ and $|\beta|<1$. The boundary of this
 component seems very difficult  
to visualize, although it is easily described as a semi-algebraic set.
Certainly we must have $|\alpha|\le 1$ and $|\beta|\le 1$ for every
$(\alpha,~\beta,~\gamma)$ in the closure $\overline H$, and using the
holomorphic index formula it is not hard to show that $\Re(\gamma)\ge 1$.
In fact, $\overline H$ is precisely equal to the set of all
$(\alpha,~\beta,~\gamma)$  in $\overline\D\times\overline\D\times\{\gamma~;~
\Re(\gamma)\ge 1\}$ which satisfy the equation (\ref{e-m2fm}). Since $|\gamma|$
is unbounded, it follows that ${\overline H}_0$ is non-compact,
 as noted in \cite{Ep}.

This closure ${\overline H}_0$
 contains the singular point $\alpha=\beta=\gamma=1$,
which leads to rather bad behavior. For example, although ${\overline H}_0$
is simply connected, it acquires a free cyclic fundamental group if
we remove the singular point. This awkward behavior disappears if
we eliminate the singularity by passing to the 2-fold branched covering
space which is branched only over the singular point. (Compare the discussion
of the ``totally marked'' moduli space following Remark \ref{r-slice} below.)
In fact, the corresponding hyperbolic component in the covering space
has boundary homeomorphic to an open 3-cell. 
\smallskip

Epstein has pointed out that there is a completely analogous example
for cubic polynomial maps. If we look at the space of fixed point
multipliers $(\alpha,\beta,~\gamma)$ in this case, there is again a single
relation, which now takes the form
\begin{equation}\label{e-cub}
 3-2(\alpha+\beta+\gamma) +(\alpha\beta+\alpha\gamma+\beta\gamma)~=~0~.
\end{equation}
If $|\alpha|<1$ and $\beta|<1$, then it follows from the holomorphic
index formula that $|\gamma-3/2|<1/2$. The closure $\overline H$ of the
hyperbolic component described in this way, can be obtained by intersecting the
locus (\ref{e-cub}) with the product
 $~\overline\D\times\overline\D\times\overline\D_{1/2}(3/2)~$
of the corresponding closed disks.
Again this closure has very peculiar behavior, associated with the singularity
at the triple fixed point. However, if we eliminate the singularity by
passing to the 2-fold covering branched at this point,
 corresponding to the family of monic
polynomials $~z\mapsto z^3+az^2+\gamma z\,$,
 then $\overline H$ becomes a closed topological 4-cell.
 (For a similar example in the
context of Kleinian groups, see \cite[Appendix A]{Mc3}.) 
\end{ex}

\begin{figure}[!ht]
\centerline{\psfig{figure=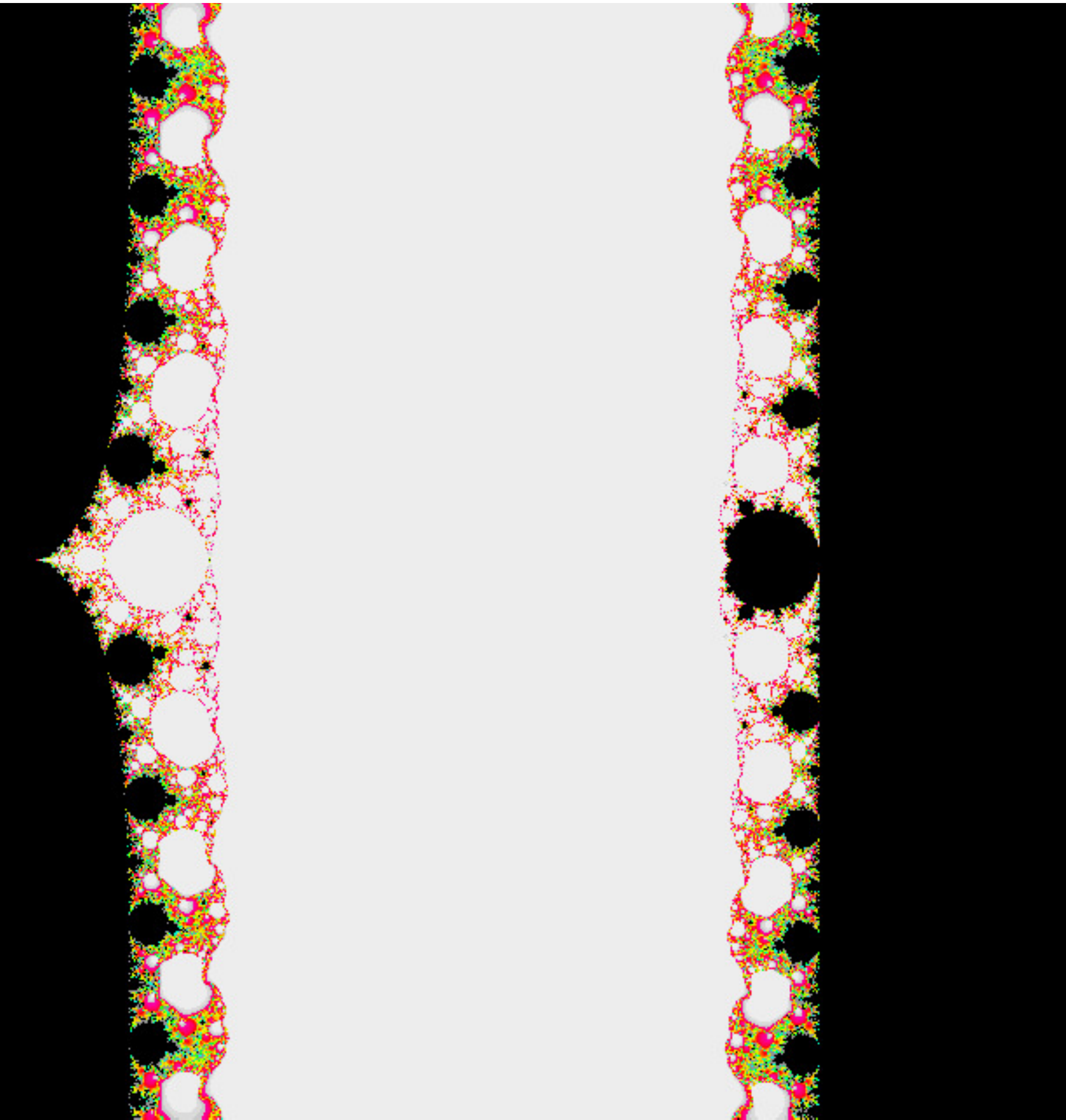,height=2.8in}}
\caption{\it \label{f-rat}
The $a$-plane for the family $~z\mapsto(1-a)/(z(z-a))$.}
\end{figure}


\begin{ex}\label{ex-frac} 
Such examples with semi-algebraic boundary are presumably very rare.
 Here is an example of a hyperbolic component with fractal boundary.
Figure \ref{ex-frac} shows the $a$-parameter plane for the family of maps
\begin{equation}\label{e-frat}
f_a(z)~=~\frac{a-1}{z(a-z)}
\end{equation}
with critical points at $\infty$ and $a/2$, normalized so that $f_a(1)=1$.
Note the superattractive
period two orbit $\infty\leftrightarrow 0$.
(Compare  \cite{Du} and \cite{T}, who use a different parametrization.)
Here the value $a=1$ (the cusp point at right center) must be excluded.
In fact the multiplier $f_a'(1)=(2-a)/(a-1)$ tends to infinity as $a\to 1$.
 The central white
region in this figure is part of the hyperbolic component centered at
$f_0(z)=1/z^2$, and consisting of maps with a bitransitive
 attracting orbit of period two.
This component $H$ in the $a$-parameter plane
 has fractal boundary, and has non-compact closure (for example, in the
neighborhood of $a=1$).
The small white regions in the figure
 correspond to capture components where both critical orbits converge
to the same attracting orbit, but only one critical point lies in a periodic
Fatou component. The black regions correspond to everything else.

 The fixed points
are not marked in this family (\ref{e-frat}), but the fractal nature of
$\partial H$ would remain if we replace $H$ by the corresponding hyperbolic
component $H_{\rm fm}$ in the full moduli space
${\mathcal M}^2_{\rm fm}$ with marked fixed points.  For all maps
in this hyperbolic component, note that the Julia set
is a simple closed curve separating the two Fatou components.

It is interesting to note that the closure
${\overline H}_{\rm fm}\subset {\mathcal M}^2_{\rm fm}$
 is not only non-compact, but has at least three distinct ends
 (conjecturally exactly three).
In fact, if we pass to infinity within moduli space, then at least one of
the multipliers $~\alpha,\,\beta,\,\gamma~$ must have norm tending to infinity.
Since none of the fixed points are attracting, it follows from equation
(\ref{e-m2fm}) that the other two multipliers must
remain bounded, tending towards conjugate points
on the unit circle. (For example as $a\to 1$ along the real axis in the family
(\ref{e-frat}), the other two multipliers tend to $-1$.)
Thus there are at least three essentially different
ways of tending to infinity within $\overline H_{\rm fm}$.
\end{ex}

\begin{rem}\label{r-slice}
One would like to be able to form a 
submanifold of the moduli space ${\mathcal M}^d_{\rm fm}$ by requiring
 one or more critical points to be periodic of specified
period. 
Since our fixed points
are already numbered, we can easily specify  a subvariety of
 ${\mathcal M}^d_{\rm fm}$ by requiring one or more of the {\it fixed}
 points to be critical of specified multiplicity. 
However, if we want a family of maps such that
some critical point is periodic of period $p>1$, then we must either
restrict attention to one hyperbolic component as in Remark~\ref{r-subcurve},
or forget the
fixed point marking as in Example~\ref{ex-frac}, or else
pass to a branched covering space of ${\mathcal M}^d_{\rm fm}$
by also marking this critical point.\smallskip

Define a rational map of degree $d$ to be \textbf{\textit{\,totally marked\,}}
 if we have specified an ordered list, not only of its $d+1$ fixed points,
but also of its $2d-2$ critical points.
I will describe only the quadratic case, which is easier to deal with.

\begin{theorem}\label{t-qtm}
The moduli space ${\mathcal M}^2_{\rm tm}$ for totally marked quadratic maps
is biholomorphic to the smooth simply connected affine variety $V$
consisting of all $~(\,x_1,\,x_2,\,x_3)\in\C^3~$
which satisfy the equation
\begin{equation}\label{e-qtm}
~x_1+x_2+x_3~+~x_1\,x_2\,x_3~=~0~.
\end{equation}
In terms of these coordinates, the multiplier of a representative
 rational map $f$ at the $k$-th fixed point of $f$ is given by
\begin{equation}\label{e-mu}
 \lambda_k~=~1+x_hx_j~,
\end{equation}
where $\{h,\,j,\,k\}$ can be any permutation of $\{1,\,2,\,3\}$.
If we switch the numbering
of the two critical points, then all of the $x_j$ change sign; while if we
renumber the three fixed points, then the $x_j$ are permuted and multiplied by
either $+1$ or $-1$ according as the permutation is even or odd.
\end{theorem}

 The proof will be given below.\smallskip

It follows that the projection from $~{\mathcal M}^2_{\rm tm}~$ to
 $~{\mathcal M}^2_{\rm fm}~$ is a smooth 2-fold branched covering, branched
only over at the point $~x_1=x_2=x_3=0\,$, which maps to $~\lambda_1=\lambda_2=\lambda_3=1$.
It follows easily that any
hyperbolic component $H\subset{\mathcal M}^2_{\rm tm}$ projects
diffeomorphically onto a corresponding hyperbolic component
$H'\subset{\mathcal M}^2_{\rm fm}$.

Given any integer $p\ge 1$, we can define a subvariety 
  ${\rm Per}_p(0)\subset{\mathcal M}^2_{\rm tm}$ by requiring that
the first
marked critical point should be periodic of period $p$, and a ``dual''
subvariety ${\rm Per}^*_p(0)$ by requiring that the second marked critical
point has period $p$. As in Remark \ref{r-subcurve}, it follows that for each
hyperbolic component $H$ which intersects one of these curves, the intersection
is a topological 2-cell with a unique critically finite point.

\begin{lem}\label{l-per0}
Each periodic curve ${\rm Per}_p(0)\subset{\mathcal M}^2_{\rm tm}$
is a smooth complex\break
 1-manifold. Furthermore, for each $p,\,q\ge 1$, the curves
${\rm Per}_p(0)$ and ${\rm Per}^*_q(0)$ intersect transversally in a finite
number of points.
\end{lem}

The proof is completely analogous to the corresponding proof for cubic
polynomial maps, as given in \cite[\S5]{M5}. Details will be omitted.\qed
\smallskip

 However, this discussion leaves three unanswered questions:
\begin{quote}
$\bullet$ How can we count the number of points in this transverse intersection?
(Compare \cite{M5}.)\smallskip

\noindent$\bullet$ How can we compute the Euler characteristic of the curve
  ${\rm Per}_p(0)$? (Compare \cite{BKM}.)\smallskip

\noindent$\bullet$ Is this curve always connected?
\end{quote}

\begin{proof}[Proof of Theorem \ref{t-qtm}]
A totally marked quadratic map
$$~(f,\,z_1,\,z_2,\,z_3,\,c_1,\,c_2)~$$
 is uniquely determined by the 5-tuple
$~(z_1,\,z_2,\,z_3,\,c_1,\,c_2)~$ of fixed points and critical points.
(Compare \cite[\S6]{M3}.) In fact, if we put $c_1$ at
the origin and $c_2$ at infinity, then the map $f$ takes the form
\begin{equation}\label{e-nf1} f(z)~=~\frac{az^2+b}{cz^2+d}~.
\end{equation}
Thus there is a fixed point at infinity if and only if $c=0$, and
 the finite fixed points satisfy the equation $~cz^3-az^2+dz-b=0\,$.
Clearly the collection of roots $z_1,\,z_2,\,z_3$
uniquely determines
 the point $(a:b:c:d)$ in projective 3-space, and hence uniquely determines
the map $f$.


Now consider the three cross-ratios
$$ r_h~=~ \frac{(c_1-z_j)(c_2-z_k)}{(c_1-z_k)(c_2-z_j)}~,$$
where $~(h,\,j,\,k)~$ is to be any cyclic permutation of $~(1,\,2,\,3)\,$.
These are clearly invariant under conformal conjugacy. In fact, they form
a complete conjugacy invariant. Still using the
normal form (\ref{e-nf1}), we see that $r_h=z_j/z_k$. If another triple
has the same ratios $r_j$ then we can get one triple from the other by
 multiplying by some constant $\lambda\ne 0$. This corresponds
 to a conjugation of the form
$$f(z)~\mapsto~ \lambda f(z/\lambda)~.$$
 (Here one has to take care with the
 special case that some $z_j$ is zero and/or some $z_k$ is infinity, but the
 conclusion still follows, using the fact that a double fixed point can never
be a critical point. Details will be left to the reader.)

Recall that a point in the moduli space ${\mathcal M}^2_{\rm fm}$ is determined
by the multipliers at the three fixed points, which we now denote by
 $\lambda_1,\,\lambda_2,\,\lambda_3$. Define a map from the variety $V$ onto
 ${\mathcal M}^2_{\rm fm}$ by setting
$$ \lambda_h~=~1+x_j\,x_k~.$$
Then the required identity (\ref{e-m2fm}),
$$ \lambda_1\,\lambda_2\,\lambda_3\,-\,\lambda_1\,-\,\lambda_2\,-\,\lambda_3\,+\,2~=~0~, $$
is easily verified. Note the identity
\begin{equation}\label{e-xsquare}
x_h^2~=~ 1-\lambda_j\,\lambda_k~.
\end{equation}
In fact $~ 1-\lambda_j\lambda_k\,=\,1-(1+x_hx_k)(1+x_hx_j)~$ which simplifies easily to
\break$~-x_h(x_k+x_j+x_hx_jx_k)~=~x_h^{\,2}~. $\medskip

If $z_j\ne z_k$, then we can put the fixed point $z_j$ at zero and
 $z_k$ at infinity, and write the map as
$$f(z)~=~z\frac{z+\lambda_j}{\lambda_k\,z+1}~,$$
as in equation (\ref{e-nf2}). A brief computation then shows that the two
critical points $~f'(z)=0~$ satisfy the equation $~\lambda_kz^2+2z+\lambda_j\,=\,0\,$,
 with solution
$$ c_i~=~\frac{-1\pm\sqrt{1-\lambda_j\lambda_k}}{\lambda_k}~=~\frac{-1\pm x_h}{\lambda_k}~.$$
To fix our ideas, suppose that 
$$ c_1~=~(-1+x_h)/\lambda_k\,,\quad c_2~=~(-1-x_h)/\lambda_k ~.$$
Then the cross-ratio $r_h$ is given by
$$ r_h~=~\frac{c_1}{c_2}~=~\frac{-1+x_h}{-1-x_h}~=~\frac{1-x_h}{1+x_h}~.$$
Thus the cross ratios $r_h$ and hence the
conjugacy class in ${\mathcal M}^2_{\rm tm}$ are uniquely determined by the
coordinates $x_1,\,x_2,\,x_3$, yielding a holomorphic map from $V$ to
${\mathcal M}^2_{\rm tm}$.

Conversely, given the $r_h$, we can solve for  $~x_h=(1-r_h)/(1+r_h)\,$.
(Using this correspondence $~r_h\leftrightarrow x_h\,$,
 note that the defining identity (\ref{e-qtm})
for the variety $V$ is completely equivalent to the relation 
$$ 1~=~r_1r_2r_3~=~\frac{(1-x_1)(1-x_2)(1-x_3)}{(1+x_1)(1+x_2)(1+x_3)}~,$$
again taking special care with the cases where some $r_h$ is zero or infinity.)
This completes the proof that ${\mathcal M}^2_{\rm tm}$ is biholomorphic to the
affine variety $V$.

As an example, if we take $c_1=0$ and $c_2=\infty$, with $h,\,j,\,k$ in
positive cyclic order, then
$$ x_h~=~\frac{z_k-z_j}{z_k+z_j}~,$$
where the denominator can never be zero.
\medskip

It is easy to prove that $V$ is smooth and connected. In fact it is covered
by three coordinate neighborhoods
$$ V_h~=~\{~(x_1,\,x_2,\,x_3)~\in ~V~~;~ x_j\,x_k~\ne~-1~\}~. $$
Here $~V_1\cup V_2\cup V_3~=~V~$ since the equations $~x_1x_2=x_1x_3=x_2x_3=-1~$
have no simultaneous solution within $V$. For $~(x_1,\,x_2,\, x_3)\in V_h$,
we can solve uniquely for
$$ x_h~=~\frac{-x_j-x_k}{1+x_j\,x_k}$$
as a holomorphic function of the other two variables.\smallskip

We first show that each coordinate neighborhood
 $V_h$ has fundamental group $\pi_1(V_h)~\cong~\Z$. To simplify notation,
note that $V_h$ is biholomorphic to the complement of the quadratic curve
$$ W~=~\{~(x,\,y)~\in~\C^2~~;~~ x\,y\,=\, -1~\}~.$$
Let $~{\mathcal S}\subset\C^2~$ be the real hypersurface consisting of all
 products $~(\xi,\,\eta)=(tx,\,ty)$ with
$~(x,\,y)\in W~$ and $t>1$. 
 Then the complement $~\C^2\ssm(W\cup {\mathcal S})~$
is star-shaped. That is, the line segment joining any point to the origin
is completely contained in \hbox{$~\C^2\ssm(W\cup {\mathcal S})\,$}. Any loop
 in \hbox{$~\C^2\ssm W~$} can be perturbed until it intersects the 
hypersurface ${\mathcal S}$ transversally
 in finitely many points. The homotopy class of such 
 a loop is determined by the number of
transverse intersection points, counted with a sign of $-1$ or $+1$ according
as the imaginary part $~\Im(\xi\,\eta)~$ is increasing or decreasing as
$(\xi,\,\eta)$
passes through ${\mathcal S}$. In fact we can use the star shaped property
to drag the loop $L$ down to the origin except in a small neighborhood of each
 intersection point $(\xi_j,\,\eta_j)$ The part of this loop within this
small neighborhood can then be deformed to a triangular loop $T(\xi_j,\,\eta_j)$
consisting of a line segment from the origin to
 $(1\pm i\epsilon)(\xi_j,\,\eta_j)$, followed by a line segment to
 $(1\mp\epsilon) (\xi_j,\, \eta_j)$, and then followed by a line segment back
to the origin. Since $\mathcal S$ is connected, the homotopy class of this
triangular loop does not depend on the particular choice of $(\xi_j,\,\eta_j)$.
Finally, the composition of two consecutive loops of opposite orientation
is homotopic to the zero loop. Since a standard topological argument
shows that the number of intersections, 
counted with sign, is a homotopy invariant, this proves that
 \hbox{$\pi_1(\C\ssm W)\cong \pi_1(V_h)\cong\Z$.}\smallskip

As an example, consider the loop $L$ in $V$ which is given by
$$\theta~\mapsto (x_1,\,x_2,\,x_3)\quad{\rm with}\quad
x_1~=~1+\epsilon\, e^{i\theta}\,,~~~ x_2=-x_1\,,~~~x_3=0~. $$
Then $L$ is homotopic to a constant in $V_1$
or in $V_2$, since we can simply let $\epsilon$ tend to zero;
and yet it represents a generator of $\pi_1(V_3)$.
It follows that the variety
 $V$ is simply-connected. In fact, each inclusion $V_h\subset V$ induces a
 homomorphism from $\pi_1(V_h)$ onto
$\pi_1(V)$, since it is easy to homotop any loop in $V$ away from the locus
$x_jx_k=-1$. \smallskip

 Further details of the proof of Theorem
\ref{t-qtm} are straightforward.
\end{proof}
\end{rem}\smallskip

\begin{rem}
The field $F\subset \C$ generated by the coordinates $x_j\in\C$
 is an interesting
invariant of the conjugacy class of $f$. It can be characterized as the
 smallest field such that some M\"obius conjugate of $f$ has all critical points
and all fixed points within $F$.
\end{rem}

\appendix
\section{Realizing Reduced Schemes (by Alfredo Poirier\protect\footnotemark[6])\label{sa}}   
\bigskip

The purpose of this appendix is to prove that every reduced mapping scheme can
 be realized by a postcritically finite polynomial. 
In order to do so, we will construct an appropriate Hubbard tree that mimics 
the dynamics of the scheme. 
For the benefit of the reader we recall briefly the main concepts involved in 
the construction of Hubbard trees  following closely \cite{P3}.\medskip

\footnotetext[6]
{Departmento de Ciencias, Secci\'on Matem\'aticas, Pontificia Universidad Cat\'olica del Per\'u, Apartado 1761, Lima 100, Per\'u;  \hspace{.5cm}
{\it email\/}: {apoirie@pucp.edu.pe}}
\setcounter{footnote}{6}

Given a degree $d \ge 2$ postcritically finite polynomial $f$, 
we know that its \textbf{\textit{filled Julia set}} $K(f)$, besides being 
connected, is locally connected.  Call a periodic orbit that contains a 
critical point a \textbf{\textit{critical cycle}}. 
In the postcritically finite setting, a periodic orbit belongs to the Fatou 
set $F(f)$ if and only if it is a critical cycle 
(for details we refer to \cite[Corollary 14.5]{M4}). 

In this postcritically finite case, the polynomial $f$ when restricted to the 
interior of $K(f)$ 
(which happens to be nonempty only when there exists a critical cycle) 
maps each bounded Fatou component 
---always simply connected by the maximum modulus principle--- 
onto some other as a branched covering map. 
Furthermore, all of them are eventually periodic (see \cite[Theorem 16.4]{M4}). 
And also, each component can be uniformized so that in local charts 
$f$ reads $z \mapsto z^m$ for some $m \ge 1$ (see \cite[Theorem 9.1]{M4}). 
More is true. 
If $U$ is a periodic bounded Fatou component, 
then the first return map is conjugate to $z \mapsto z^k$, this time with 
$k \ge 2$. 
In particular, loops of components are in perfect correspondence with critical
 cycles. 
Also, in each component there is a unique point which eventually maps to a 
critical point 
(precisely the one marked as $0$ in local coordinates), its 
\textbf{\textit{center}}. 

It is well known (see for instance \cite[Corollary VII.4.2]{DH1}) 
that given a degree $d \ge 2$ postcritically finite polynomial $f$,  
for any $z \in K(f)$ the sets $K(f)-\{z\}$ and $J(f)-\{z\}$ consist both of a 
finite number of connected components.  
In this way, the filled Julia set can be thought of as arranged in a tree-like
 fashion. 

To get rid of inessentials, we pick a finite invariant set $M$ containing all 
critical points. 
Within $K(f)$ we interconnect $M$ by arcs subject to the extra condition that 
when a Fatou component is met, 
then this intersection consists of radial segments in the associated 
coordinate. 
Douady and Hubbard proved that this construction defines a finite topological 
tree $T(M)$ 
when $M$ together with the intrinsic branching points are considered vertices. 

The vertex dynamics is invariant and carries the endpoints of any edge to 
distinct elements, 
so that it can be extended to a function from $T(M)$ to itself which is one to
 one on each edge 
and is isotopic to $f$, the original map. 
We also keep record of the local degree at every vertex $v$ as $d(v)$. 
In addition, if three or more edges meet at a vertex, 
then their cyclic order should be remembered. 
In other words, we specify how this tree is embedded in the complex plane, 
again, up to isotopy. 

Unfortunately, this data alone is not enough to determine the affine conjugacy
 class of $f$. 
However, if we append enough information to recover the inverse tree, 
then different postcritically finite polynomials yield different structures. 
To formally deal with this condition we introduce angles around vertices 
(this is to be credited again to \cite{DH1}). 
In what follows we measure angles in turns, so that 1 degree measures 1/360 of
 a turn. 
At the center of a component the angle between edges is measured using the 
local chart. 
Near Julia vertices, where $m$ components of $K(f)$ intersect, 
the angle is naturally defined as a multiple of $1/m$. 

These angles satisfy two obvious conditions. 
First, they are compatible with the embedding of the tree. That is, as we
go around a vertex in the positive direction, the successive angles are between
zero and one, and add up to $+1$. Second, they satisfy the identity 
$$
\angle_{f(v)}(f(e),f(e'))=d(v)\angle_{v}(e,e') \mod 1,
$$ 
where $d(v)$ is the local degree at $v$, and $e,e'$ are edges incident at $v$. 
When this further structure is provided, 
we are in front of a \textbf{\textit{Hubbard tree}}. 

\subsection*{Abstract Hubbard Trees}
Now we move in reverse: 
we start with an abstract dynamical tree and we reconstruct the appropriate 
postcritically finite polynomial.

An \textbf{\textit{angled tree $H$}} is a finite simplicial tree  
together with a function\break $e,e' \mapsto \angle_{v}(e,e') \in \mathbb Q/ 
\mathbb Z$ 
which assigns a rational number modulo $1$ to each pair of edges $e,e'$ 
incident at a vertex $v$. 
The \textbf{\textit{angle}} $\angle_{v}(e,e')$ is skew symmetric 
with $\angle_{v}(e,e')=0$ if and only if $e=e'$, 
and satisfies $\angle_{v}(e,e'')=\angle_{v}(e,e')+\angle_{v}(e',e'')$ 
whenever the three edges meet at $v$.  This angle function 
 determines a preferred isotopy 
class of embeddings of $H$ into $\mathbb C$. 

Let $V$ be the set of vertices in $H$.  We specify a 
\textbf{\textit{vertex dynamics}} $f:V \to V$ subject to $f(v) \ne f(v')$ 
whenever $v,v'$ are  end-points of a common 
 edge $e$. 
We consider also a \textbf{\textit{local degree}} 
\hbox{$d:V \to  \{1,2,\dots\}$}. We require that the 
\textbf{\textit{total degree}} \hbox{$d_H=1 + \sum_{v \in V} (d(v)-1)$} must be
 greater than $1$.  By definition a vertex is \textbf{\textit{critical}} if 
$d(v)>1$ and non-critical otherwise.  The \textbf{\textit{critical set}} is 
thus non-empty.  

We require $f$ and the degree $d$ to be related to the angles 
 as follows. 
Extend $f$ to a map \hbox{$f:H \to H$} that carries each edge homeomorphically
 onto 
the shortest path joining the images of its endpoints. We  then need
\hbox{$\angle_{f(v)}(f(e),f(e'))=d(v)\angle_{v}(e,e')$} whenever $e,\, e'$ are 
incident at $v$ 
(so that $f(e),\, f(e')$ intersect at the vertex $f(v)$
 where the angle is measured). 

A vertex $v$ is \textbf{\textit{periodic}} if \hbox{$f^{\circ k}(v)=v$} for some 
\hbox{$k \ge 1$}.  The orbit of a periodic critical point is a 
\textbf{\textit{critical cycle}}.  A vertex is of \textbf{\textit{Fatou type}}
 if it eventually maps to a critical cycle;  else, it is of 
\textbf{\textit{Julia type}} or a \textbf{\textit{Julia vertex}}. 

The distance \hbox{$dist_{H}(v,v')$} between vertices in $H$ counts the number
 of  edges in the shortest path joining $v$ to $v'$.  We call $H$ 
\textbf{\textit{expanding}} if for every edge $e$ whose endpoints $v,v'$ are 
Julia vertices  there is $n \ge 1$ for which we have 
$$dist_{H}(f^{\circ n}(v),f^{\circ n}(v')) > 1.$$ 
(In practice this property must be tested only for adjacent Julia vertices.)

Angles at Julia vertices are rather artificial, so it is better to normalize 
them.  
If $m$ edges \hbox{$e_1,\,\dots,\,e_m$} meet at a periodic Julia vertex $v$, 
then each \hbox{$\angle_{v}(e_i,e_j)$} should be a multiple of $1/m$. 
(Therefore, angles around a periodic Julia vertex convey no information beyond
 the cyclic order of the edges.) 
An angled tree that satisfies this 
 condition around all periodic Julia vertices
 is said to be \textbf{\textit{normalized}}. 

By an \textbf{\textit{abstract Hubbard tree}} ---or simply a Hubbard tree--- 
we mean a normalized angled tree that obeys the expanding condition. 
The basic existence and uniqueness theorem is stated now as follows.
\medskip 
 
\begin{theorem}[Poirier \cite{P3}]\label{t-A1P}
A normalized dynamical angled tree can be realized as the tree associated to a 
postcritically finite polynomial if and only if it is expanding,
or in other words if and only if it is an abstract Hubbard tree. 
Such a realization is unique up to affine conjugation. 
\end{theorem}

Note that there many cases where we can apply this result directly: 
any tree which has no adjacent Julia vertices is trivially expanding. 
For instance, a star-shaped dynamical tree in which a critical cycle
 pivots around a fixed vertex can always be realized. 
\medskip 

Now we are ready to realize a given reduced scheme \hbox{$\overline{S}=(
|\overline{S}|,\, F,\, \w)$} and  
settle the existence of a postcritically finite polynomial $f$ of degree 
\hbox{$d(\overline{S})=\w(\overline{S})+1$}
whose associated reduced scheme \hbox{$\overline{S}(f)$} is isomorphic to
$\overline{S}$. 

\begin{theorem} 
Every reduced scheme can be  realized by a postcritically finite polynomial.
\end{theorem} 

\begin{proof}
First we construct a non-reduced scheme $S$ which reduces to $\overline{S}$. 
This is done by adding new vertices of weight zero 
in such a way that the associated graph $\Gamma(S)$ can be obtained from 
\hbox{$\Gamma(\overline{S})$}
by plotting an extra vertex within each old edge.
(Compare Figures \ref{f-A.1} and \ref{f-A.2}, where the new vertices are 
indicated by small circles.)
More explicitly, 
starting with a reduced scheme $\overline{S}$ with associated map 
$\overline{F}:|\overline{S}| \to |\overline{S}|$, 
construct a non-reduced scheme $S$ with associated map $F:|S| \to |S|$, where
 $|\overline{S}| \subset |S|$, as follows. 
The difference set $|S|-|\overline{S}|$ is to consist of one vertex 
$s^{\sharp}$ for each $s \in |\overline{S}|$, and the map $F:|S| \to |S|$ is 
defined by
$$
F(s)=s^{\sharp} \qquad \hbox{and} \qquad 
F(s^\sharp)=s', \qquad \hbox{where} \qquad s'=\overline{F}(s). 
$$

\begin{figure}[!ht]
\centerline{\psfig{figure=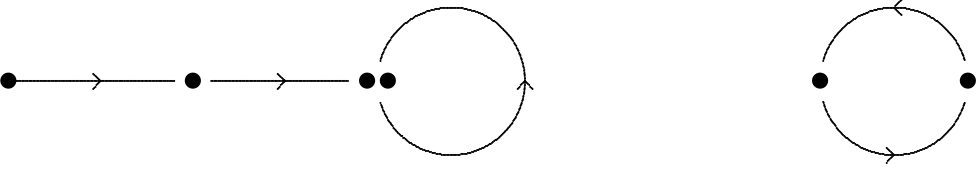,height=.8in}}
\vspace{-.4cm}
\caption{\it \label{f-A.1} A reduced scheme.}
\bigskip

\end{figure}

\begin{figure}[!ht]
\bigskip

\centerline{\psfig{figure=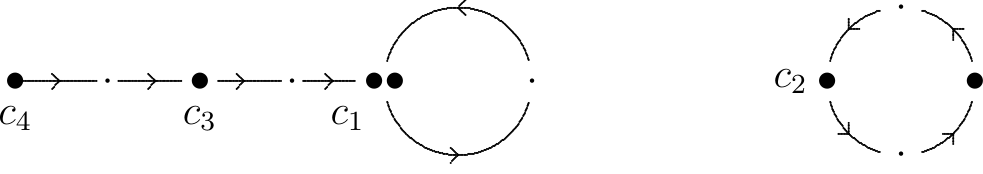,height=.8in}}
\caption{\it \label{f-A.2} The associated non reduced mapping scheme. 
From each cycle we pick a representative ($c_1,\,c_2$). 
All non-periodic critical points are also named ($c_3,\,c_4$).}
\end{figure}

In simple words, $s^{\sharp}$ lies in the middle of $s$ and 
$s'=\overline{F}(s)$, 
so that the main difference between $F$ and $\overline{F}$ is that a vertex 
$s \in |\overline{S}|$ 
now takes an intermediate (artificial) step before reaching 
$s'=F(s) \in |\overline{S}|$. 
Formally, we have $F^{\circ 2}=\overline{F}$ when restricted to 
$|\overline{S}|$. 

By its very definition, every vertex of the form $s^\sharp$ has $s$ as its only 
preimage. 
From this construction it follows readily that $S$ has $\overline{S}$ as its
reduced scheme. 
This scheme is the one that we will bring to life with the help of a suitable 
expanding Hubbard tree. 

\begin{figure}[!ht]
\centerline{\psfig{figure=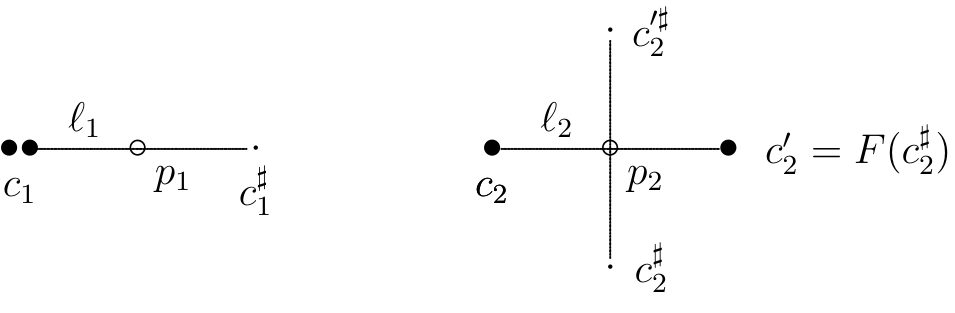,height=1.5in}}
\caption{\it \label{f-A.3} The dynamical graphs corresponding to the two 
cycles in Figure \ref{f-A.2}.}
\end{figure}

Let ${\mathcal C}_i$ be a cycle 
$s_0 \mapsto s_0^{\sharp} \mapsto s_{1} \mapsto \cdots \mapsto s_{n-1}^{\sharp} 
\mapsto s_{0}$ in $S$. 
We join all these $2n$ vertices consecutively around a new vertex $p_i$ 
in order to get a star-shaped symmetrical graph (compare Figure \ref{f-A.3}). 
Mapping $p_i$ to itself, we have a Hubbard tree.
(Here and elsewhere the degree  \hbox{$d(s)=\w(s)+1$} is copied from the
 scheme.)
All these $s_k$ and $s_k^{\sharp}$ belong to a critical cycle and as such are of
 Fatou type.  Hence, 
the dynamics in this graph is expanding since there is only one Julia type 
vertex in sight (the fixed point $p_i$). 
For future reference we pick a critical vertex in the loop 
(for instance $s_0$, which is critical because it belongs to $|\overline{S}|$
where
$\overline{S}$ is reduced) and call it $c_i$. 
Also, the edge between $c_i$ and $p_i$ will be referred to as $\ell_i$.  

\begin{figure}[!ht]
\centerline{\psfig{figure=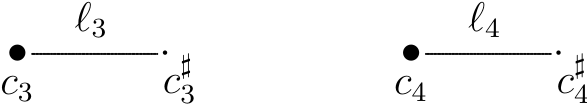,height=.6in}}
\caption{\it \label{f-A.4} The remaining non-cyclic edges.}
\end{figure}

Let $c_1\,,\,\ldots\,,\,c_m$ be the vertices lying within a critical cycle,
 as numbered above. We now have to perform an extra auxiliary construction. 
Let 
$$c_{m+1}~=~s_{m+1}\,,~\dots\,,~c_{m+r}~=~s_{m+r}~ \in ~\overline{|S|}$$
be the non-periodic critical vertices, i.e, those outside a critical cycle.
For \break \hbox{$i=m+1,\dots, m+r$} consider segments $\ell_i$ between 
$c_i=s_{i}$ and  $s_i^{\sharp}$ (compare Figure~\ref{f-A.4}). 

\begin{figure}[!ht]
\centerline{\psfig{figure=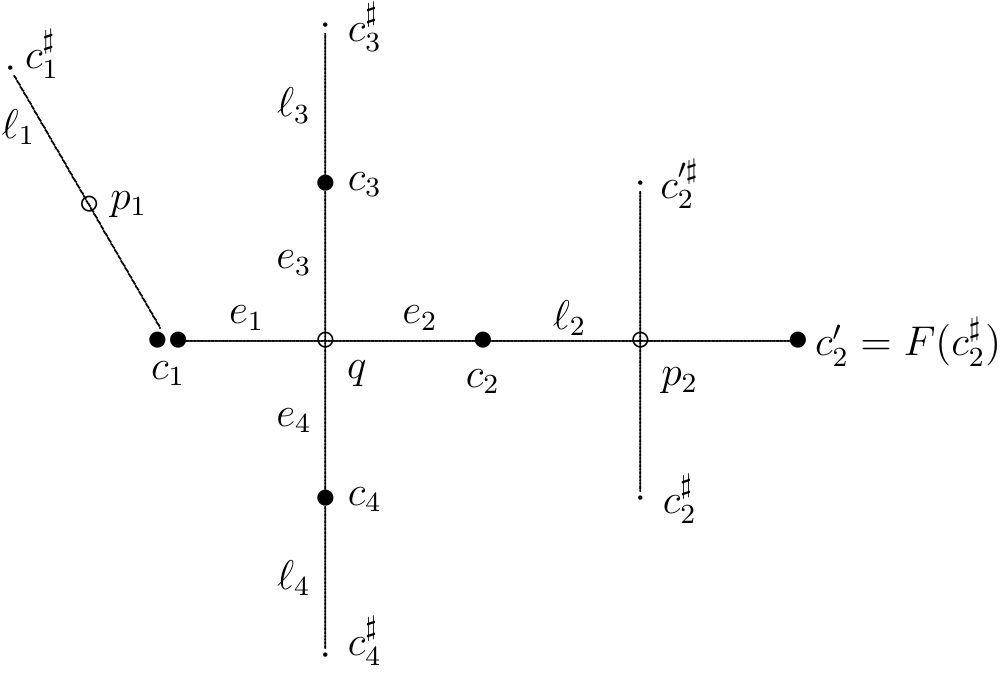,height=2.7in}}
\caption{\it \label{f-A.5} Assembling the pieces together. 
Notice that at $c_1$ the angle should be $1/3~({\rm mod}~1)$. 
By construction,
the Julia type fixed points $p_1,~q,~p_2$ $($represented here by 
 circles$)$
have rotation numbers $1/2,~1,~1/4$ respectively. 
Also recall that $F$ maps $c_4^{\sharp}=F(c_4)$ to $c_3$ 
and   $c_3^{\sharp}=F(c_3)$  to $c_1$.}
\end{figure}

Now we are ready to construct the Hubbard tree. 
For \hbox{$i =1,\cdots, m+r$},  add a segment $e_i$ at $s_i$  
making an angle of $1/d(s_i)$ units with $\ell_i$.
(This is to guarantee a complete folding at $s_i$ when we iterate.)
Merge the $e_i$'s at a new vertex $q$ making a uniform angle of $1/(m+r)$ 
between consecutive edges  (the order here is irrelevant). 
Map the Julia vertex $q$ to itself (compare Figure \ref{f-A.5}).
Since there are no adjacent Julia vertices, we have an expanding Hubbard tree. 
 The essentially unique postcritically finite polynomial that realizes
this Hubbard tree, whose existence is guaranteed by Theorem~\ref{t-A1P}, clearly
has the required reduced scheme $\overline S$. 
\end{proof}
\bigskip

\section{Census of Reduced Schemes}\label{sab}

In this appendix, all schemes are to be reduced. The object is to count
the number $N(\w)$ of distinct isomorphism classes of schemes of weight
\hbox{$\w(S)=\w$} for each small value of $\w$. The counting process can be 
broken down  into a number of smaller steps as follows.\bigskip


{\bf (a)\;}~ Every scheme is uniquely a disjoint sum of connected schemes.
If  \hbox{$~S_1 ,~ S_2,~$} \hbox{$~ S_3,~\ldots~$}  is a list of all connected
 schemes,  then every scheme can be expressed uniquely as a sum
$$   S~=~ S_{i(1)} \,+\,S_{i(2)} \,+\,\cdots\,+\,S_{i(k)} \,, $$
where $k$ is the number of connected components, and where 
$$~i(1)~\le~ i(2)~\le~\cdots~\le~ i(k)~.$$
 The total weight of such a sum is \hbox{$~\w(S)=\sum_j \w(S_{i(j)})$}.
Thus, in order to compute the total number $N(\w)$ of schemes of weight $\w$,
it suffices to know the smaller number $N_c(\w')$ of connected schemes of
weight $\w'$, for every $\w'\le \w$. 
\bigskip

\begin{definition}\label{d-tree}
By a \textbf{\textit{~weighted tree~}} will be meant an acyclic simplicial
complex of dimension $\le 1$ with a preferred \textbf{\textit{~root~}} vertex,
together with a weight function which assigns a positive integer to each
non-root vertex. By definition, the root vertex always has weight zero.
The \textbf{\textit{~trivial~}} weighted tree consists of the root vertex alone,
with no edges. \end{definition}

\begin{figure}[!ht]
\centerline{\psfig{figure=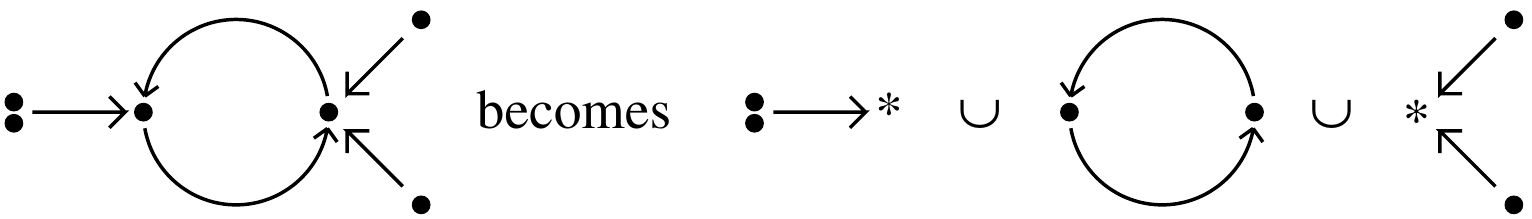,height=.55in}}
\caption{\it \label{fB1} Splitting off the trees from a scheme of weight
 $~\w(S)=6$.}
\end{figure}

{\bf(b)}~ Every connected scheme $S$ consists of
 a central cycle $C$ of weight\break \hbox{$\w(C)\ge 1$}, together with a
(possibly trivial) weighted 
tree $T(s)$ which is pasted onto each vertex \hbox{$s\in |C|$}.
Here the root point of $T(s)$ is to be identified with $s$.
(Compare Figure~\ref{fB1}, where each root point is represented by the symbol
$*\,$.) Thus
$$ \w(S)~=~\w(C)\,+\,\sum_{s\in|C|} \w\Big(T(s)\Big)\,, \qquad{\rm where}\qquad
\w(T)~=~\sum_{t\in|T|,~t\ne *} \w(t)\,. $$
Note that the cycle $C$ can be economically described by
 a symbol of the form\break
\hbox{$(\w_1,\,\w_2,\,\ldots,\,\w_n)$} which is well defined up to cyclic
 permutation. Here the $\w_i$ are positive integers with sum $\w(C)$. Each 
$\w_i$ corresponds to a vertex $s_i$ of weight $\w_i$ which maps to $s_{i+1}$, 
where the subscript $i$ varies over $\Z/n$.
\bigskip

{\bf(c)}~ By a \textbf{\textit{trunk\/}} of a tree will be meant an edge 
incident to the root
point. Thus every non-trivial tree has at least one trunk. Let
$~T_1,\,T_2,\,T_3,\,\ldots ~$ be a list of all trees with only one trunk.
Then a tree with $k\ge 2$ trunks is isomorphic to a unique wedge sum
\begin{equation}\label{e-ws}
   T~=~T_{i(1)}\,\vee\,\cdots\,\vee\,T_{i(k)}
\end{equation}
 of non-trivial trees pasted together at the root point, where
 $i(1)\le\cdots\le i(k)$. Just as in paragraph
 {\bf(a)} above, the total weight is the
sum $\w(T)=\sum_j \w\big(T_{i(j)}\big)\,.$
(Remember 
 that the weight of the root point
 is always zero.)
\bigskip

{\bf(d)}~ If $N_{\rm tree}(\w)$ is the number of distinct trees
 of weight $\w$, and  $N_1(\w)$ is the number of such trees with only one trunk,
then 
\begin{equation}\label{e-sab*}
N_1(\w) ~=~ N_{\rm tree}(0)+N_{\rm tree}(1)+\cdots+
 N_{\rm tree}(\w-1) \,, 
\end{equation}
In fact if $T$ is any tree with just one trunk $T_0\subset T$, then by
collapsing $T_0$ to a point we obtain a tree $T/T_0$ with weight
$$\w(T/T_0)~=~\w(T)~-~\w(T_0)~<~\w(T)\,.$$
 Conversely, $T$ can be reconstructed by pasting
$T/T_0$ onto $T_0$. The identity (\ref{e-sab*}) follows easily.
Note that $N_{\rm tree}(0)=1$, since there is a unique (trivial) tree
 of weight zero.\medskip

\medskip

\begin{table}[!ht]
\begin{center}
\begin{tabular}{|l|cccccc|}
           \hline
           $\w$ & 0 & 1 & 2 & 3 & 4& 5\cr
           \hline
           $N_1(\w)$ & 0 & 1 & 2 & 5 & 13 & 37\cr
           $N_{\rm tree}(\w)$ & 1 & 1 & 3 & 8& 24& 71\cr
           \hline
\end{tabular}
\vspace{.3cm}
\caption{\it Numbers of trees with given total weight $\w$.\label{tB1}}
\end{center}
\end{table}
The values for $\w\le 5$ are shown in Table~\ref{tB1}.
This table can be constructed inductively as follows. Suppose 
that we know the values
 \hbox{$N_{\rm tree}(\w')$} for $\w'<\w$. Then $N_1(\w)$ can be computed
immediately from equation (\ref{e-sab*}). On the other hand, if we know
\hbox{$N_1(\w')$} for all \hbox{$\w'\le \w$}, then \hbox{$N_{\rm tree}(\w)$} can 
be computed as follows.

 Note that any wedge sum
expression (\ref{e-ws}) gives rise to a partition of the total weight $\w$,
that is a sequence  of positive integers which can be ordered so that

$$ \w(T_{i(1)})~\le~ \w(T_{i(2)})~\le~ \cdots~\le~ \w(T_{i(k)}), $$
with sum   equal to $\w$.
First consider the special case where all $k$ of these wedge summands
have the same weight $\w_0=\w/k$. Then there are $N_1(\w_0)$
possible choices for each of these summands, where their order
doesn't matter. The total number of possibilities in this case is equal to the
binomial coefficient 
\begin{equation}\label{e-bc}
\left(\begin{matrix}
N_1(\w_0)+k-1\\ k
\end{matrix}\right)\,.
\end{equation}
 To see this, let  $a_h$ be the number of
copies of the $h^{\rm th}$ tree in this $k$-fold wedge sum,
so that \hbox{$~a_h\ge 0~$} with \hbox{$~a_1+a_2+\cdots +a_{N_1(\w_0)}~=~k\,$.}
Then the partial sums \hbox{$~a_1+a_2+\cdots+a_i+i~$} with \hbox{$~1\le i<
N_1(\w_0)~$} can be any increasing sequence of \hbox{$N_1(\w_0)-1$} distinct
integers between one and \hbox{$~N_1(\w_0)+k-1$}. Hence the number of
 possibilities is given by the binomial coefficient (\ref{e-bc}).\smallskip

More generally,
suppose that for each $n$ between $1$ and $\w$  there are $~k_n\ge 0~$ summands
of weight $n$, so that
$$ \w~=~ k_1+2\,k_2+\cdots+ \w\,k_\w\,.$$
Then the number of possibilities is equal to the product
$$ \prod_{k_n>0} \left(\begin{matrix}
N_1(n)+k_n-1\\ k_n
\end{matrix}\right)\,,$$
taken over all $n$ with $k_n>0\,.$
Thus, all together, the number of possibilities is given by
\begin{equation}\label{e-tree-ct}
 N_{\rm tree}(\w)~=~\sum_{\rm partitions}~\prod_{k_n>0}
\left(\begin{matrix}
N_1(n)+k_n-1\\ k_n
\end{matrix}\right)\,,
\end{equation}
to be summed over all partitions \hbox{$~\sum n\,k_n\,=\,\w$}.\medskip

As an example, suppose that $N_1(\w)$ is known for \hbox{$\w\le 4$}. Using the
values shown in Table \ref{tB1}, since the 
integer 4 has 5 different partitions
$$ 1+1+1+1~~=~~1+1+2~~=~~1+3~~=~~2+2~~=~~4$$
(using a different notation for partitions),
it follows that \hbox{$~N_{\rm tree}(4)~$} 
 can be expressed as a 5-fold sum
$$\quad\qquad\left(\begin{matrix}4 \\ 4\end{matrix}\right)~+~
 \left(\begin{matrix} 2\\2 \end{matrix}\right)
 \left(\begin{matrix}2 \\ 1\end{matrix}\right)~+~
 \left(\begin{matrix}1 \\1 \end{matrix}\right)
 \left(\begin{matrix}5 \\ 1\end{matrix}\right)~+~
 \left(\begin{matrix}3 \\ 2\end{matrix}\right)~+~
 \left(\begin{matrix}13 \\ 1\end{matrix}\right)
\,,$$
yielding $$~N_{\rm tree}(4)~=~1+2+5+3+13~=~24\,.$$ Other entries in
Table \ref{tB1} can be computed similarly.
\bigskip

The following table lists the 
number of connected schemes for each given value of the cyclic weight $\w(C)$
together with the tree weight $\sum_{s\in|C|} \w(T(s))$, within the range
 \hbox{$~\w(S)\,=\,\w(C) + \sum \w(T)\,\le\, 6$}. 

\begin{table}[!ht]
\begin{center}
\begin{tabular}{cc|cccccc|c|}
\hline
\multicolumn{2}{|c|}{$\sum \w(T)=\!\!\!$}& 0 & 1 & 2 & 3 & 4 &5\\ 
\hline
\multicolumn{1}{|c|}{\multirow{6}{*}{$\w(C)=$}} &
\multicolumn{1}{|c|}{1} & 1 & 1 & 3 & 8 & 24 & 71   \\ 
\multicolumn{1}{|c|}{}                        &
\multicolumn{1}{|c|}{2} &  2 & 2& 7 & 19&62&  \\ 
\multicolumn{1}{|c|}{}                       &
\multicolumn{1}{|c|}{3} & 3 & 4 & \underline{14} &45& & \\ 
\multicolumn{1}{|c|}{}                        &
\multicolumn{1}{|c|}{4} & 5 & 8& 31 & & &\\ 
\multicolumn{1}{|c|}{}                        &
\multicolumn{1}{|c|}{5} & 7 &16& & & &\\ 
\multicolumn{1}{|c|}{}                        &
\multicolumn{1}{|c|}{6}& 13& & & & &\\ \cline{1-8}
\end{tabular}
\vspace{.3cm}
\caption{\it Number of connected schemes $S$ with given $\w(C)$ and 
$\sum \w(T)$. 
\label{tB2}}
\end{center}
\end{table}

 Rather than explaining each entry in this table, let me simply give a detailed
explanation for the one typical entry which is underlined in the table,
 corresponding to cyclic weight
\hbox{$\w(C)=3$},~ tree weight \hbox{$\sum \w(T)=2$},~ and hence total weight 
\hbox{$3+2=5$}.
 For this example,  we need to know the numbers
\hbox{$N_{\rm tree}(1)=1$} and \hbox{$N_{\rm tree}(2)=3$}, and we need to study
 each of the three cyclic
schemes of weight \hbox{$\w(C)=3$} separately. Here are the three cases, with
 notation as in {\bf(b)} above.

\begin{itemize}
\item For the cycle $(3)$, with a single vertex of weight 3, we can paste
any one of the three trees of weight 2 onto the unique cyclic vertex,
 so we get a total of 3 possible schemes.\medskip

\item  For the cycle $(1,2)$, we can paste a tree of weight 2
 onto either one of the two vertices, yielding 6 distinct possibilities.
 But we can also paste a  tree of weight one onto each vertex, yielding a
 7-th possibility.\medskip

\item For the cycle $(1,1,1)$, note that there is a cyclic group of
 symmetries. Again we can paste a tree of weight 2 onto any vertex, but
because of the symmetries, it doesn't matter which vertex we choose, so
 there are three distinct possibilities.
Similarly, we can paste a tree of weight one onto each of two vertices,
yielding a 4-th possibility. (Again, because of the symmetries,
it doesn't matter which two we choose.)
\end{itemize}
\smallskip

\noindent Thus, all together, we get $~3+7+4=14~$
 distinct schemes, as listed in Table \ref{tB2}. Other entries in this
table can be computed similarly.\bigskip


The number $N_{\rm c}(\w)$ of connected schemes of weight
$\w$ can be obtained by adding entries along the diagonal 
\hbox{$~\w(C)+\sum \w(T)=\w~$}
in Table~\ref{tB2}. For example
$$  N_{\rm c}(5)~=~ 7+8+14+19+24~=~72\,.$$
For $n\le 6$, the total number $N_{\rm c}(\w)$  of connected
schemes of  weight $\w$, computed in this way,
 is listed in the middle row of Table~\ref{tB3}.
The total number $N(\w)$ of all schemes, connected or not, can then be computed
by a formula completely analogous to (\ref{e-tree-ct}) above,
and is listed in the last row below (as well as in Table \ref{t-numbers}).

\begin{table}[!ht]
\begin{center}
\begin{tabular}{|l|cccccc|}
           \hline
           $\w$  & 1 & 2 & 3 & 4& 5 &6\cr
\hline
           $N_{\rm c}(\w)$ & 1 & 3 & 8 & 24 & 72 & 238\cr
           $N(\w)$ & 1 & 4 & 12 & 42& 138& 494\cr
\hline
\end{tabular}
\vspace{.3cm}
\caption{\it The total count.\label{tB3}}\end{center}
\end{table}

\bibliographystyle{amsplain}

\end{document}